\documentclass[12pt,twoside]{article}
\usepackage[hmargin=0.8in,vmargin=0.9in]{geometry}
\usepackage{fancyhdr}
\usepackage{graphicx}
\usepackage{subfigure}
\usepackage{amssymb}
\usepackage{amsmath}
\usepackage{amsfonts}
\usepackage{theorem}
\usepackage{mathrsfs}
\usepackage{mathtools}
\usepackage{bm}
\usepackage{color}
\usepackage{setspace}
\usepackage{exscale}
\usepackage{relsize}
\usepackage{epstopdf}
\usepackage{float}
\usepackage{cite}

\usepackage{nicefrac}
\usepackage{setspace}

\usepackage{booktabs,multirow} 
\usepackage{array} 
\usepackage{paralist} 
\usepackage{verbatim} 
\usepackage{subfigure} 

\theoremstyle{plain}			
\newtheorem{thm}{Theorem}[section]

\newtheorem{rmk}[thm]{Remark}
{\theorembodyfont{\rmfamily}}

\renewcommand{\thefigure}{\arabic{section}.\arabic{figure}}
\renewcommand{\thetable}{\arabic{section}.\arabic{table}}

\usepackage{tikz}
\usetikzlibrary{positioning}
\usepackage{xcolor}

\allowdisplaybreaks[1]

\numberwithin{equation}{section}
\numberwithin{figure}{section}
\numberwithin{table}{section}

\newcommand\eref[1]{(\ref{#1})}

\newcommand*\xbar[1]{%
  \hbox{%
    \vbox{%
      \hrule height 0.5pt 
      \kern0.4ex
      \hbox{%
        \kern-0.05em
        \ensuremath{#1}%
        \kern-0.00em
      }%
    }%
  }%
}

\newcommand{\bmF}{\bm{\mathcal{F}}}

\newcommand{\bmK}{\bm{\mathcal{K}}}
\newcommand{\bmL}{\bm{\mathcal{L}}}
\newcommand{\mF}{\bm{F}}
\newcommand{\mK}{\bm{K}}

\newcommand{\mG}{\bm{G}}
\newcommand{\mL}{\bm{L}}

\newcommand{\mU}{\bm{U}}
\newcommand{\mE}{\bm{E}}

\newcommand{\mo}{\bm{0}}

\newcommand{\dx}{\Delta x}
\newcommand{\dy}{\Delta y}

\newcommand{\hf}{{\frac{1}{2}}}

\newcommand{\jph}{{j+\frac{1}{2}}}
\newcommand{\jmh}{{j-\frac{1}{2}}}
\newcommand{\kph}{{k+\frac{1}{2}}}
\newcommand{\kmh}{{k-\frac{1}{2}}}

\def\softd{{\leavevmode\setbox1=\hbox{d}%
          \hbox to 1.05\wd1{d\kern-0.4ex{\char039}\hss}}}

\title{Novel Local Characteristic Decomposition Based Path-Conservative Central-Upwind Schemes}
\author{Shaoshuai Chu\thanks{Department of Mathematics and Shenzhen International Center for Mathematics, Southern University of Science and
Technology, Shenzhen, 518055, China; {\tt chuss@mail.sustech.edu.cn}}, Michael Herty\thanks{Department of Mathematics, RWTH Aachen
University, 52056 Aachen, Germany; {\tt herty@igpm.rwth-aachen.de}}, and Alexander Kurganov\thanks{Department of Mathematics, Shenzhen
International Center for Mathematics, and Guangdong Provincial Key Laboratory of Computational Science and Material Design, Southern
University of Science and Technology, Shenzhen, 518055, China; {\tt alexander@sustech.edu.cn}}}

\begin{document}

\date{}
\maketitle
\begin{abstract}
We introduce local characteristic decomposition based path-conservative central-upwind schemes for (nonconservative) hyperbolic systems of
balance laws. The proposed schemes are made to be well-balanced via a flux globalization approach, in which source terms are incorporated
into the fluxes: This helps to enforce the well-balanced property when the resulting quasi-conservative system is solved using the local
characteristic decomposition based central-upwind scheme recently introduced in [{\sc A. Chertock, S. Chu, M. Herty, A. Kurganov, and M.
Luk\'{a}\v{c}ov\'{a}-Medvi{\softd}ov\'{a}}, J. Comput. Phys., 473 (2023), Paper No. 111718]. Nonconservative product terms are also
incorporated into the global fluxes using a path-conservative technique. We illustrate the performance of the developed schemes by applying
them to one- and two-dimensional compressible multifluid systems and thermal rotating shallow water equations.
\end{abstract}



\section{Introduction}
This paper is focused on the development of highly accurate and robust numerical methods for the nonconservative hyperbolic systems of
balance laws, which in the one-dimensional (1-D) case read as
\begin{equation}
\mU_t+\mF(\mU)_x=B(\mU)\mU_x+\bm S(\mU),
\label{1.1}
\end{equation}
where $x$ is a spatial variable, $t$ is the time, $\mU\in\mathbb R^d$ is a vector of unknown functions, $\mF:\mathbb R^d\to\mathbb R^d$ are
nonlinear fluxes, $ B(\mU)\in\mathbb R^{d\times d}$, and $\bm S(\mU):\mathbb R^d \to \mathbb R^d$ are source terms.

Designing a numerical scheme for \eref{1.1} is a challenging task for a number of reasons. First, it is well-known that solutions of
\eref{1.1} may develop complicated nonsmooth structures even for smooth initial data. Second, the presence of nonconservative products (when
$B\ne0$) brings many challenges for studying these systems both theoretically and numerically. In general, when $\mU$ is discontinuous, weak
solutions of \eref{1.1} cannot be understood in the sense of distributions. Instead, one can introduce Borel measure solutions
\cite{maso1995,lefloch02,lefloch2012}. This concept of weak solutions was utilized to develop path-conservative finite-volume (FV) methods;
see, e.g., \cite{pares09,CMP2017,CP20,CKM,PGCCMP,SGBNP,BDGI,Cha} and reference therein. Third, a novel scheme should be well-balanced (WB)
in the sense that it should be capable of exactly preserving some of the physically relevant steady-state solutions satisfying
\begin{equation*}
\mF(\mU)_x-B(\mU)\mU_x-\bm S(\mU)=M(\mU)\bm E(\mU)_x=\mo.
\end{equation*}
Here, $M\in\mathbb R^{d\times d}$ and $\bm E$ is the vector of equilibrium variables, which are constant at the steady states provided $M$
is invertible; see \cite{KLX_21}.

In this paper, we focus on one of the popular classes of FV methods, the Riemann-problem-solver-free central-upwind (CU) schemes, which were
originally introduced in \cite{Kurganov01,Kurganov00,Kurganov07} as a ``black-box'' solver for general multidimensional hyperbolic systems
of conservation laws and then successfully applied to a variety of hyperbolic systems. In \cite{CKM}, the CU schemes were first reformulated
in a path-conservative framework and then extended to the 1-D nonconservative systems \eref{1.1}, resulting in the second-order
path-conservative CU (PCCU) schemes. These schemes are WB in the sense that they are capable of preserving simple steady-state solutions
like ``lake-at-rest" of several shallow water models.

In order to design PCCU schemes capable of preserving more complicated steady states, path-conservative techniques were incorporated into
the flux globalization framework: This results in a new class of flux globalization based WB PCCU schemes introduced in \cite{KLX_21}. A
flux globalization approach, which was introduced in \cite{CHO_18,CDH2009,DDMA11,GC2001,MGD11}, relies on the following quasi-conservative
form of \eref{1.1}:
\begin{equation}
\bm U_t+\bm K(\bm U)_x=\bm0,
\label{1.4}
\end{equation}
where $\bm K(\bm U)$ is a global flux
\begin{equation}
\bm K(\bm U)=\bm F(\bm U)-\bm W(\bm U),\quad\bm W(\bm U)=\int\limits^x_{\widehat x}B(\bm U){\bm U}_\xi(\xi,t)+S(\mU(\xi,t))\,{\rm d}\xi,
\label{1.5}
\end{equation}
and ${\widehat x}$ is an arbitrary number. The flux globalization based WB PCCU schemes are obtained by applying the CU numerical flux to
the semi-discretization of \eref{1.4} and by using the path-conservative technique for the evaluation of the integrals in \eref{1.5}. These
schemes have been applied to a variety of hyperbolic systems in \cite{CKL23,CKLX_22,CKLZ_23,CKN_23,KLX_21,CKX_23a}.

Even though the flux globalization based WB PCCU schemes from \cite{KLX_21} are accurate, efficient, and robust, higher resolution of the
numerical solutions can be achieved by reducing numerical dissipation present in those schemes. This can be done, for example, by using the
local characteristic decomposition (LCD) based CU (LCD-CU) numerical fluxes instead of the original CU ones employed in
\cite{CKL23,CKLX_22,CKLZ_23,CKN_23,KLX_21}. The LCD-CU schemes have been recently introduced in \cite{CCHKL_22} for hyperbolic systems of
conservation laws and we use the same LCD based framework to enhance the performance of the flux globalization based WB PCCU schemes.
The proposed scheme combines for the first time the modifications enhancing the performance of the CU schemes for both conservative and
nonconservative hyperbolic systems in multiple space dimensions. We derive the new LCD-PCCU schemes and apply them to two particular
systems. First, we consider the $\gamma$-based compressible multifluid systems, which can be, in principle, written in the conservative
form, but we use its nonconservative formulation to achieve high resolution of material interfaces. Second, we consider a challenging
thermal rotating shallow water (TRSW) model, which contains geometric source term due to a nonflat bottom topography, source terms due to
Coriolis forces, and generically nonconservative momentum exchange terms.

The rest of the paper is organized as follows. In \S\ref{sec2}, we briefly describe the 1-D second-order LCD-CU scheme. In \S\ref{sec3}, we
propose the novel scheme by considering the nonconservative system \eref{1.1} and design a 1-D flux globalization based WB LCD-PCCU scheme.
In \S\ref{sec4}, we extend the 1-D LCD-PCCU scheme to the two-dimensional (2-D) case. Finally, in \S\ref{sec5}, we apply the developed
LCD-PCCU schemes to the $\gamma$-based compressible multifluid systems and TRSW equations. We illustrate the performance of the scheme on a
number of numerical examples.

\section{1-D LCD-CU Scheme}\label{sec2}
In this section, we briefly describe the second-order semi-discrete LCD-CU scheme, which is a slightly modified version of the scheme from
\cite{CCHKL_22}. For the difference between the schemes, we refer the reader to Remarks \ref{rmk2.1}--\ref{rmk2.3} below.

We consider a 1-D hyperbolic system of conservative laws
\begin{equation}
\mU_t+\mF(\mU)_x=\mo,
\label{2.1}
\end{equation}
cover the computational domain with the uniform cells $C_j:=[x_\jmh,x_\jph]$ of size $\dx$ centered at $x_j=(x_\jmh+x_\jph)/2$, and denote
by $\xbar\mU_j(t)$ the computed cell averages of $\mU(\cdot,t)$, 
\begin{equation*}
\xbar\mU_j(t):\approx\frac{1}{\dx}\int\limits_{C_j}\mU(x,t)\,{\rm d}x,\qquad j=1,\ldots,N,
\end{equation*}
which are assumed to be available at a certain time $t\ge0$. Note that both $\xbar\mU_j(t)$ and many other indexed quantities below are
time-dependent, but from here on, we suppress their time-dependence for the sake of brevity.

The computed cell averages \eref{2.1} are evolved in time by numerically solving the following system of ODEs:
\begin{equation*}
\frac{{\rm d}\xbar\mU_j}{{\rm d}t}=-\frac{\bmF^{\rm LCD}_\jph-\bmF^{\rm LCD}_\jmh}{\dx},
\end{equation*}
where the LCD-CU numerical fluxes are given by 
\begin{equation}
\bmF^{\rm LCD}_\jph=R_\jph P^{\rm LCD}_\jph R^{-1}_\jph\mF^-_\jph+R_\jph M^{\rm LCD}_\jph R^{-1}_\jph\mF^+_\jph+
R_\jph Q^{\rm LCD}_\jph R^{-1}_\jph\big(\mU^+_\jph-\mU^-_\jph\big).
\label{2.3}
\end{equation}
Here, $\mF^\pm_j=\mF(\mU^\pm_j)$ and $\mU^\pm_\jph$ are the one-sided values of a piecewise linear interpolant
\begin{equation}
\widetilde\mU(x)=\,\xbar\mU_j+(\mU_x)_j(x-x_j),\quad x\in C_j,
\label{2.3a}
\end{equation}
at the cell interface $x=x_\jph$, namely
\begin{equation}
\mU^-_\jph=\,\xbar\mU_j+\frac{\dx}{2}(\mU_x)_j,\quad\mU^+_\jph=\,\xbar\mU_{j+1}-\frac{\dx}{2}(\mU_x)_{j+1}.
\label{2.5a}
\end{equation}
In order to ensure the non-oscillatory nature of this reconstruction, we compute the slopes $(\mU_x)_j$ in \eref{2.5a} using a generalized
minmod limiter \cite{lie03,Nessyahu90,Sweby84}:
\begin{equation}
(\mU_x)_j={\rm minmod}\left(\theta\frac{\,\xbar\mU_j-\,\xbar\mU_{j-1}}{\dx},\,\frac{\,\xbar\mU_{j+1}-\,\xbar\mU_{j-1}}{2\dx},\,
\theta\frac{\,\xbar\mU_{j+1}-\,\xbar\mU_j}{\dx}\right),\quad\theta\in[1,2],
\label{2.5b}
\end{equation}
applied in the component-wise manner. Here, the minmod function is defined as
\begin{equation}
{\rm minmod}(z_1,z_2,\ldots):=\begin{cases}
\min_j\{z_j\}&\mbox{if}~z_j>0\quad\forall\,j,\\
\max_j\{z_j\}&\mbox{if}~z_j<0\quad\forall\,j,\\
0            &\mbox{otherwise},
\end{cases}
\label{2.7b}
\end{equation}
and the parameter $\theta$ can be used to control the amount of numerical viscosity present in the resulting scheme. We have taken
$\theta=1.3$ in all of the numerical experiments reported in \S\ref{sec5}.

The diagonal matrices $P^{\rm LCD}_\jph$, $M^{\rm LCD}_\jph$, and $Q^{\rm LCD}_\jph$ are given by
\begin{equation}
\begin{aligned}
&P^{\rm LCD}_\jph={\rm diag}\big((P^{\rm LCD}_1)_\jph,\ldots,(P^{\rm LCD}_d)_\jph\big),\quad
M^{\rm LCD}_\jph={\rm diag}\big((M^{\rm LCD}_1)_\jph,\ldots,(M^{\rm LCD}_d)_\jph\big),\\
&Q^{\rm LCD}_\jph={\rm diag}\big((Q^{\rm LCD}_1)_\jph,\ldots,(Q^{\rm LCD}_d)_\jph\big)
\end{aligned}
\label{2.7a}
\end{equation}
with
\begin{align}
&\hspace*{-0.2cm}\big((P^{\rm LCD}_i)_\jph,(M^{\rm LCD}_i)_\jph,(Q^{\rm LCD}_i)_\jph\big)\label{2.7}\\
&=\left\{\begin{aligned}
&\frac{1}{(\lambda^+_i)_\jph-(\lambda^-_i)_\jph}\Big((\lambda^+_i)_\jph,-(\lambda^-_i)_\jph,(\lambda^+_i)_\jph(\lambda^-_i)_\jph\Big)&&
\mbox{if}~(\lambda^+_i)_\jph-(\lambda^-_i)_\jph> \varepsilon_0,\\[0.3ex]
&\frac{1}{a^+_\jph-a^-_\jph}\Big(a^+_\jph, -a^-_\jph, a^+_\jph a^-_\jph\Big)&&
\hspace*{-3.5cm}\mbox{if}~(\lambda^+_i)_\jph-(\lambda^-_i)_\jph\le\varepsilon_0~~\mbox{and}~~a^+_\jph-a^-_\jph>\varepsilon_0,\\[0.3ex]
&\Big(\hf,\hf,0\Big)&&\hspace*{-3.5cm}\mbox{if}~a^+_\jph-a^-_\jph\le\varepsilon_0,
\end{aligned}\right.\nonumber
\end{align}
where
\begin{equation*}
\begin{aligned}
&(\lambda^+_i)_\jph=\max\left\{\lambda_i\big(A(\mU^-_\jph)\big),\,\lambda_i\big(A(\mU^+_\jph)\big),\,0\right\},\\
&(\lambda^-_i)_\jph=\min\left\{\lambda_i\big(A(\mU^-_\jph)\big),\,\lambda_i\big(A(\mU^+_\jph)\big),\,0\right\},
\end{aligned}\quad i=1,\ldots,d.
\end{equation*}
Here, $A(\mU)=\frac{\partial\mF}{\partial\mU}(\mU)$ is the Jacobian and the matrices $R_\jph$ and $R^{-1}_\jph$ are such that
$R^{-1}_\jph\widehat A_\jph R_\jph$ is diagonal, where $\widehat A_\jph=A(\widehat\mU_\jph)$ and $\widehat\mU_\jph$ is either a simple
average $(\xbar\mU_j+\xbar\mU_{j+1})/2$ or another type of average of the $\xbar\mU_j$ and $\xbar\mU_{j+1}$ states. The one-sided local
speeds of propagation $a^\pm_\jph$ can be estimated using the largest and the smallest eigenvalues of $A$, for example, by taking 
\begin{equation*}
\begin{aligned}
 a^+_\jph=\max\left\{\lambda_d\big(A(\mU^+_\jph)\big),\lambda_d\big(A(\mU^-_\jph)\big),0\right\},\quad
 a^-_\jph=\min\left\{\lambda_1\big(A(\mU^+_\jph)\big),\lambda_1\big(A(\mU^-_\jph)\big),0\right\}.
\end{aligned}
\end{equation*}
Finally, $\varepsilon_0$ in \eref{2.7} is a very small desingularization constant, which is taken $\varepsilon_0=10^{-18}$ in all of the
numerical examples reported in \S\ref{sec5}.
\begin{rmk}\label{rmk2.1}
The numerical flux \eref{2.3} is identical to the LCD-CU numerical flux in \cite{CCHKL_22}, which we have simplified using the fact that
\begin{equation}
P^{\rm LCD}_\jph+M^{\rm LCD}_\jph=I.
\label{2.9}
\end{equation}
\end{rmk}
\begin{rmk}\label{rmk2.2}
We perform a piecewise linear reconstruction for the conservative variables $\mU$ rather than the local characteristic variables as it was
done in \cite{CCHKL_22}. As no substantial spurious oscillations have been observed in any of the conducted numerical examples, there is no
need to switch to the local characteristic variables in this case. This allows the developed scheme to be more efficient than the original
LCD-CU scheme from \cite{CCHKL_22}.
\end{rmk}
\begin{rmk}\label{rmk2.3}
The desingularization in \eref{2.7} is performed differently than in \cite{CCHKL_22}. Our numerical experiments indicate that the current
desingularization is more robust compared with the previous desingularization strategy.
\end{rmk}

\section{1-D Flux Globalization Based LCD-PCCU Scheme}\label{sec3}
In this section, we propose an extension of the LCD-CU scheme from \S \ref{sec2} to the nonconservative system \eref{1.1} and design a 1-D
flux globalization based LCD-PCCU scheme.

We apply the LCD-CU scheme to the quasi-conservative system \eref{1.4}--\eref{1.5} to obtain
\begin{equation}
\frac{{\rm d}\xbar\mU_j}{{\rm d}t}=-\frac{\bmK^{\rm LCD}_\jph-\bmK^{\rm LCD}_\jmh}{\dx},
\label{3.1}
\end{equation}
where the numerical fluxes $\bmK^{\rm LCD}_\jph$ are given by (compare with \eref{2.3})
\begin{equation}
\bmK^{\rm LCD}_\jph=R_\jph P^{\rm LCD}_\jph R^{-1}_\jph\mK^-_\jph+R_\jph M^{\rm LCD}_\jph R^{-1}_\jph\mK^+_\jph+
R_\jph Q^{\rm LCD}_\jph R^-_\jph\big(\breve\mU^+_\jph-\breve\mU^-_\jph\big).
\label{3.1a}
\end{equation}
Here, $\mK^\pm_\jph$ are the global fluxes, which are computed as the same as in \cite[\S2]{CKLX_22}. 

The global fluxes $\mK^\pm_\jph$ are computed using the point values $\mU^\pm_\jph$, which, in order to ensure the WB property of the
resulting scheme, are to be obtained with the help of a piecewise linear reconstruction \eref{2.3a}--\eref{2.7b} applied to the equilibrium
variables $\mE$ (instead of $\mU$):
\begin{equation*}
\bm E^-_\jph=\,\xbar{\bm E}_j+\frac{\dx}{2}(\bm E_x)_j,\quad\bm E^+_\jph=\,\xbar{\bm E}_{j+1}-\frac{\dx}{2}(\bm E_x)_{j+1},
\end{equation*}
where $\bm E_j:=\bm E(\,\xbar{\bm U}_j)$. We then obtain the corresponding values $\bm U^\pm_\jph$ by solving the (nonlinear systems of)
equations
\begin{equation}
\bm E(\bm U^+_\jph)=\bm E^+_\jph\quad\mbox{and}\quad\bm E(\bm U^-_\jph)=\bm E^-_\jph
\label{3.3}
\end{equation}
for $\bm U^+_\jph$ and $\bm U^-_\jph$, respectively. We emphasize that the values $\mU^\pm_\jph$ cannot be used in \eref{3.1a} instead of
$\breve{\bm U}^\pm_\jph$ as $\mU^+_\jph$ may not be equal to $\mU^-_\jph$ at discrete steady states. We therefore use another set of point
values $\breve\mU^\pm_\jph$, which are obtained by solving modified versions of \eref{3.3}:
\begin{equation*}
\breve{\bm E}(\breve{\bm U}^+_\jph)=\bm E^+_\jph\quad\mbox{and}\quad\breve{\bm E}(\breve{\bm U}^-_\jph)=\bm E^-_\jph.
\end{equation*}
Here, the functions $\breve{\bm E}$ are the modifications of the functions $\bm E$, which are made in such a way that
$\breve{\bm U}^+_\jph=\breve{\bm U}^-_\jph$ as long as $\bm E^+_\jph=\bm E^-_\jph$ and, at the same time,
$\breve{\bm U}^+_\jph=\bm U^+_\jph+{\cal O}((\dx)^2)$ and $\breve{\bm U}^-_\jph=\bm U^-_\jph+{\cal O}((\dx)^2)$ for smooth solutions; see
\cite{CKN_23,CKLX_22,CKL23,CKLZ_23,KLX_21} for examples of such modifications.

The diagonal matrices $P^{\rm LCD}_\jph$, $M^{\rm LCD}_\jph$, and $Q^{\rm LCD}_\jph$ in \eref{3.1a} are defined by \eref{2.7a}--\eref{2.7},
but with 
\begin{equation*}
\begin{aligned}
(\lambda^+_i)_\jph&=\max\left\{\lambda_i\big({\cal A}(\mU^-_\jph)\big),\,\lambda_i\big({\cal A}(\mU^+_\jph)\big),\,0\right\},\\
(\lambda^-_i)_\jph&=\min\left\{\lambda_i\big({\cal A}(\mU^-_\jph)\big),\,\lambda_i\big({\cal A}(\mU^+_\jph)\big),\,0\right\}, 
\end{aligned}\quad i=1,\ldots,d,
\end{equation*}
where ${\cal A}(\mU)=\frac{\partial\mF}{\partial\mU}(\mU)-B(\mU)$ and the matrices $R_\jph$ and $R^{-1}_\jph$ are such that
$R^{-1}_\jph\widehat{\cal A}_\jph R_\jph$ is diagonal for $\widehat{\cal A}_\jph={\cal A}(\widehat\mU_\jph)$. Finally, the one-sided local
speeds of propagation $a^\pm_\jph$ are estimated by
\begin{equation*}
\begin{aligned}
a^+_\jph=\max\left\{\lambda_d\big({\cal A}(\mU^+_\jph)\big),\lambda_d\big({\cal A}(\mU^-_\jph)\big),\,0\right\},\quad
a^-_\jph=\min\left\{\lambda_1\big({\cal A}(\mU^+_\jph)\big),\lambda_1\big({\cal A}(\mU^-_\jph)\big),\,0\right\}.
\end{aligned}
\end{equation*}
\begin{rmk}
The introduced 1-D semi-discrete flux globalization based LCD-PCCU scheme is WB. In order to prove this, we only need to show that if the
discrete data are at a steady state, namely, if $\xbar E_j\equiv\widehat E=\texttt{Const}$, then the RHS of \eref{3.1} vanishes. Following
the proof of \cite[Theorem 2.1]{CKLX_22}, we obtain
\begin{equation}
\mK^+_\jph=\mK^-_\jph=\mK^+_\jmh=\widehat\mK=\texttt{Const},\quad\forall j.
\label{3.4}
\end{equation}
Substituting \eref{3.4} and \eref{2.9} into \eref{3.1a} and taking into account the fact that at the steady states
$\breve\mU^+_\jph=\breve\mU^-_\jph,~\forall j$, yields 
\begin{equation*}
\bmK^{\rm LCD}_\jph=\widehat\mK+R_\jph Q^{\rm LCD}_\jph R^-_\jph\big(\breve\mU^+_\jph-\breve\mU^-_\jph\big)=\widehat\mK,\quad\forall j,
\end{equation*}
which ensures that the RHS of \eref{3.1} vanishes.
\end{rmk}

\section{2-D Flux Globalization Based LCD-PCCU Scheme}\label{sec4}
In this section, we propose a novel 2-D LCD-PCCU scheme for
\begin{equation}
\mU_t+\mF(\mU)_x+\mG(\mU)_y=B(\mU)\mU_x+C(\mU)\mU_y+\bm S^x(\mU)+\bm S^y(\mU),
\label{1.2}
\end{equation}
where $\mG:\mathbb R^d\to\mathbb R^d$, $ C(\mU)\in\mathbb R^{d\times d}$, $\bm S^x(\mU):\mathbb R^d\to\mathbb R^d$, and
$\bm S^y(\mU):\mathbb R^d\to\mathbb R^d$. To this end, we design the 2-D scheme in a ``dimension-by-dimension'' manner.

As in the 1-D case, we first rewrite the general 2-D nonconservative system \eref{1.2} in the following quasi-conservative form:
\begin{equation}
\bm U_t+\bm K(\bm U)_x+\bm L(\bm U)_y=\mo,\quad\bm K(\bm U)=\bm F(\bm U)-\bm W^x(\bm U),\quad\bm L(\bm U)=\bm G(\bm U)-\bm W^y(\bm U),
\label{5.1}
\end{equation}
where
\begin{equation*}
\bm W^x(\bm U):=\int\limits^x_{\widehat x}\left[B(\bm U)\bm U_\xi(\xi,t)+\bm S^x(\mU(\xi,t))\right]{\rm d}\xi,\quad
\bm W^y(\bm U):=\int\limits^y_{\widehat y}\left[C(\bm U)\bm U_\eta(\eta,t)+\bm S^y(\mU(\eta,t))\right]{\rm d}\eta,
\end{equation*}
with $\widehat x$ and $\widehat y$ being arbitrary numbers.

Applying the 2-D version of the semi-discrete LCD-CU scheme from \S\ref{sec2} to the quasi-conservative system \eref{5.1} results in the
following system of ODEs:
\begin{equation}
\frac{{\rm d}\xbar\mU_{j,k}}{{\rm d}t}=-\frac{\bmK^{\rm LCD}_{\jph,k}-\bmK^{\rm LCD}_{\jmh,k}}{\dx}-
\frac{\bmL^{\rm LCD}_{j,\kph}-\bmL^{\rm LCD}_{j,\kmh}}{\dy},
\label{5.2}
\end{equation}
where the numerical fluxes $\bmK^{\rm LCD}_{\jph,k}$ and $\bmL^{\rm LCD}_{j,\kph}$ are given by
\begin{equation}
\begin{aligned}
\bmK^{\rm LCD}_{\jph,k}&=R_{\jph,k}P^{\rm LCD}_{\jph,k}R^{-1}_{\jph,k}\mK^-_{\jph,k}+
R_{\jph,k}M^{\rm LCD}_{\jph,k}R^{-1}_{\jph,k}\mK^+_{\jph,k}\\
&+R_{\jph,k}Q^{\rm LCD}_{\jph,k}R^{-1}_{\jph,k}\big(\breve\mU^+_{\jph,k}-\breve\mU^-_{\jph,k}\big),\\[0.8ex]
\bmL^{\rm LCD}_{j,\kph}&=R_{j,\kph} P^{\rm LCD}_{j,\kph} R^{-1}_{j,\kph} \mL^-_{j,\kph}+
R_{j,\kph}M^{\rm LCD}_{j,\kph}R^{-1}_{j,\kph}\mL^+_{j,\kph}\\
&+R_{j,\kph}Q^{\rm LCD}_{j,\kph}R^{-1}_{j,\kph}\big(\breve\mU^+_{j,\kph}-\breve\mU^-_{j,\kph}\big).
\end{aligned}
\label{5.2a}
\end{equation}
All of the terms in \eref{5.2a} are computed as in \S\ref{sec3} separately in the $x$- and $y$-directions. 

The global fluxes $\mK^\pm_{\jph,k}$ and $\mL^\pm_{j,\kph}$ are computed using the WB path-conservative technique as in \cite{CKL23}. The
LCD matrices $P^{\rm LCD}_{\jph,k}$, $M^{\rm LCD}_{\jph,k}$, $Q^{\rm LCD}_{\jph,k}$, $P^{\rm LCD}_{j,\kph}$, $M^{\rm LCD}_{j,\kph}$, and
$Q^{\rm LCD}_{j,\kph}$ are computed as in \cite{CCHKL_22}, but based on the eigensystems of
$\widehat{\cal A}^{\,\mF}_{\jph,k}={\cal A}^{\mF}(\widehat\mU_{\jph,k}):=
\frac{\partial\mF}{\partial\mU}(\widehat\mU_{\jph,k})-B(\widehat\mU_{\jph,k})$ and
$\widehat{\cal A}^{\,\mG}_{j,\kph}={\cal A}^{\mG}(\widehat\mU_{j,\kph}):=
\frac{\partial\mG}{\partial\mU}(\widehat\mU_{j,\kph})-C(\widehat\mU_{j,\kph})$. As in \S\ref{sec2}, we do not employ local characteristic
variables in the piecewise linear reconstruction and perform the desingularization of $P^{\rm LCD}_{\jph,k}$, $M^{\rm LCD}_{\jph,k}$,
$Q^{\rm LCD}_{\jph,k}$, $P^{\rm LCD}_{j,\kph}$, $M^{\rm LCD}_{j,\kph}$, and $Q^{\rm LCD}_{j,\kph}$ as in \eref{2.7}.
\begin{rmk}
The introduced 2-D semi-discrete flux globalization based LCD-PCCU scheme is WB in the sense that it is capable of exactly preserving
discrete steady states, for which $\mK(\mU)$ is independent of $x$ and $\mL(\mU)$ is independent of $y$. We omit the proof for the sake of
brevity.
\end{rmk}

\section{Numerical Examples}\label{sec5}
In this section, we apply the developed LCD-PCCU schemes to the $\gamma$-based compressible multifluid and TRSW equations. We conduct
several numerical experiments and compare the performance of the LCD-PCCU schemes with their PCCU counterparts obtained by applying the
original CU numerical fluxes in \eref{3.1} and \eref{5.2}. For instance, in the 1-D case, such a PCCU scheme is obtained by replacing
\eref{2.7} by
$$
\big((P_i)_\jph,(M_i)_\jph,(Q_i)_\jph\big)
=\left\{\begin{aligned}
&\frac{1}{a^+_\jph-a^-_\jph}\Big(a^+_\jph,-a^-_\jph,a^+_\jph a^-_\jph\Big)&&\mbox{if}~a^+_\jph-a^-_\jph>\varepsilon_0,\\[0.3ex]
&\Big(\hf,\hf,0\Big)&&\mbox{otherwise}.
\end{aligned}\right.
$$

In all of the numerical examples, we have solved the ODE systems \eref{3.1} and \eref{5.2} using the three-stage third-order strong
stability preserving (SSP) Runge-Kutta method (see, e.g., \cite{Gottlieb12,Gottlieb11}) and used the CFL number 0.45. 

\subsection{$\gamma$-Based Compressible Multifluid Equations}
In this section, we show the applications of the developed LCD-PCCU schemes introduced in \S\ref{sec3} and \S\ref{sec4} to the 1-D and 2-D
$\gamma$-based compressible multifluid systems \cite{CZL20,Shyue98}, which in the 2-D case can be written as \eref{1.2} with
$\mU=\big(\rho,\rho u,\rho v,E,\Gamma,\Pi\big)^\top$, $\mF(\mU)=\big(\rho u,\rho u^2+p,\rho uv,u(E+p),0,0\big)^\top$,
$\mG(\mU)=\big(\rho v,\rho uv,\rho v^2+p,v(E+p),0,0)^\top$, $\bm S^x(\mU)=\bm S^y(\mU)=\mo$, $B(\mU)={\rm diag}(0,0,0,0,u,u)$, and
$C(\mU)={\rm diag}(0,0,0,0,v,v)$. Here, $\rho$ is the density, $u$ and $v$ are the $x$- and $y$-velocities, $p$ is the pressure, and $E$ is
the total energy. The system is closed through the following (stiff) equation of state (EOS) for each of the fluid components:
\begin{equation}
p=(\gamma-1)\left[E-\frac{\rho}{2}(u^2+v^2)\right]-\gamma\pi_\infty,
\label{A2}
\end{equation}
where the parameters $\gamma$ and $\pi_\infty$ represent the specific heat ratio and stiffness parameter, respectively (when
$\pi_\infty\equiv0$ for all of the components, the system reduces to the ideal gas multicomponent case). The variables $\Gamma$ and $\Pi$
are defined as $\Gamma:=\frac{1}{\gamma-1}$ and $\Pi:=\frac{\gamma\pi_\infty}{\gamma-1}$.
\begin{rmk}
We note that the studied $\gamma$-based compressible multifluid equations contain no source terms and thus there is no need to numerically
enforce the WB property. We therefore do not need to introduce any equilibrium variables---instead, we reconstruct the conservative
quantities $\mU$ and replace $\breve\mU^\pm_\jph$, $\breve\mU^\pm_{\jph,k}$, and $\breve\mU^\pm_{j,\kph}$ in \eref{3.1a} and \eref{5.2a} by
$\mU^\pm_\jph$, $\mU^\pm_{\jph,k}$, and $\mU^\pm_{j,\kph}$, respectively.
\end{rmk}
\begin{rmk}
In Appendix \ref{appa}, we provide the expressions for the matrices $\widehat{\cal A}_\jph$ and $\widehat{\cal A}^{\,\mF}_{\jph,k}$ and the
corresponding matrices $R_\jph$ and $R_{\jph,k}$ for the 1-D and 2-D $\gamma$-based compressible multifluid systems, respectively. The
matrices $\widehat{\cal A}^{\,\mG}_{j,\kph}$ and $R_{j,\kph}$ can be obtained in an analogous way.
\end{rmk}

\subsubsection{One-Dimensional Examples}
In this section, we consider two numerical examples for the 1-D version of the $\gamma$-based compressible multifluid equations.

\subsubsection*{Example 1---``Shock-Bubble'' Interaction}
In the first example taken from \cite{CKX_23}, we consider the ``shock-bubble'' interaction problem. The initial data are given by
\begin{equation*}
(\rho,u,p;\gamma,\pi_\infty)(x,0)=\begin{cases}
(13.1538,0;1,5/3,0),&|x|<0.25,\\
(1.3333,-0.3535,1.5;1.4,0),&x>0.75,\\
(1,0,1;1.4,0),&\mbox{otherwise,}
\end{cases}
\end{equation*}
which correspond to a left-moving shock, initially located at $x=0.75$, and a resting ``bubble'' with a radius of 0.25, initially located at
the origin. 

\begin{figure}[ht!]
\centerline{\includegraphics[trim=1.2cm 0.4cm 0.9cm 0.4cm, clip, width=5cm]{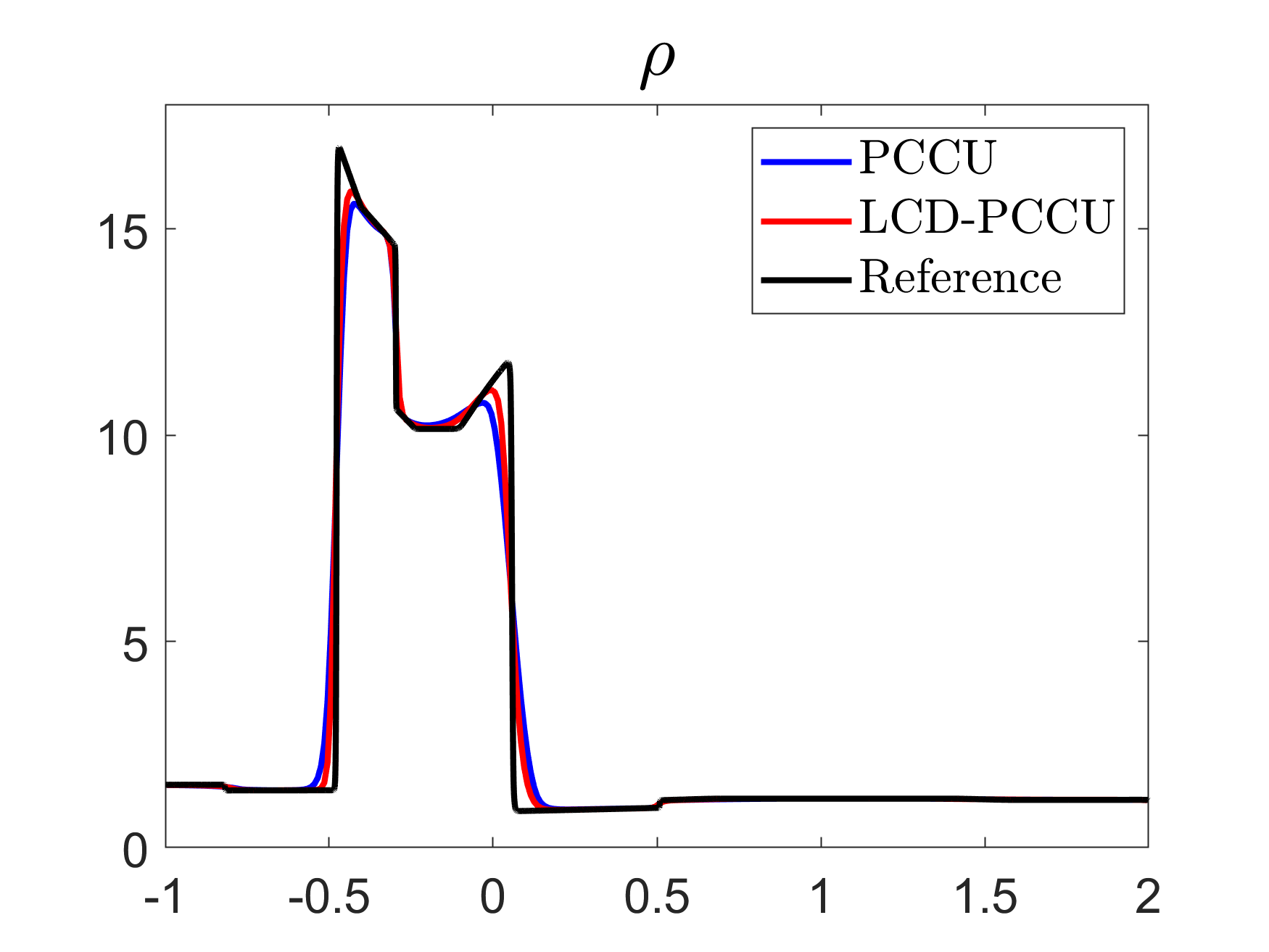}\hspace{1.cm}
            \includegraphics[trim=1.2cm 0.4cm 0.9cm 0.4cm, clip, width=5cm]{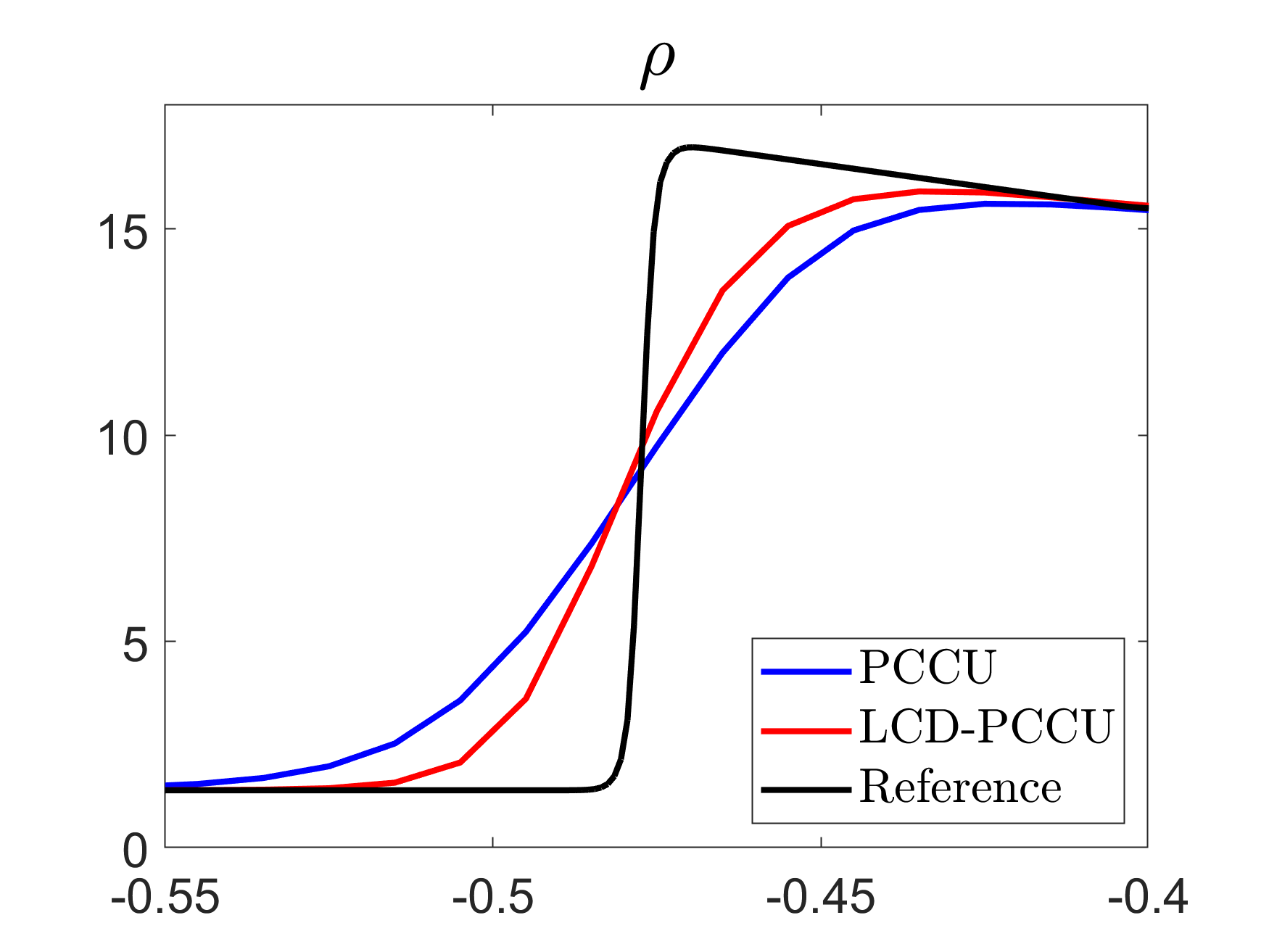}}
\caption{\sf Example 1: Density $\rho$ computed by the PCCU and LCD-PCCU schemes (left) and zoom at $x\in[-0.55,-0.4]$ (right).\label{fig1}}
\end{figure}

We impose the solid wall boundary conditions at the both ends of the computational domain $[-1,2]$ and compute the numerical solutions by
both the LCD-PCCU and PCCU schemes until the final time $t=3$ on a uniform mesh with $\dx=1/100$. The obtained densities are presented in
Figure \ref{fig1} together with the reference solution computed by the PCCU scheme on a much finer mesh with $\dx=1/1000$. As one can see,
the LCD-PCCU scheme achieves sharper resolution than the PCCU one.

\subsubsection*{Example 2---Water-Air ``Shock-Bubble'' Interaction}
In the second 1-D example also taken from \cite{CKX_23}, we consider a gas-liquid multifluid system, where the liquid component is modeled
by the stiff EOS \eref{A2} with $\pi_\infty\gg1$. The initial conditions,
\begin{equation*}
(\rho,u,p;\gamma,\pi_\infty)=\begin{cases}
(0.05,0,1;1.4,0),&|x-6|<3,\\
(1.325,-68.525,19153;4.4,6000),&x>11.4,\\
(1,0,1;4.4,6000),&\mbox{otherwise},
\end{cases}
\end{equation*}
correspond to the left-moving shock, initially located at $x=11.4$, and a resting air ``bubble'' with a radius 3, initially located at
$x=6$. 

\begin{figure}[ht!]
\centerline{\includegraphics[trim=0.9cm 0.4cm 0.7cm 0.4cm, clip, width=5cm]{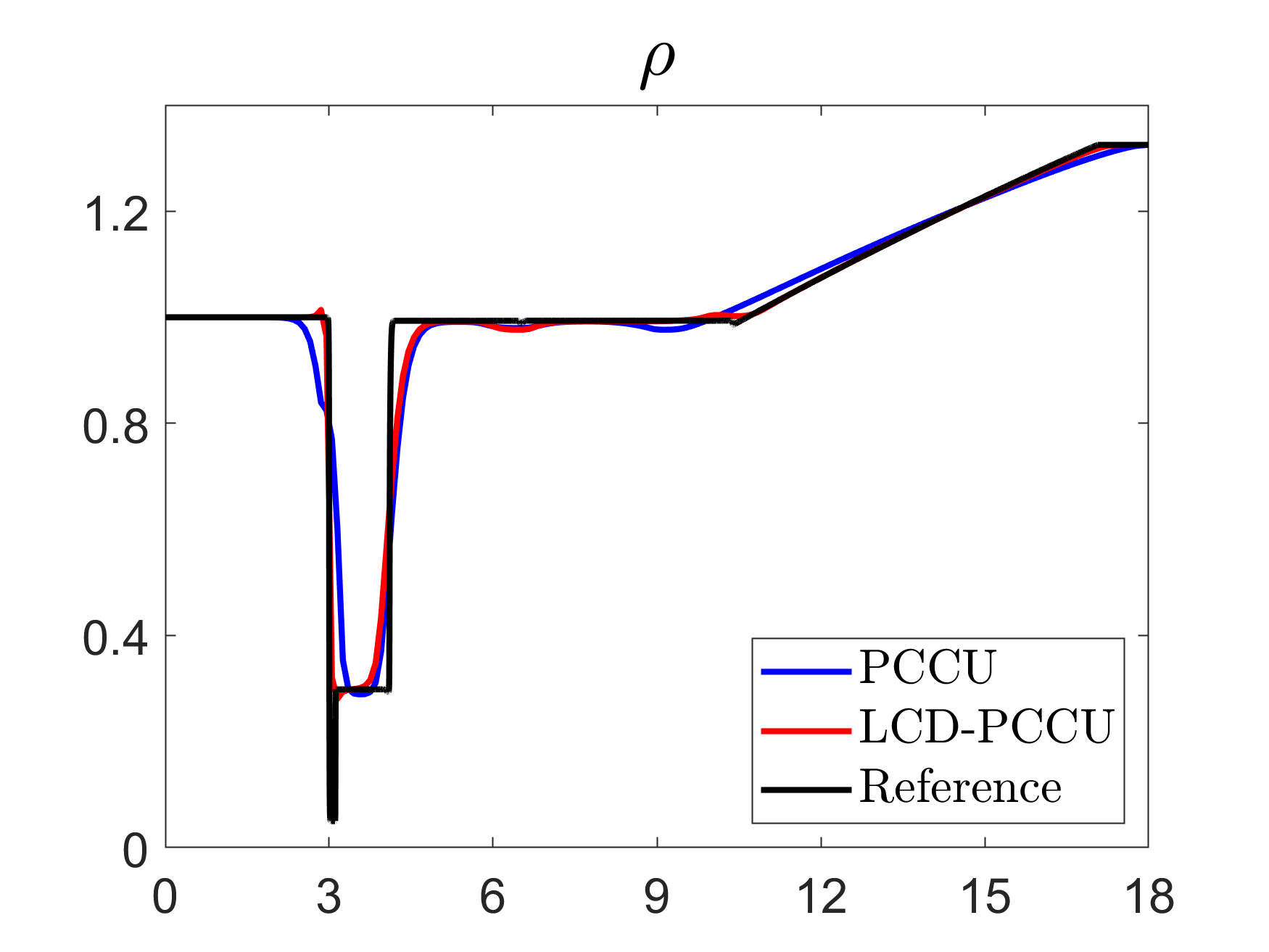}\hspace*{1.0cm}
            \includegraphics[trim=0.9cm 0.4cm 0.7cm 0.4cm, clip, width=5cm]{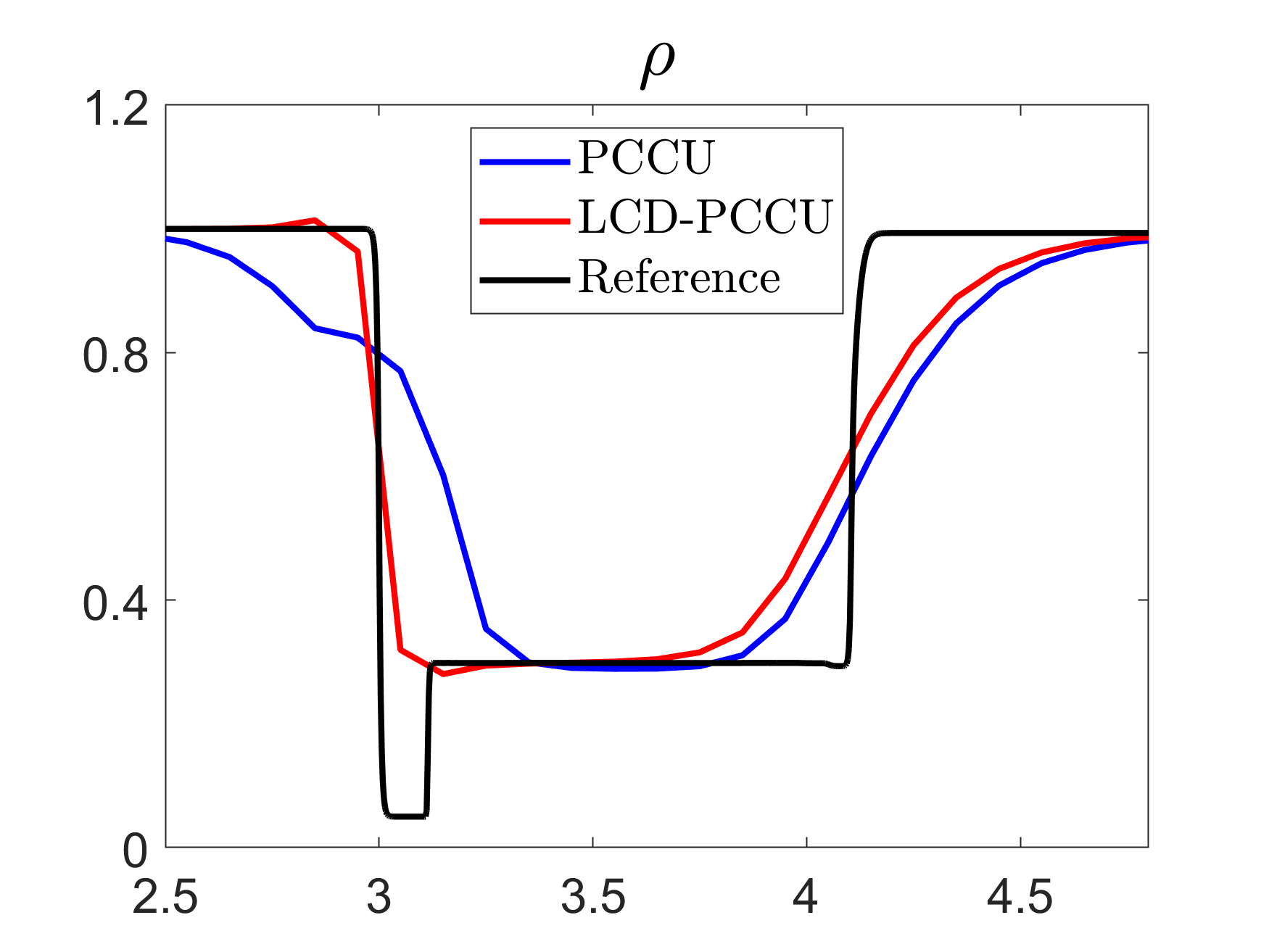}}
\caption{\sf Example 2: Density $\rho$ computed by the PCCU and LCD-PCCU schemes (left) and zoom at $x\in[2.5,4.8]$ (right).\label{fig2}}
\end{figure}

We impose the free boundary conditions at the both ends of the computational domain $[0,18]$ and compute the numerical solutions until the
final time $t=0.045$ on a uniform mesh with $\dx=1/10$ by both the LCD-PCCU and PCCU schemes. The obtained results are presented in Figure
\ref{fig2} together with the reference solution computed by the PCCU scheme on a much finer mesh with $\dx=1/400$. As one can see, the
LCD-PCCU solution achieves sharper resolution than its PCCU counterpart, especially of the contact wave located at about $x=3$. 

\subsubsection{Two-Dimensional Examples}
We now consider three 2-D examples, in which we plot the Schlieren images of the magnitude of the density gradient field, $|\nabla\rho|$. To
this end, we use the shading function
\begin{equation*}
\exp\bigg(-\frac{80|\nabla\rho|}{\max(|\nabla\rho|)}\bigg),
\end{equation*}
where the numerical derivatives of the density are computed using standard central differencing.

\subsubsection*{Example 3---Shock-Helium Bubble Interaction}
\begin{figure}[ht!]
\centerline{\includegraphics[trim=1.0cm 0.2cm 1.0cm 0.4cm, clip, width=6cm]{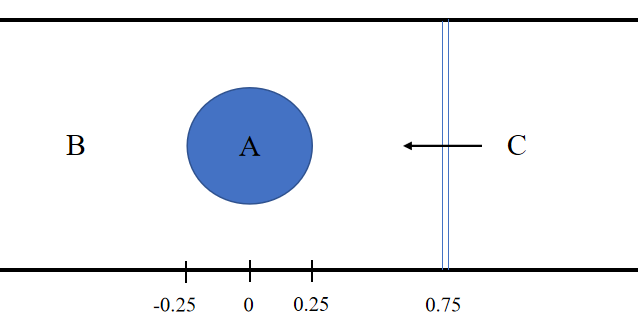}}
\caption{\sf Initial setting for the 2-D numerical examples.\label{fig44a}}
\end{figure}

In the first 2-D example taken from \cite{Quirk1996,chertock7}, a shock wave in the air hits the light resting bubble which contains helium.
We take the following initial conditions:
\begin{equation*}
(\rho,u,v,p;\gamma,p_\infty)=\begin{cases}(4/29,0,0,1;5/3,0)&\mbox{in region A},\\(1,0,0,1;1.4,0)&\mbox{in region B},\\
(4/3,-0.3535,0,1.5;1.4,0)&\mbox{in region C},\end{cases}
\end{equation*}
where regions A, B, and C are outlined in Figure \ref{fig44a}, and the computational domain is $[-3,1]\times[-0.5,0.5]$. We impose the solid
wall boundary conditions on the top and bottom and the free boundary conditions on the left and right edges of the computational domain.

\begin{figure}[ht!]
\centerline{\includegraphics[trim=0.5cm 1.9cm 1.cm 1.5cm, clip, width=6.cm]{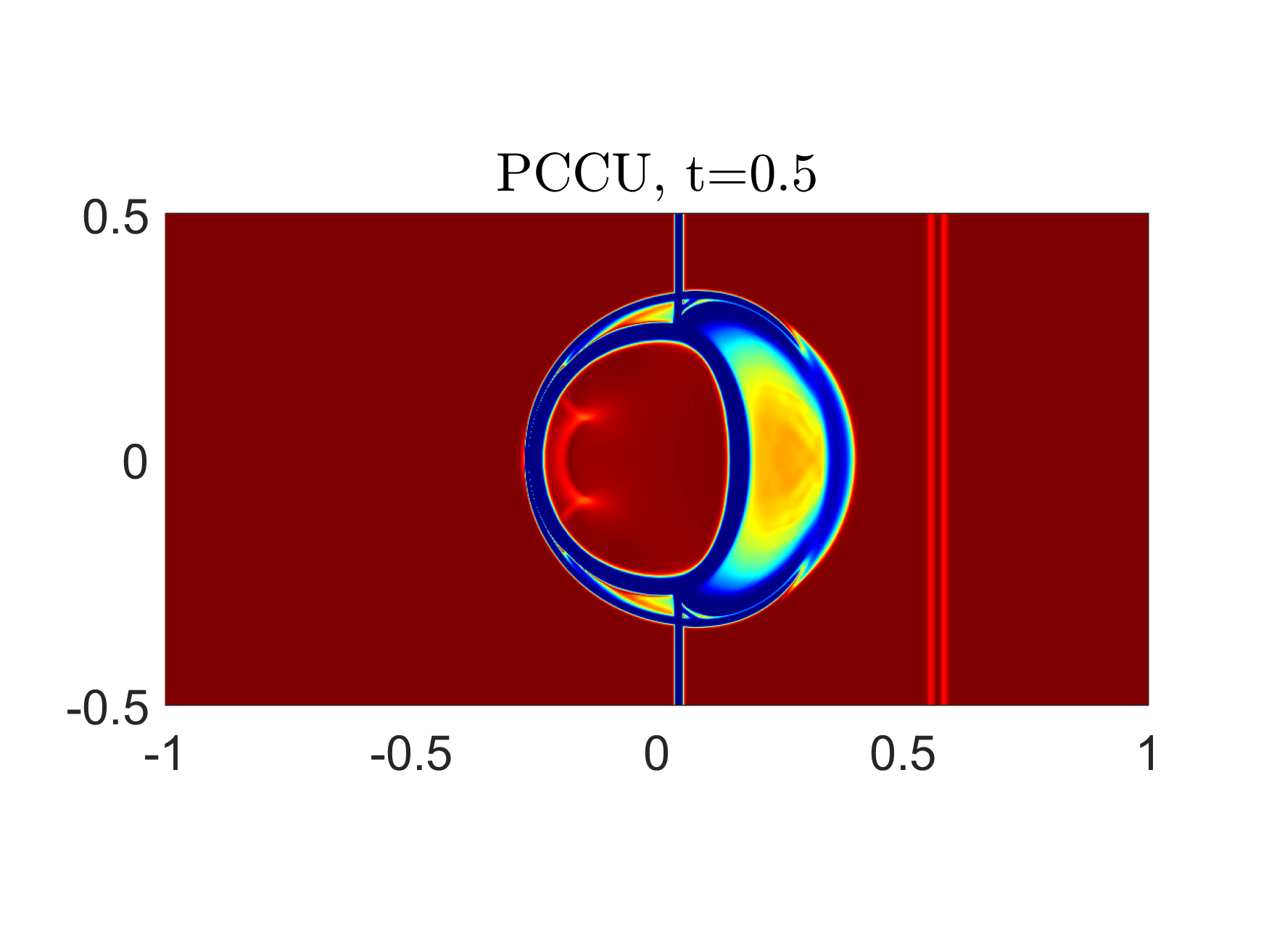}\hspace*{1.cm}
            \includegraphics[trim=0.5cm 1.9cm 1.cm 1.5cm, clip, width=6.cm]{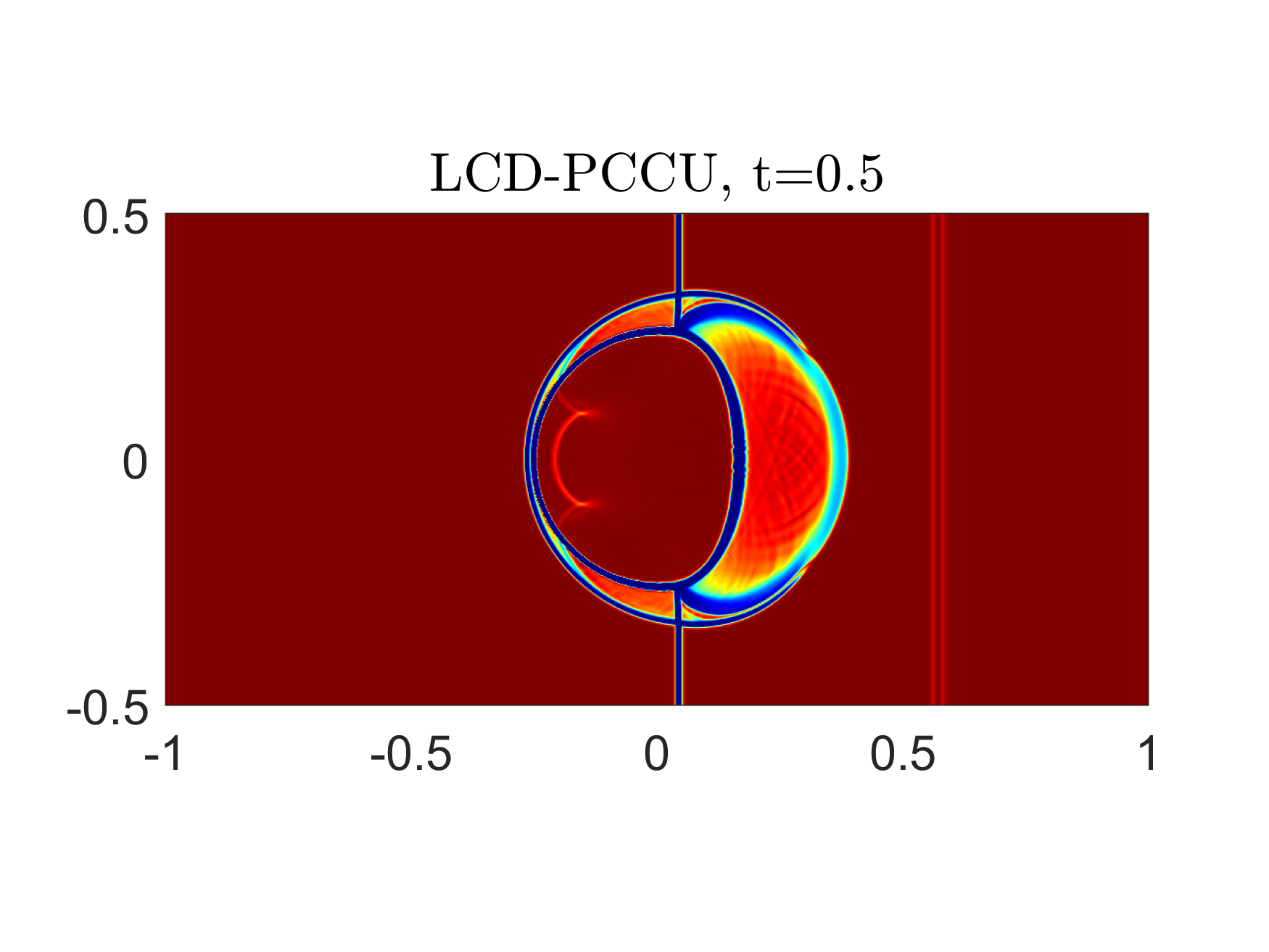}}
\centerline{\includegraphics[trim=0.5cm 1.9cm 1.cm 1.5cm, clip, width=6.cm]{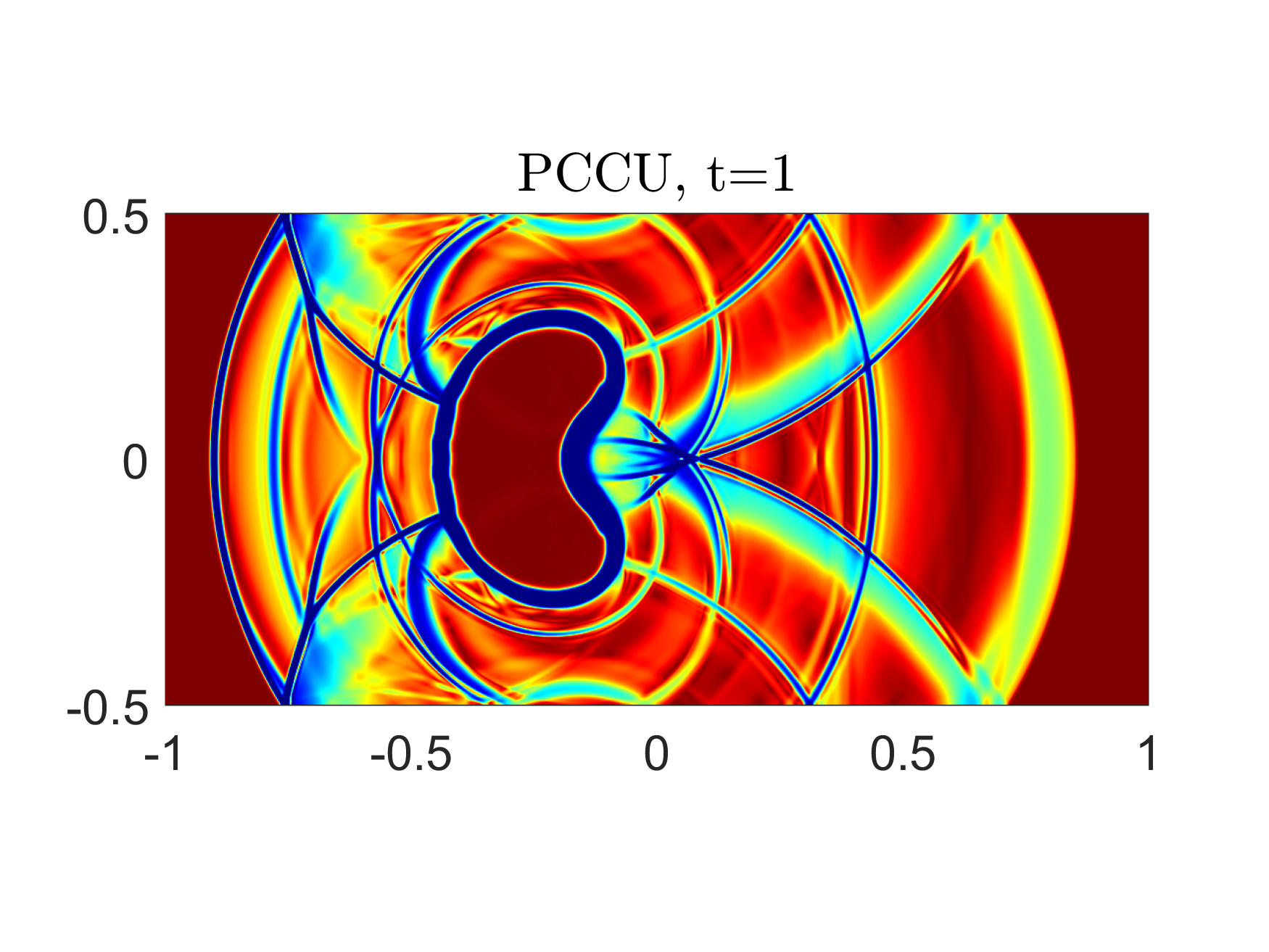}\hspace*{1.cm}
            \includegraphics[trim=0.5cm 1.9cm 1.cm 1.5cm, clip, width=6.cm]{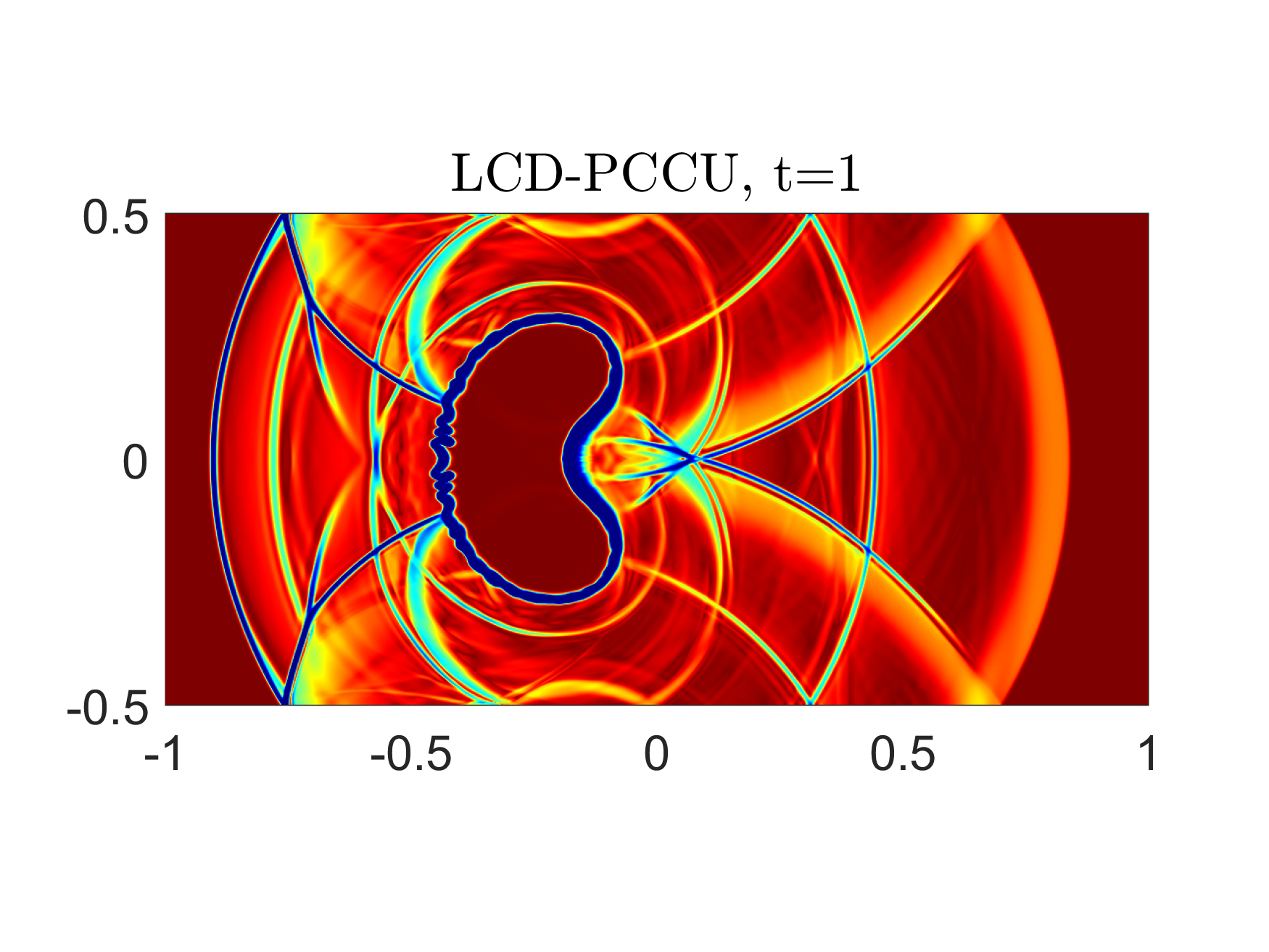}}
\centerline{\includegraphics[trim=0.5cm 1.9cm 1.cm 1.5cm, clip, width=6.cm]{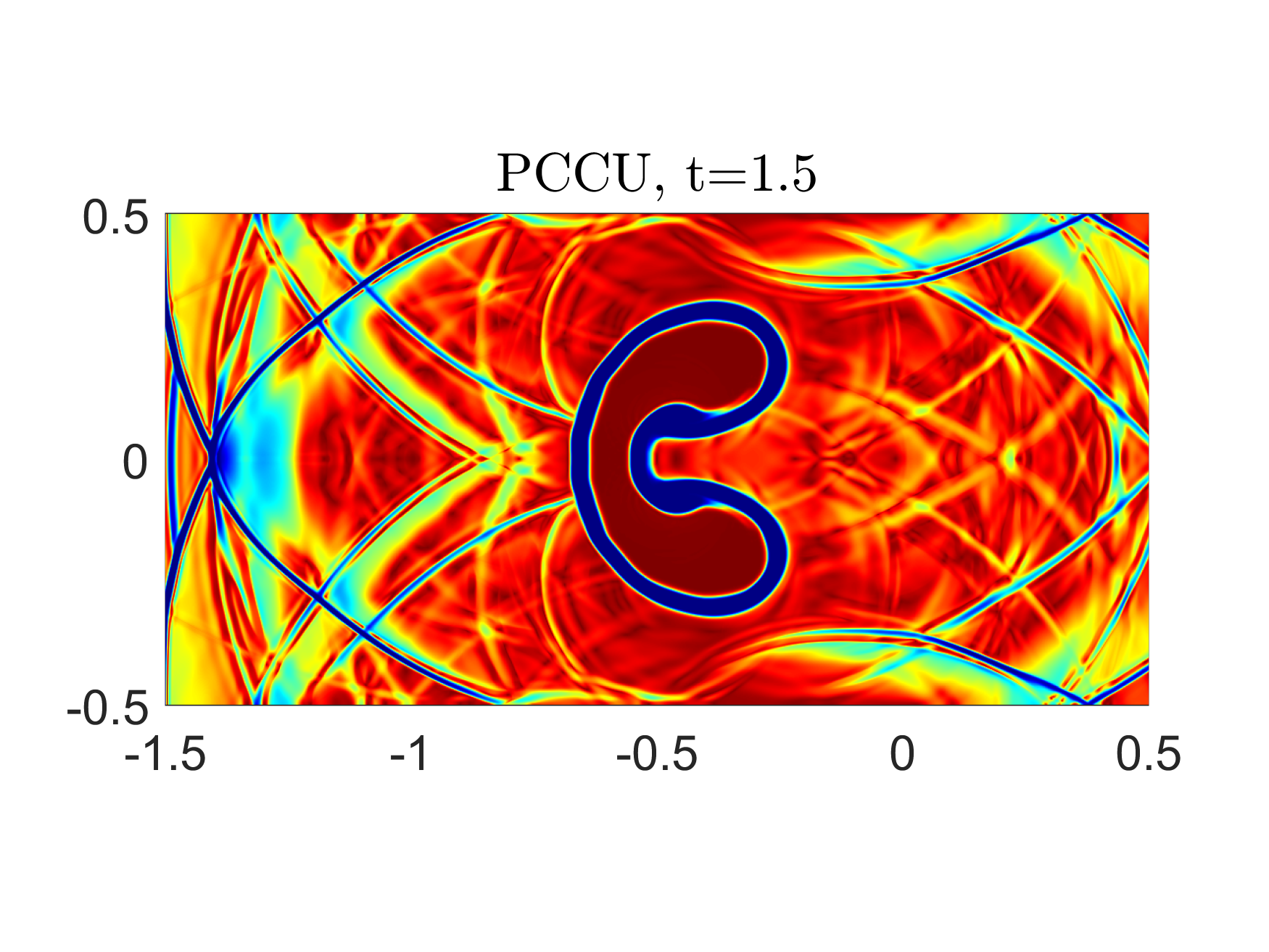}\hspace*{1.cm}
            \includegraphics[trim=0.5cm 1.9cm 1.cm 1.5cm, clip, width=6.cm]{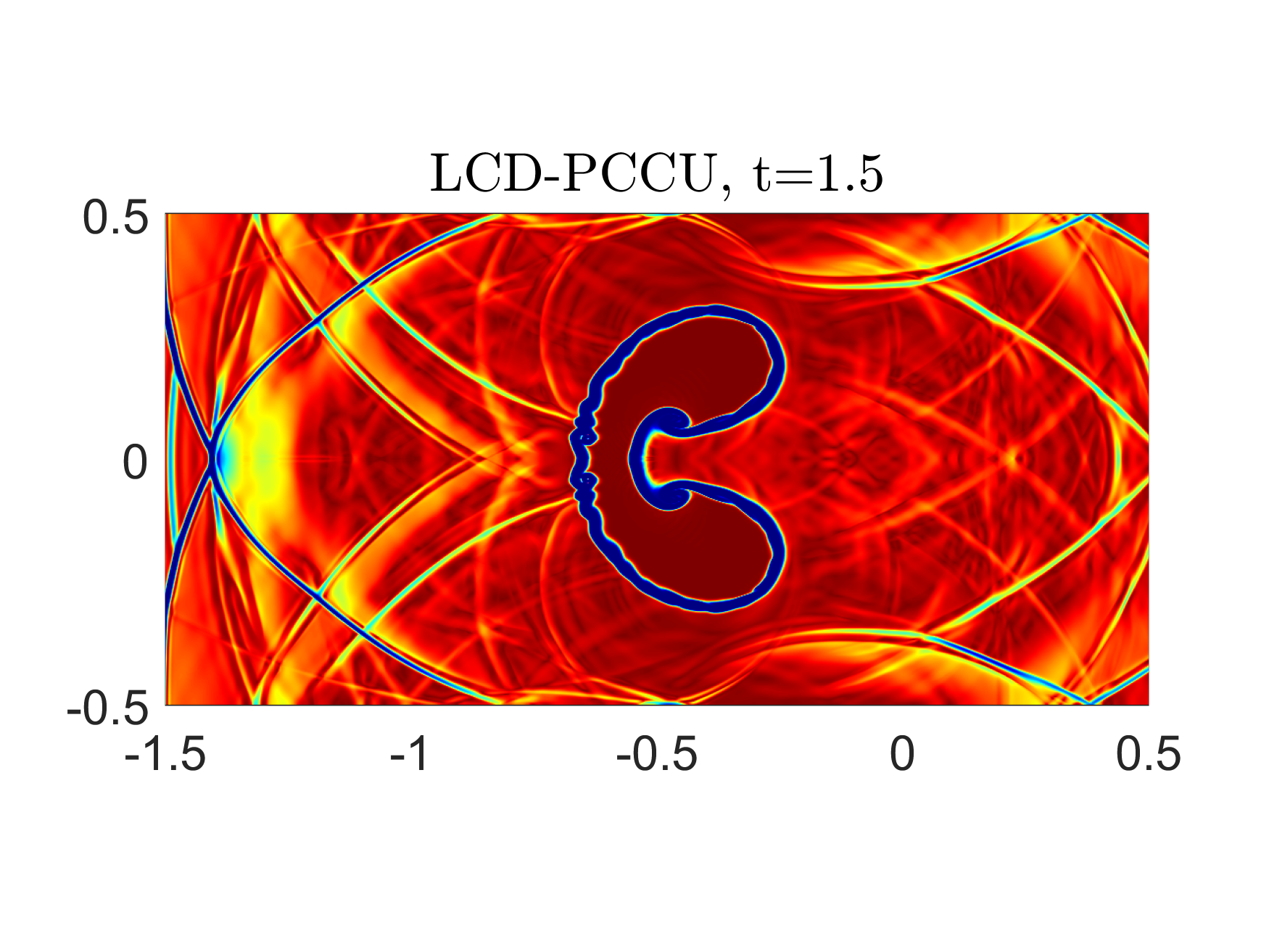}}
\centerline{\includegraphics[trim=0.5cm 1.9cm 1.cm 1.5cm, clip, width=6.cm]{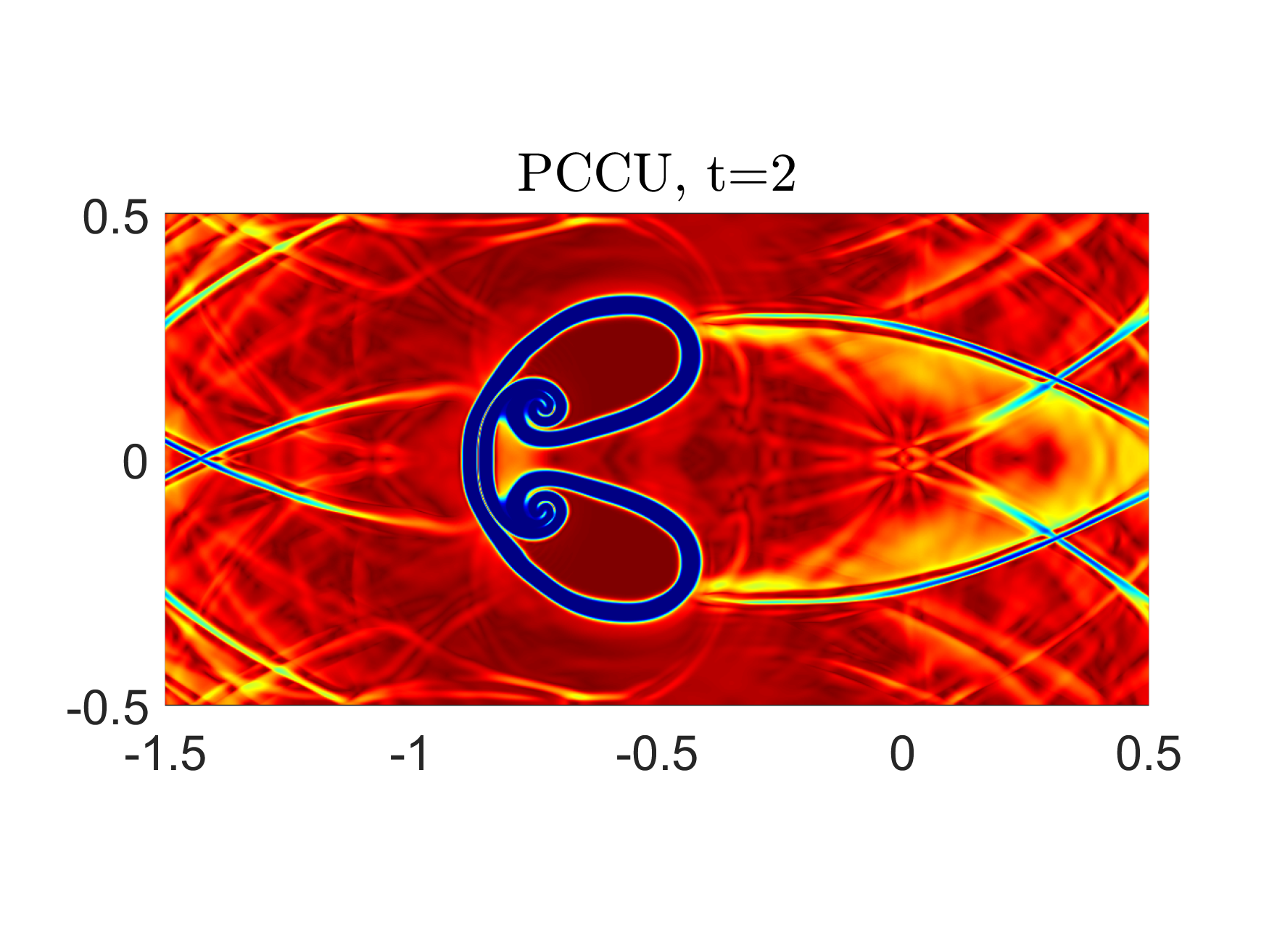}\hspace*{1.cm}
            \includegraphics[trim=0.5cm 1.9cm 1.cm 1.5cm, clip, width=6.cm]{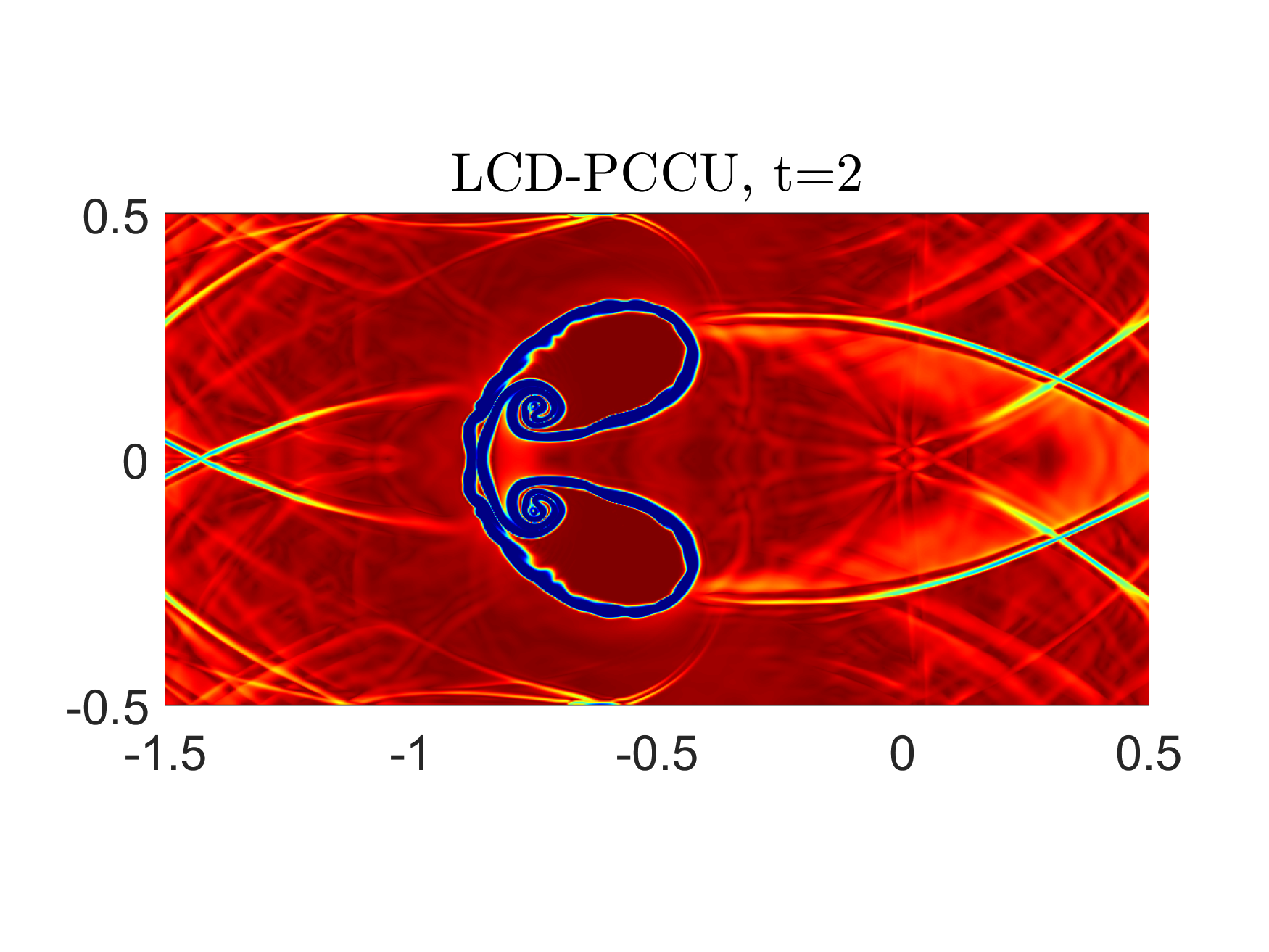}}
\centerline{\includegraphics[trim=0.5cm 1.9cm 1.cm 1.5cm, clip, width=6.cm]{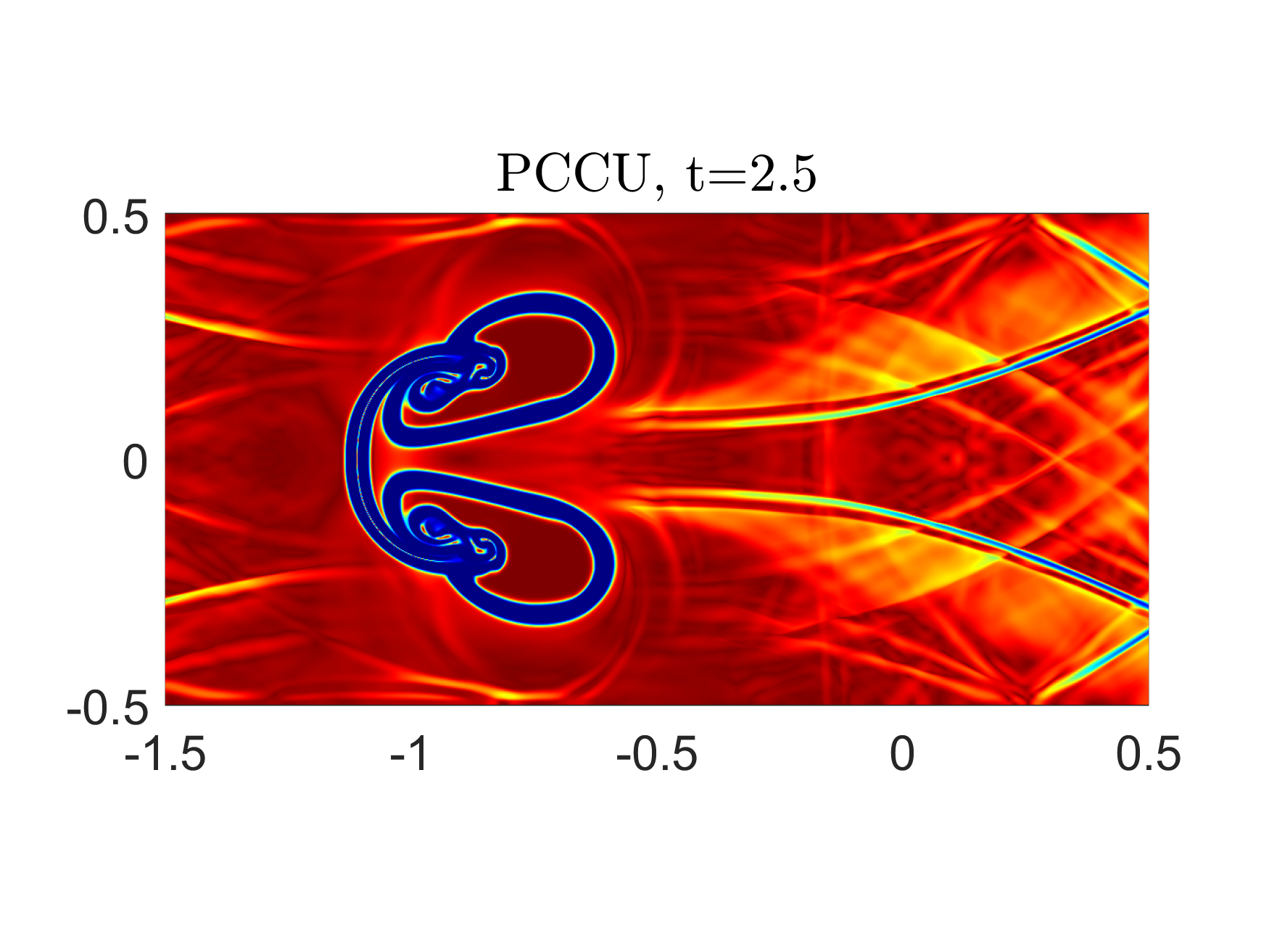}\hspace*{1.cm}
            \includegraphics[trim=0.5cm 1.9cm 1.cm 1.5cm, clip, width=6.cm]{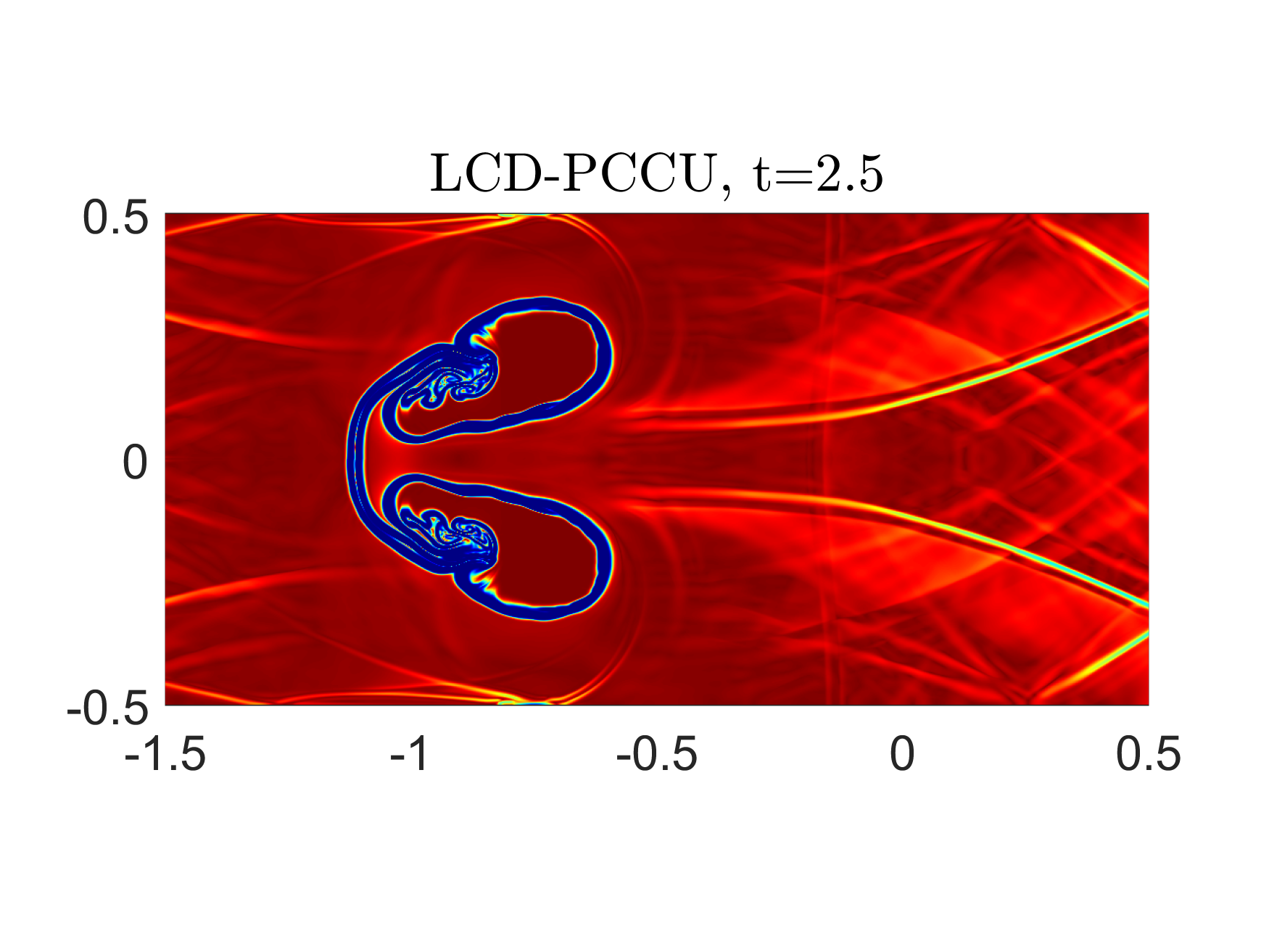}}
\centerline{\includegraphics[trim=0.5cm 1.9cm 1.cm 1.5cm, clip, width=6.cm]{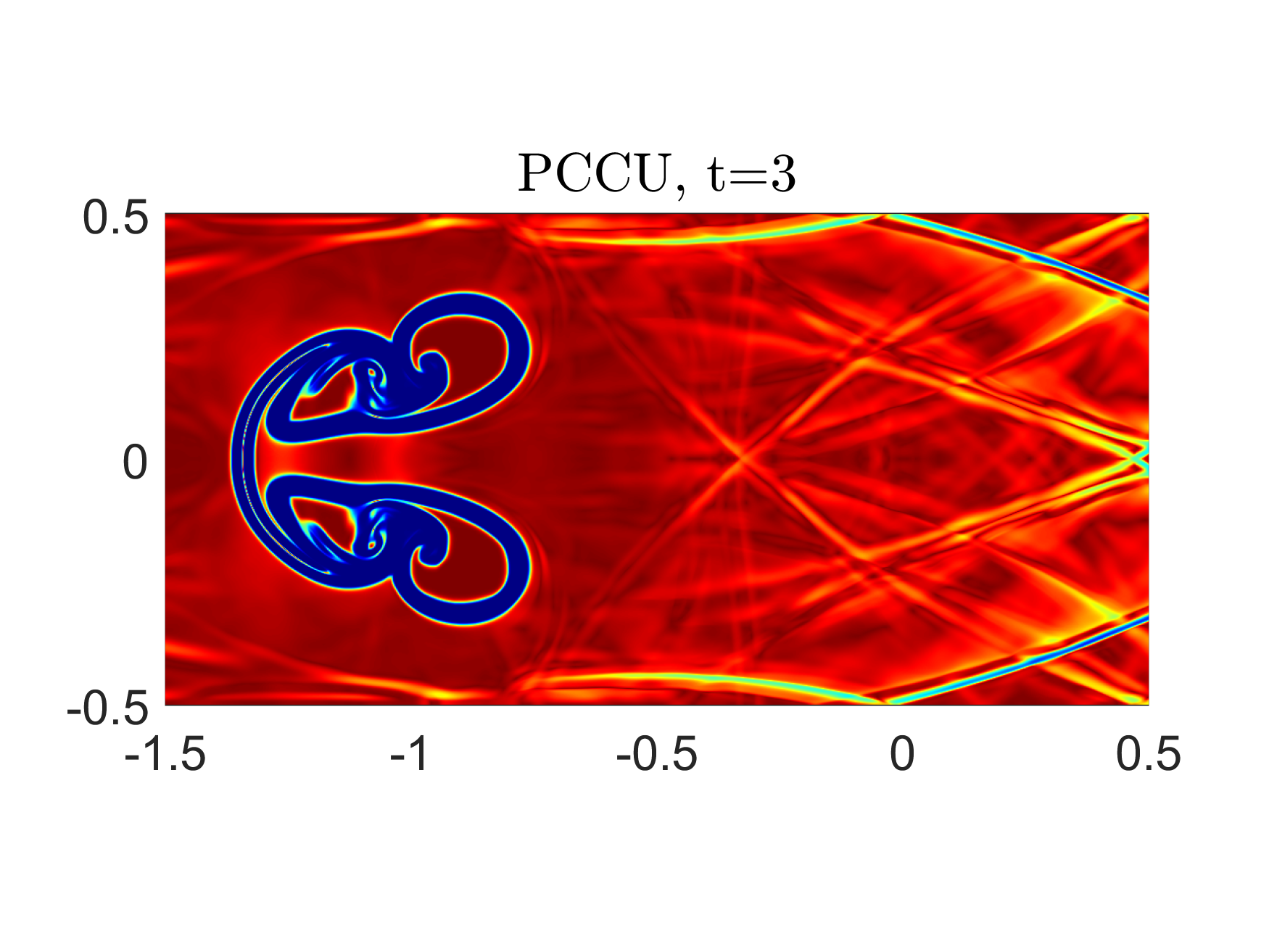}\hspace*{1.cm}
            \includegraphics[trim=0.5cm 1.9cm 1.cm 1.5cm, clip, width=6.cm]{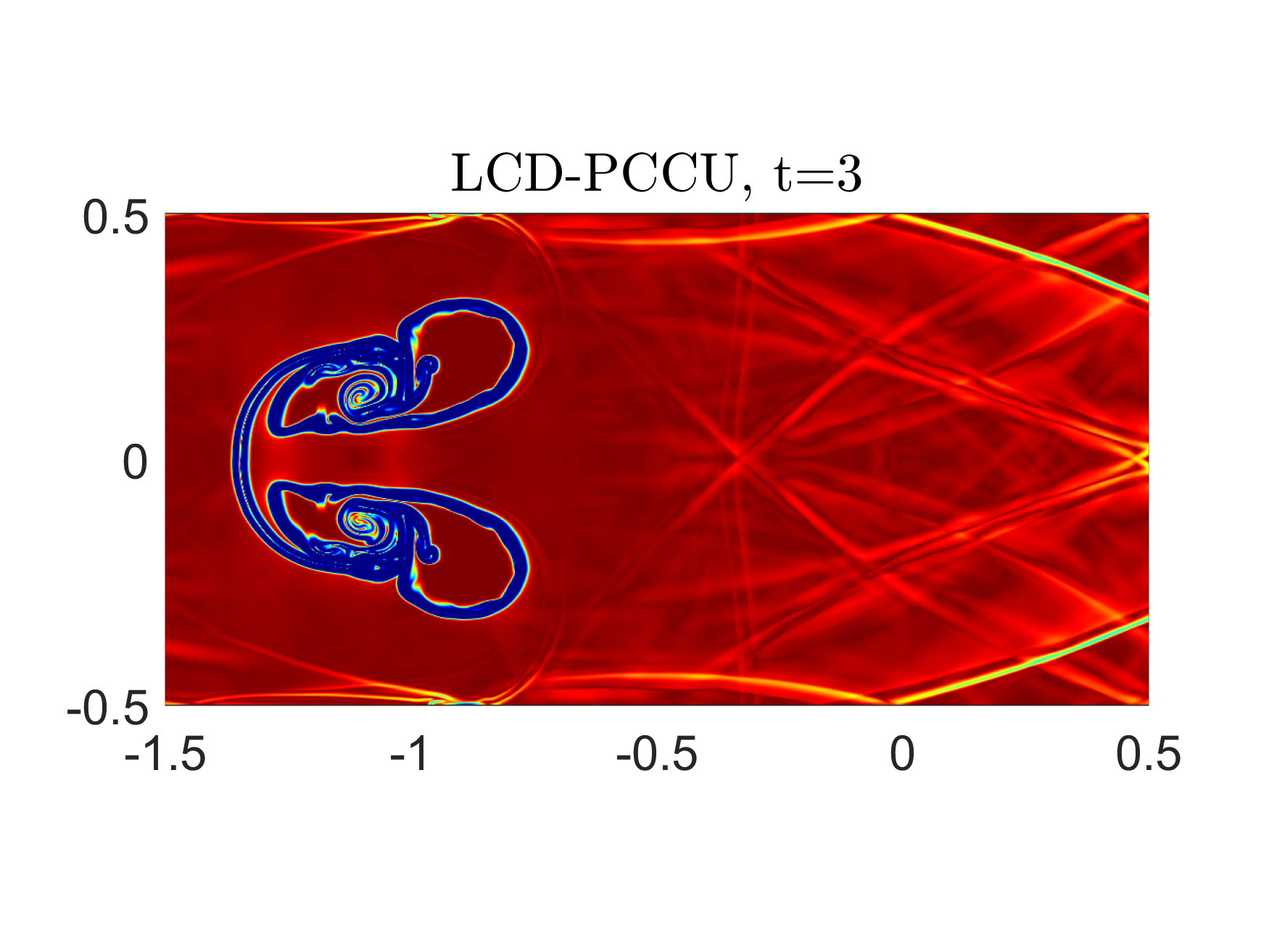}}                       
\caption{\sf Example 3: Shock-helium bubble interaction by the PCCU (left column) and LCD-PCCU (right column) schemes at different times.
\label{fig43a}}
\end{figure}

We compute the numerical solutions until the final time $t=3$ on a uniform mesh with $\dx=\dy=1/500$. In Figure \ref{fig43a}, we present
different stages of the shock-bubble interaction computed by the PCCU and LCD-PCCU schemes. One can clearly see that the bubble interface
develops very complex structures after the bubble is hit by the shock, especially at large times. We note that the obtained results are in a
good agreement with the experimental findings presented in \cite{Haas1987} and the computational results presented in
\cite{Quirk1996,chertock7, CCK_21}. Notably, at an early time $t=0.5$, the resolution of the bubble interface significantly improves with
the implementation of the LCD-PCCU scheme. As time progresses, instabilities arise within the interface, which are smeared by a more
dissipative PCCU scheme. The differences in resolution become more pronounced at later times $t=2$ and $2.5$, and especially at $t=3$,
indicating that the LCD-PCCU scheme clearly outperforms the PCCU one.

\subsubsection*{Example 4---Shock-R22 Bubble Interaction}
In the second 2-D example also taken from \cite{Quirk1996,chertock7}, a shock wave in the air hits the heavy resting bubble which contains
R22 gas. The initial conditions are
\begin{equation*}
(\rho,u,v,p;\gamma,\pi_\infty)=\begin{cases}(3.1538,0,0,1;1.249,0)&\mbox{in region A},\\(1,0,0,1;1.4,0)&\mbox{in region B},\\
(4/3,-0.3535,0,1.5;1.4,0)&\mbox{in region C},\end{cases}
\end{equation*}
and the rest of the settings are the same as in Example 3.

\begin{figure}[ht!]
\centerline{\includegraphics[trim=0.5cm 1.9cm 1.cm 1.5cm, clip, width=6.cm]{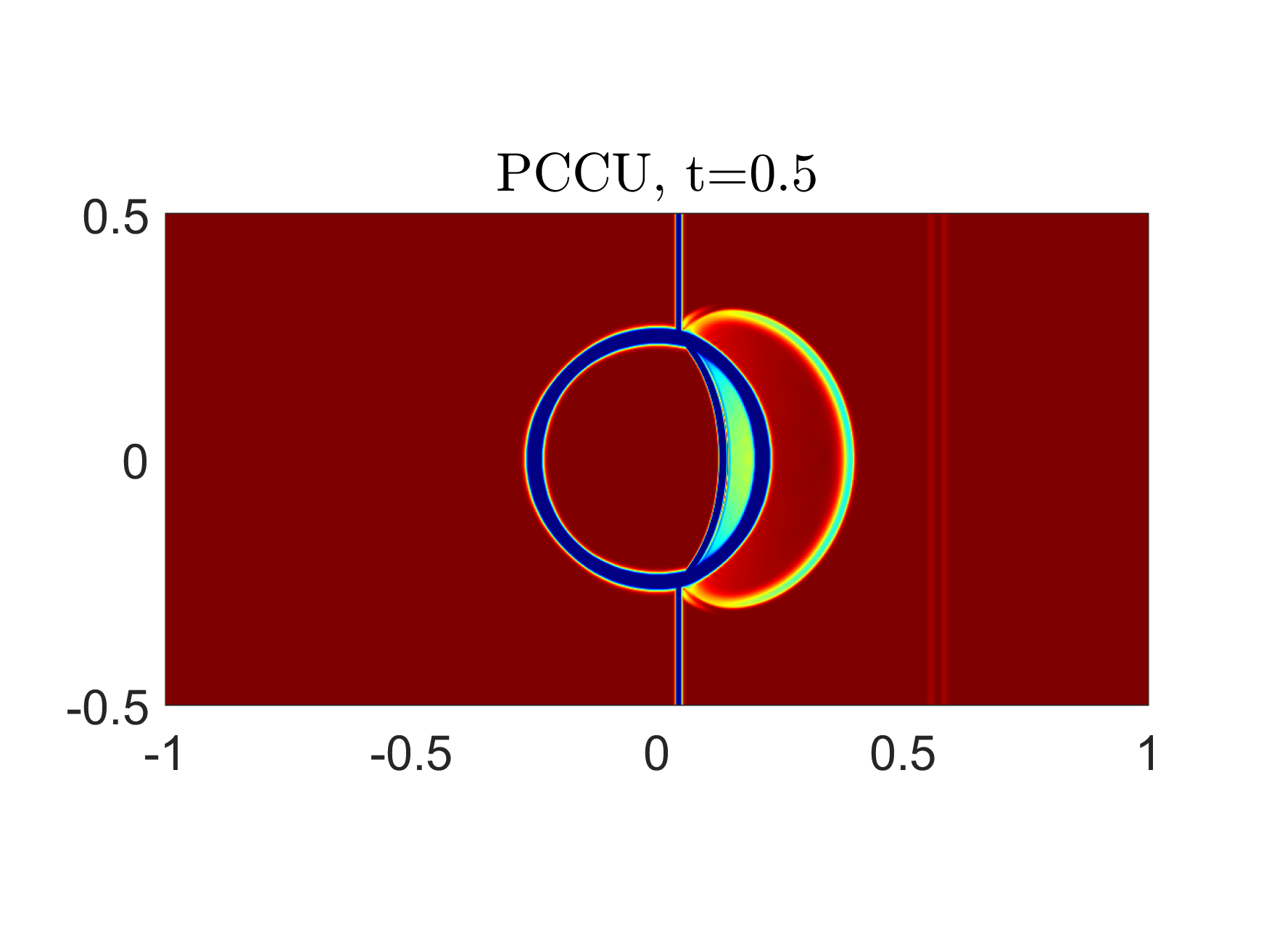}\hspace*{1.cm}
            \includegraphics[trim=0.5cm 1.9cm 1.cm 1.5cm, clip, width=6.cm]{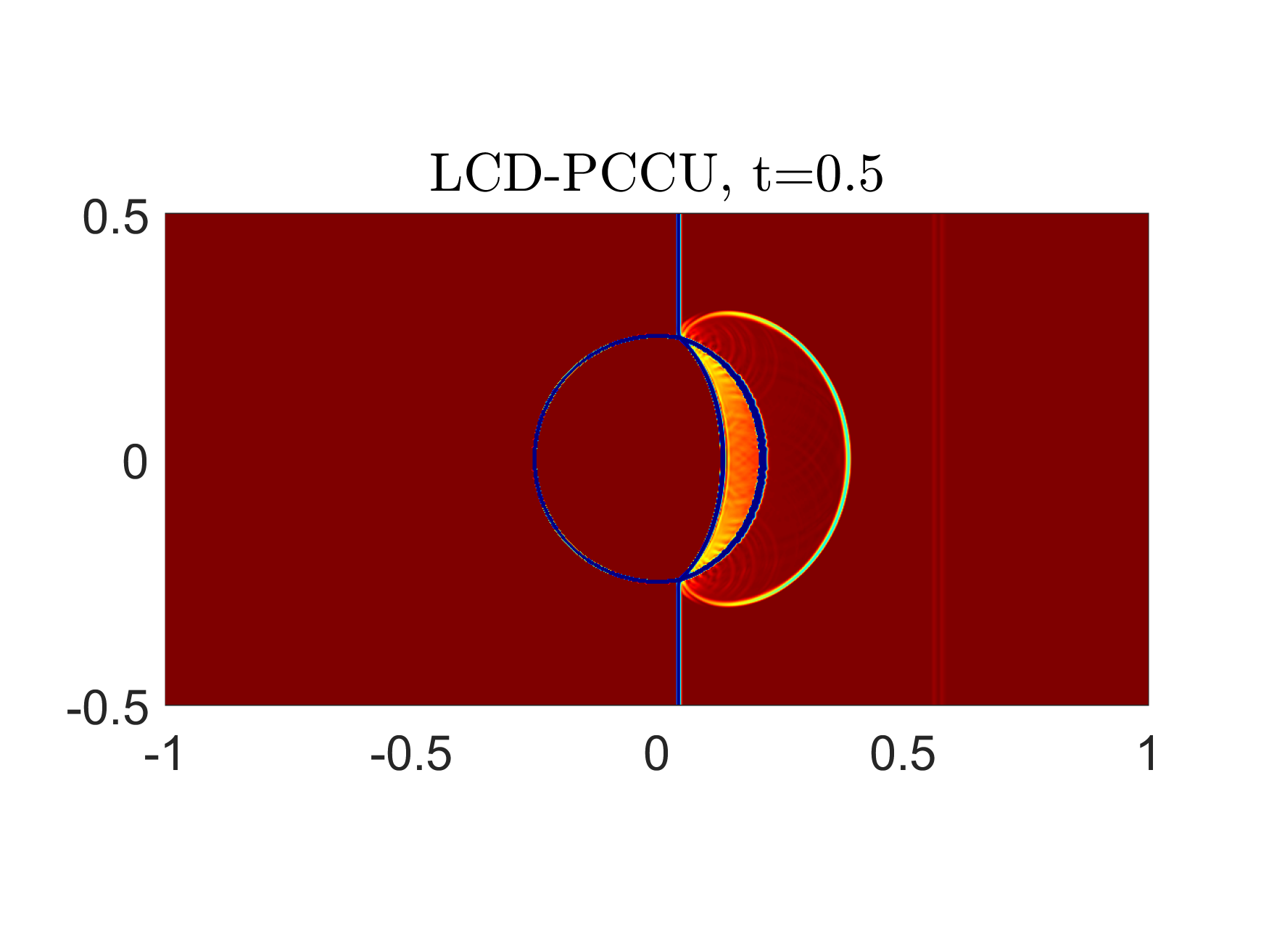}}
\centerline{\includegraphics[trim=0.5cm 1.9cm 1.cm 1.5cm, clip, width=6.cm]{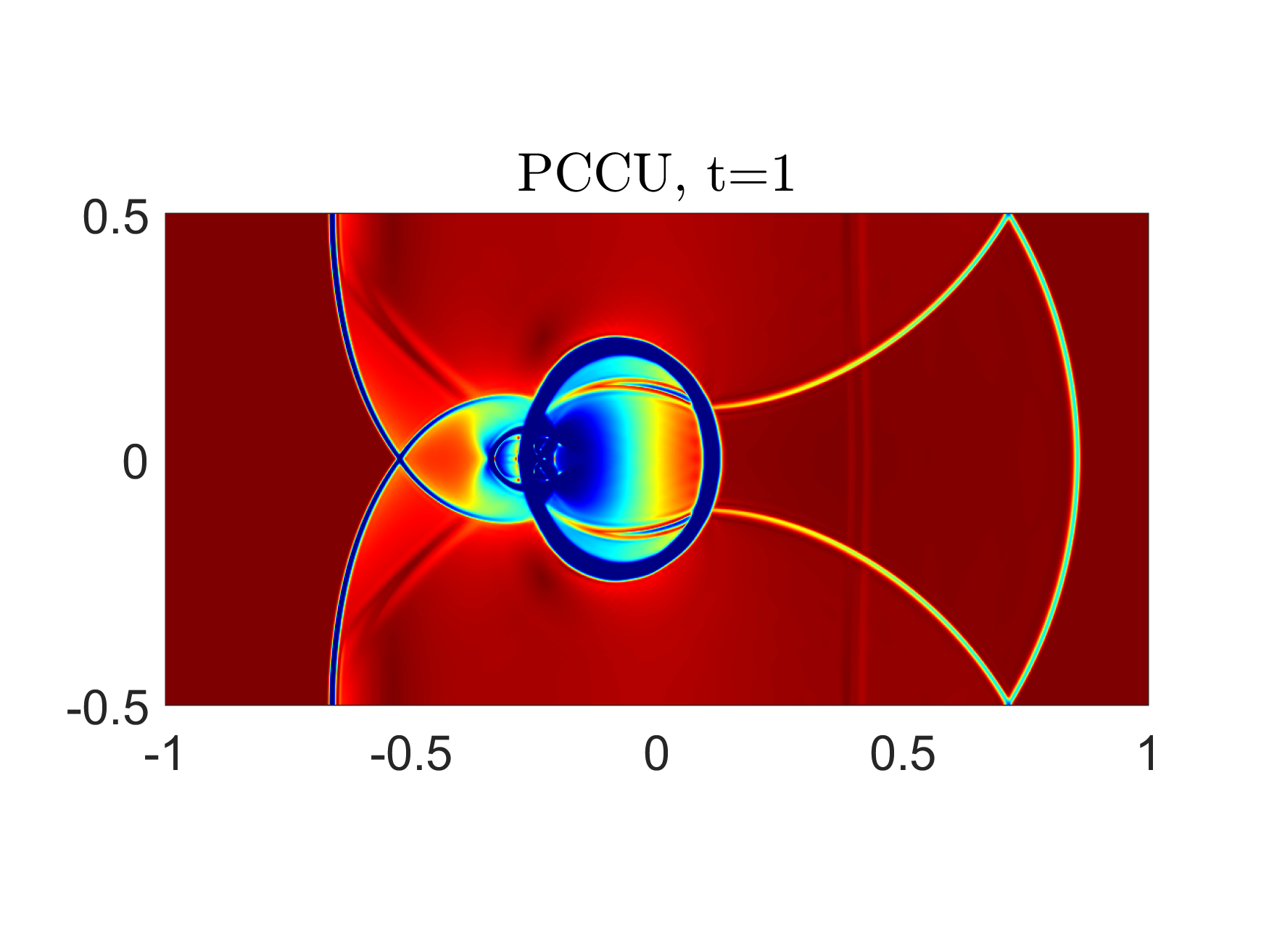}\hspace*{1.cm}
            \includegraphics[trim=0.5cm 1.9cm 1.cm 1.5cm, clip, width=6.cm]{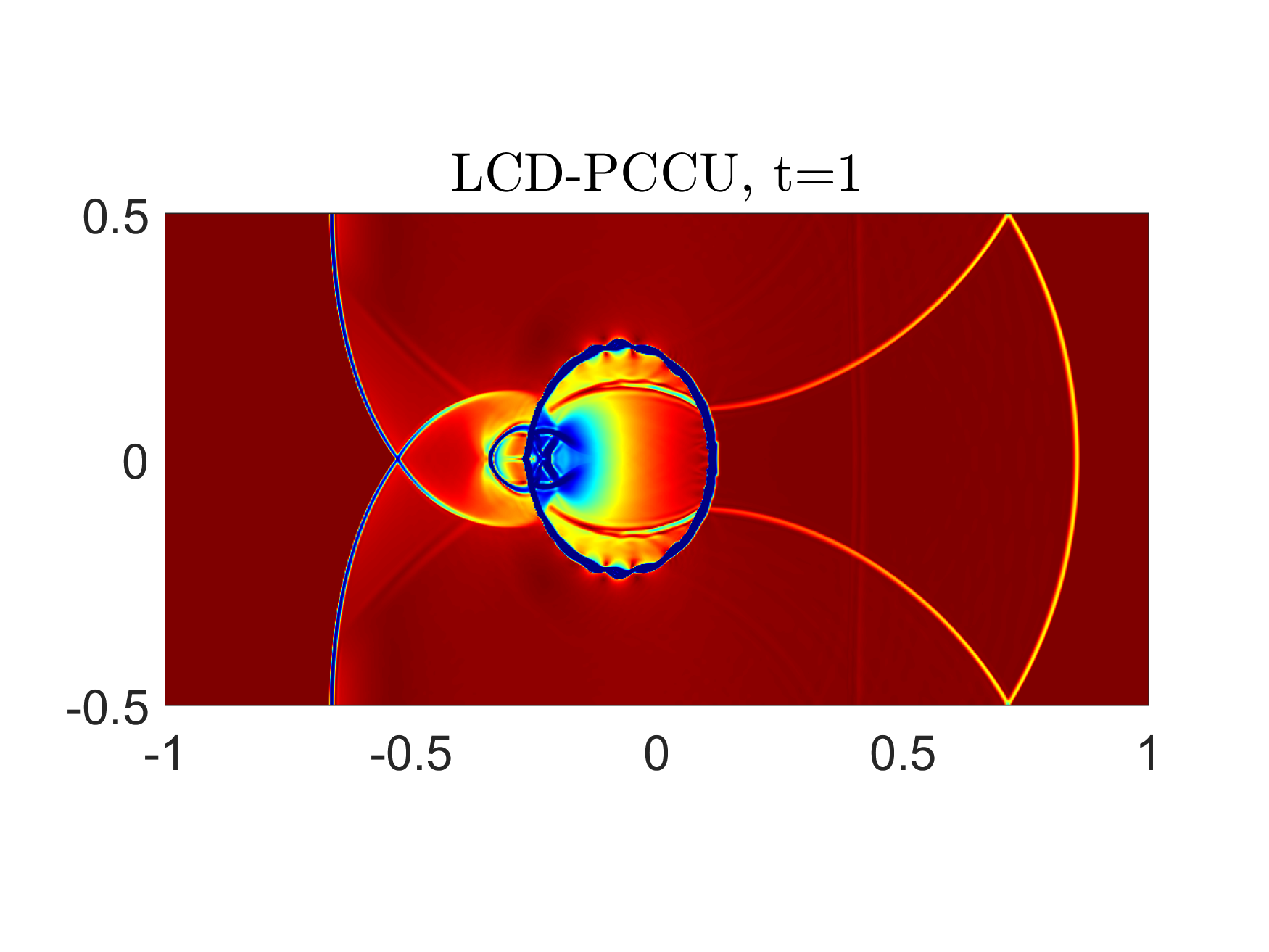}}
\centerline{\includegraphics[trim=0.5cm 1.9cm 1.cm 1.5cm, clip, width=6.cm]{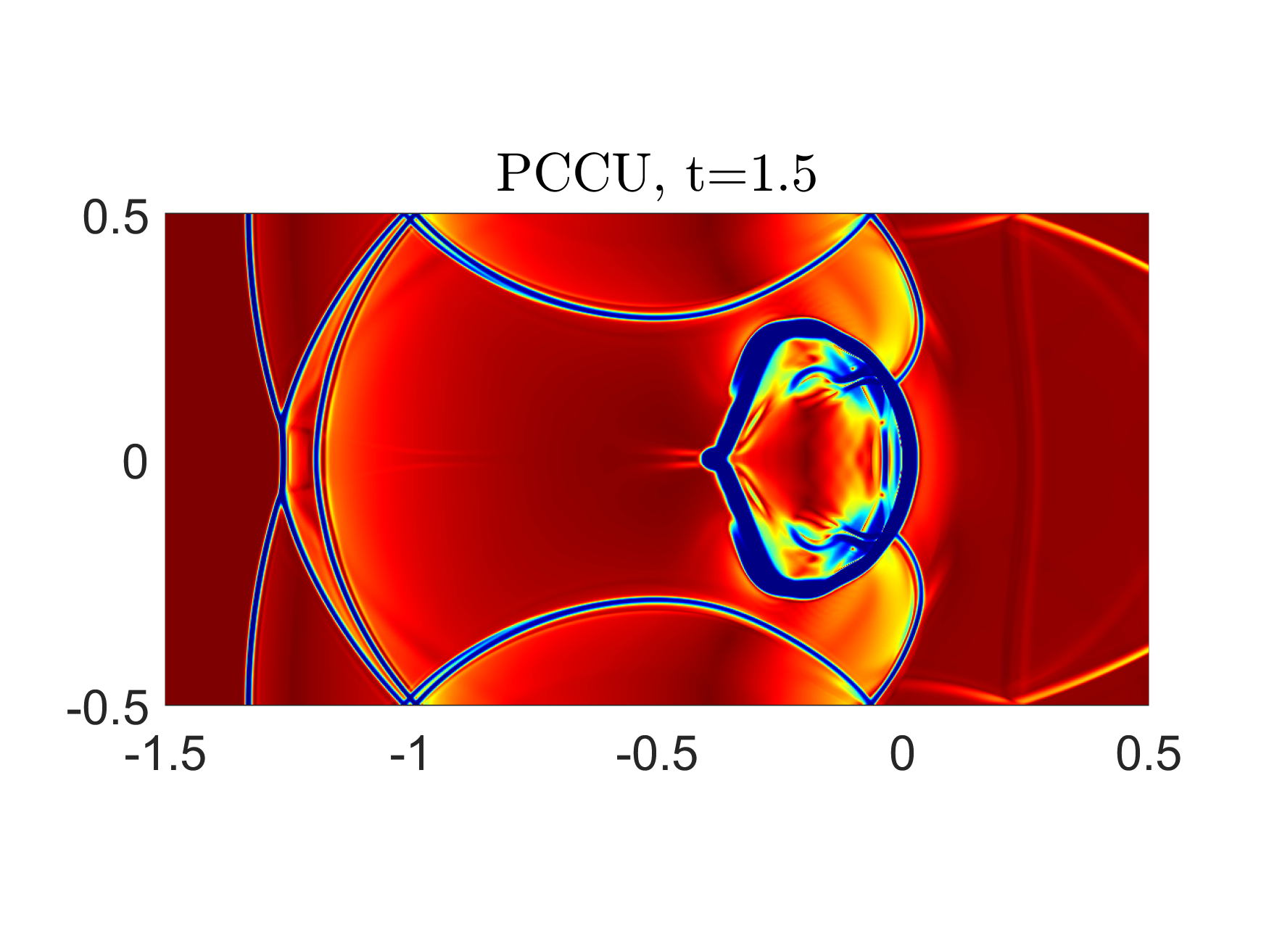}\hspace*{1.cm}
            \includegraphics[trim=0.5cm 1.9cm 1.cm 1.5cm, clip, width=6.cm]{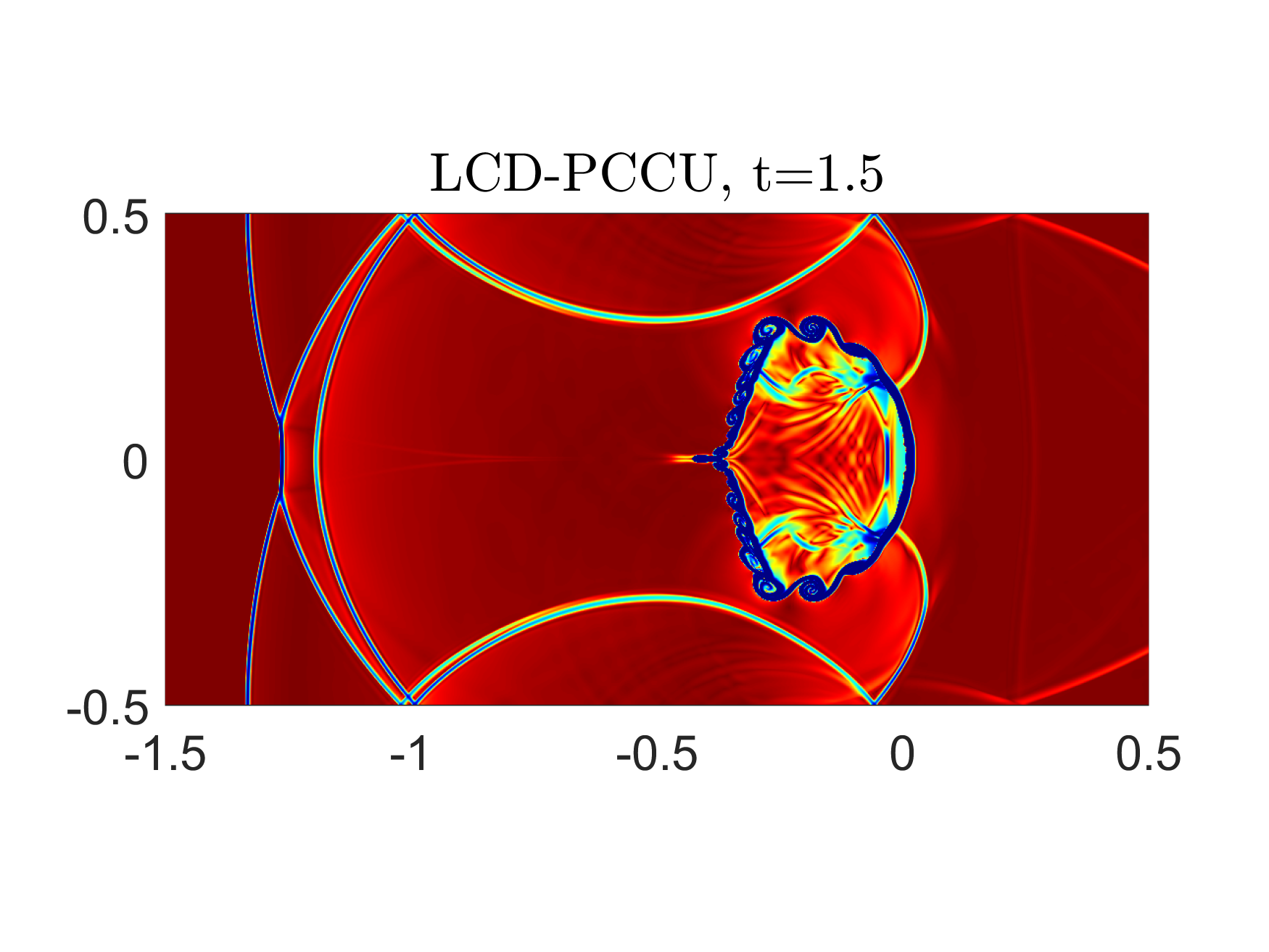}}
\centerline{\includegraphics[trim=0.5cm 1.9cm 1.cm 1.5cm, clip, width=6.cm]{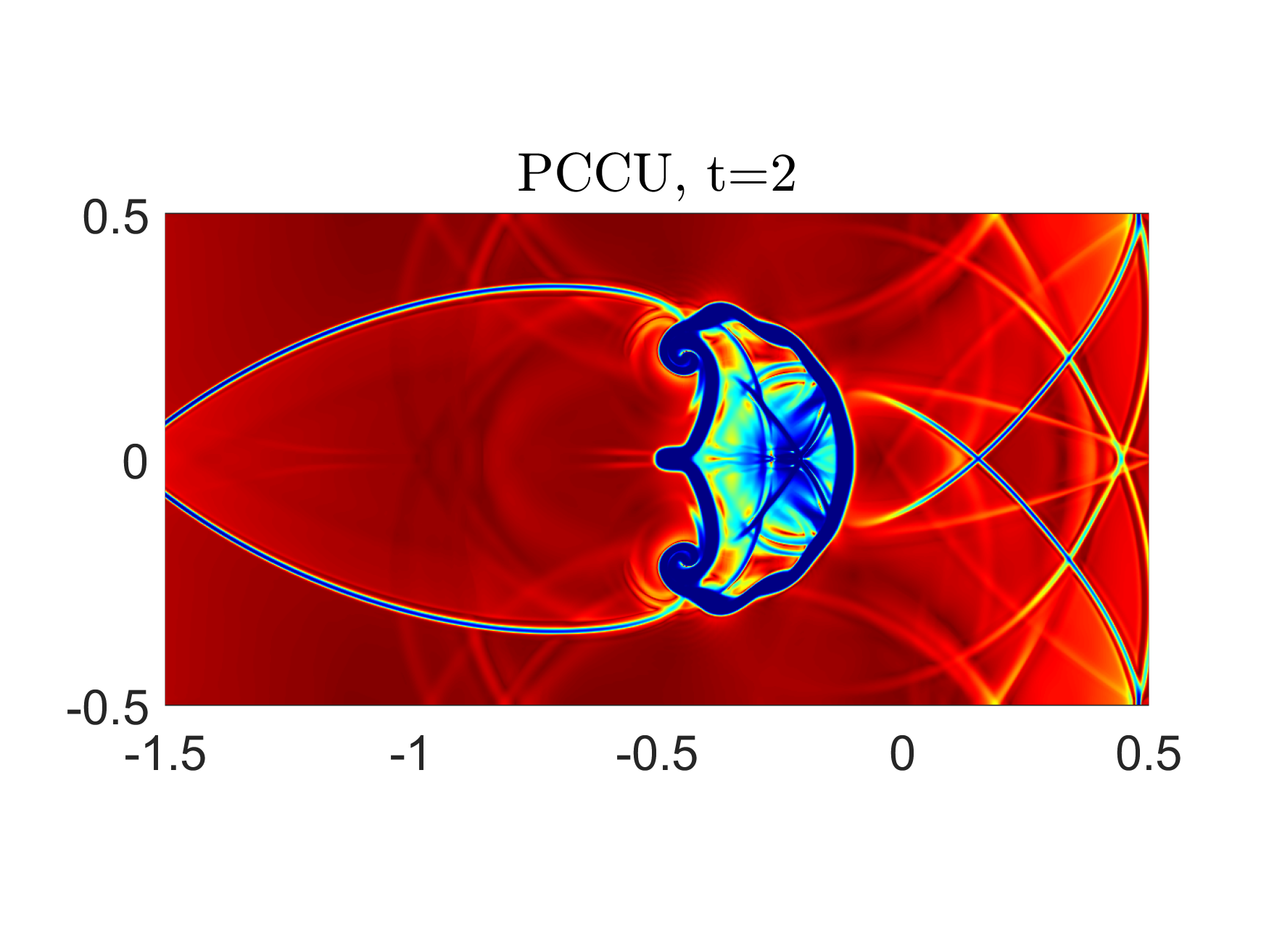}\hspace*{1.cm}
            \includegraphics[trim=0.5cm 1.9cm 1.cm 1.5cm, clip, width=6.cm]{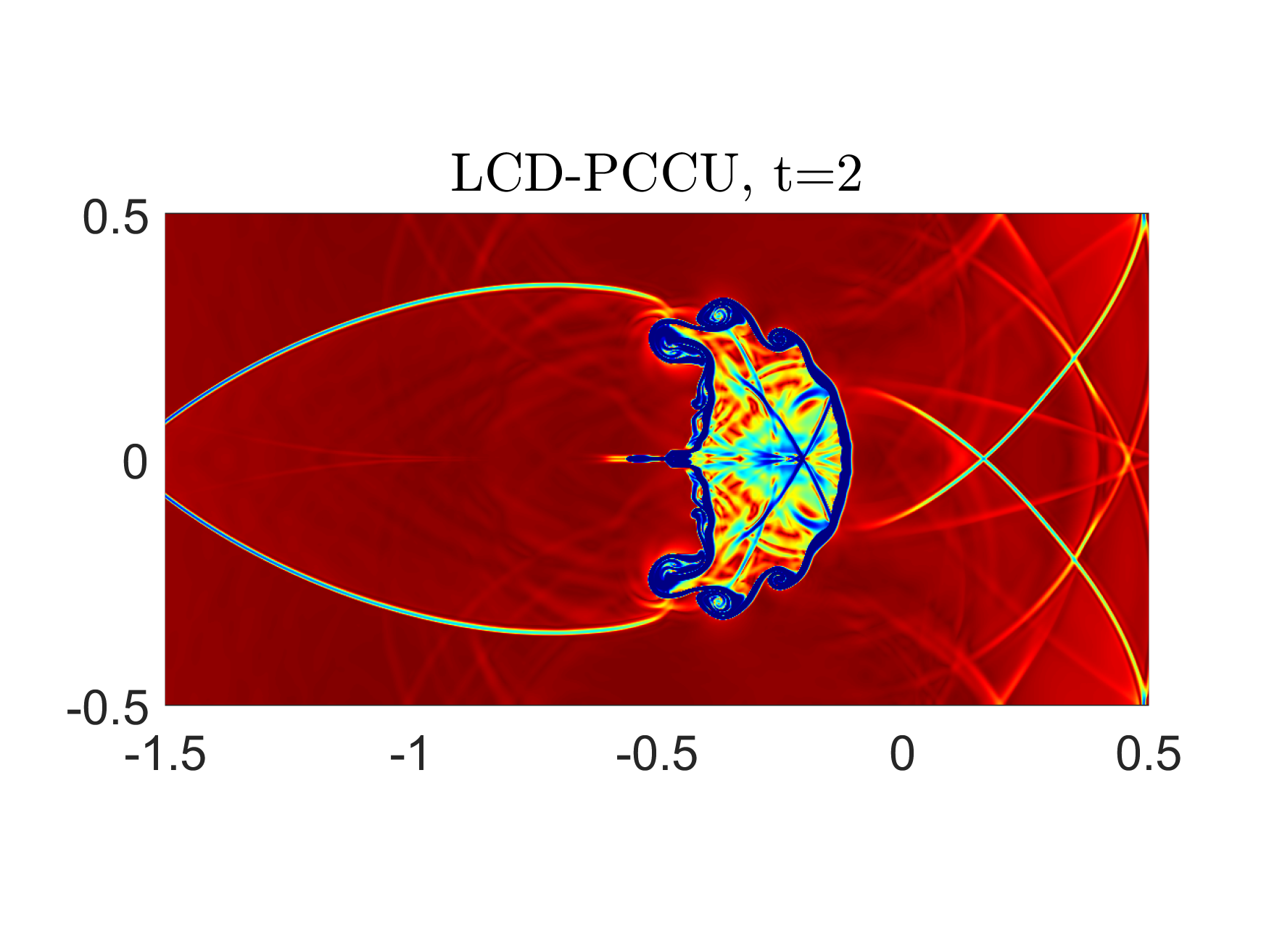}}
\centerline{\includegraphics[trim=0.5cm 1.9cm 1.cm 1.5cm, clip, width=6.cm]{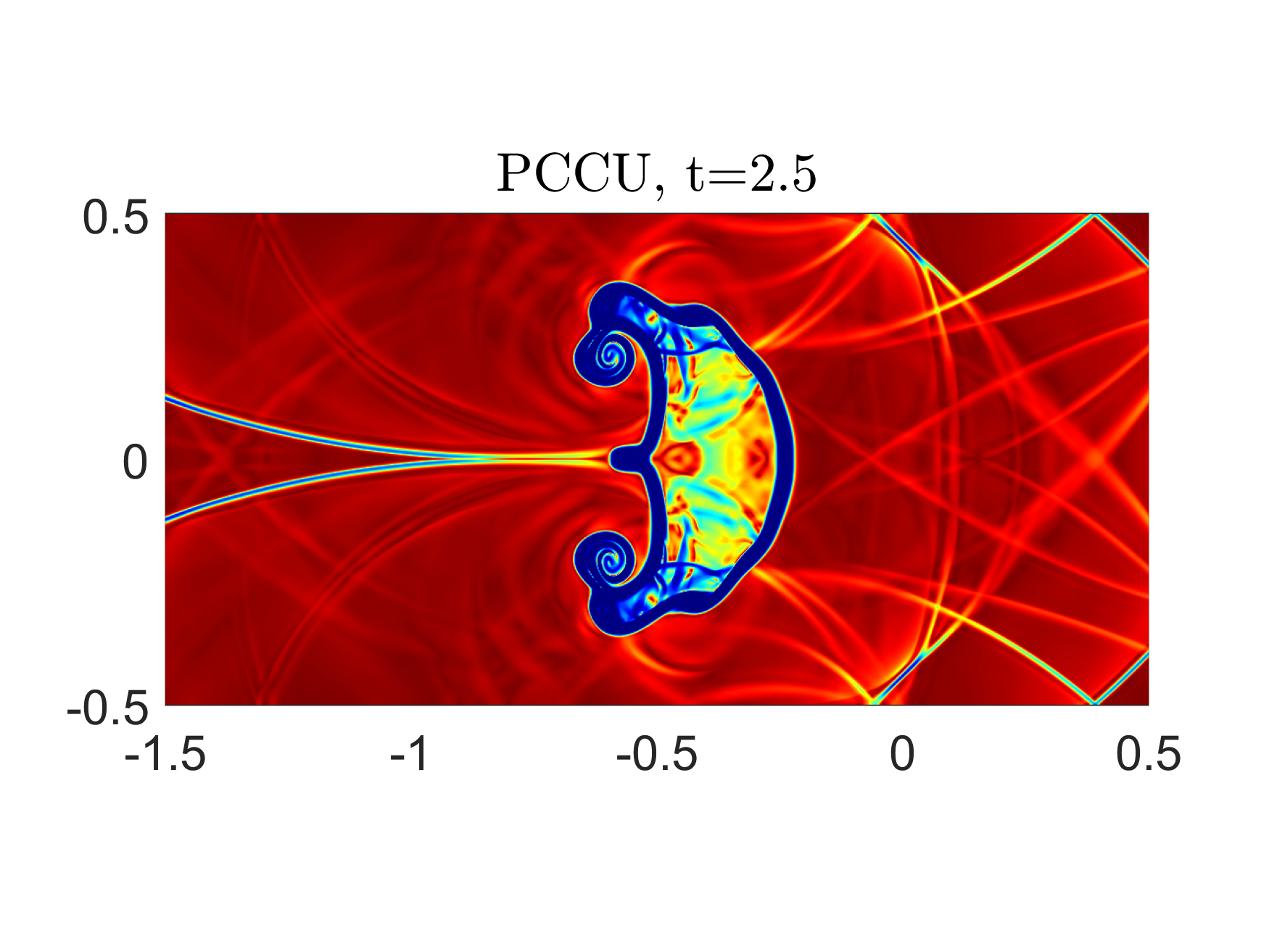}\hspace*{1.cm}
            \includegraphics[trim=0.5cm 1.9cm 1.cm 1.5cm, clip, width=6.cm]{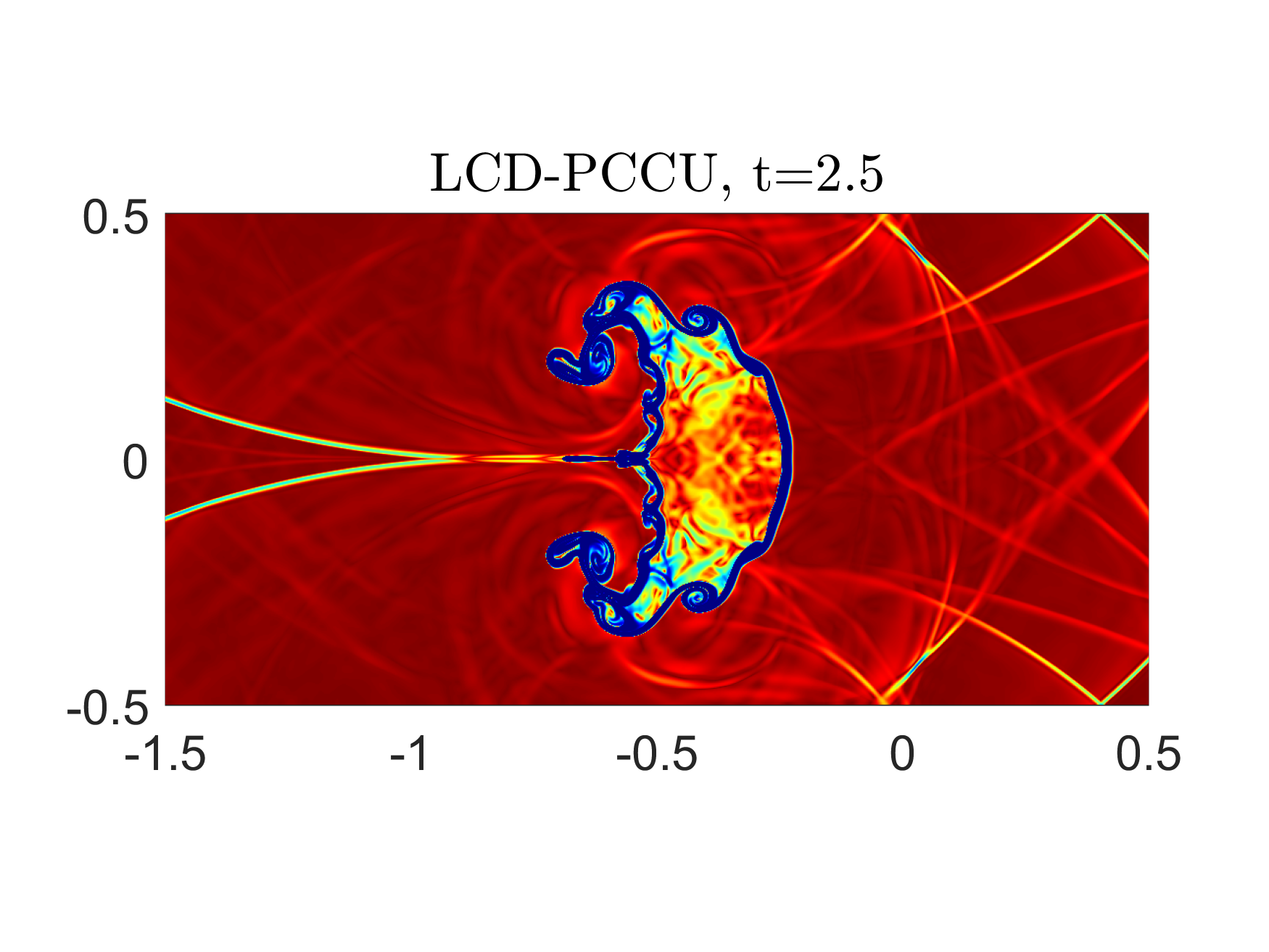}}
\centerline{\includegraphics[trim=0.5cm 1.9cm 1.cm 1.5cm, clip, width=6.cm]{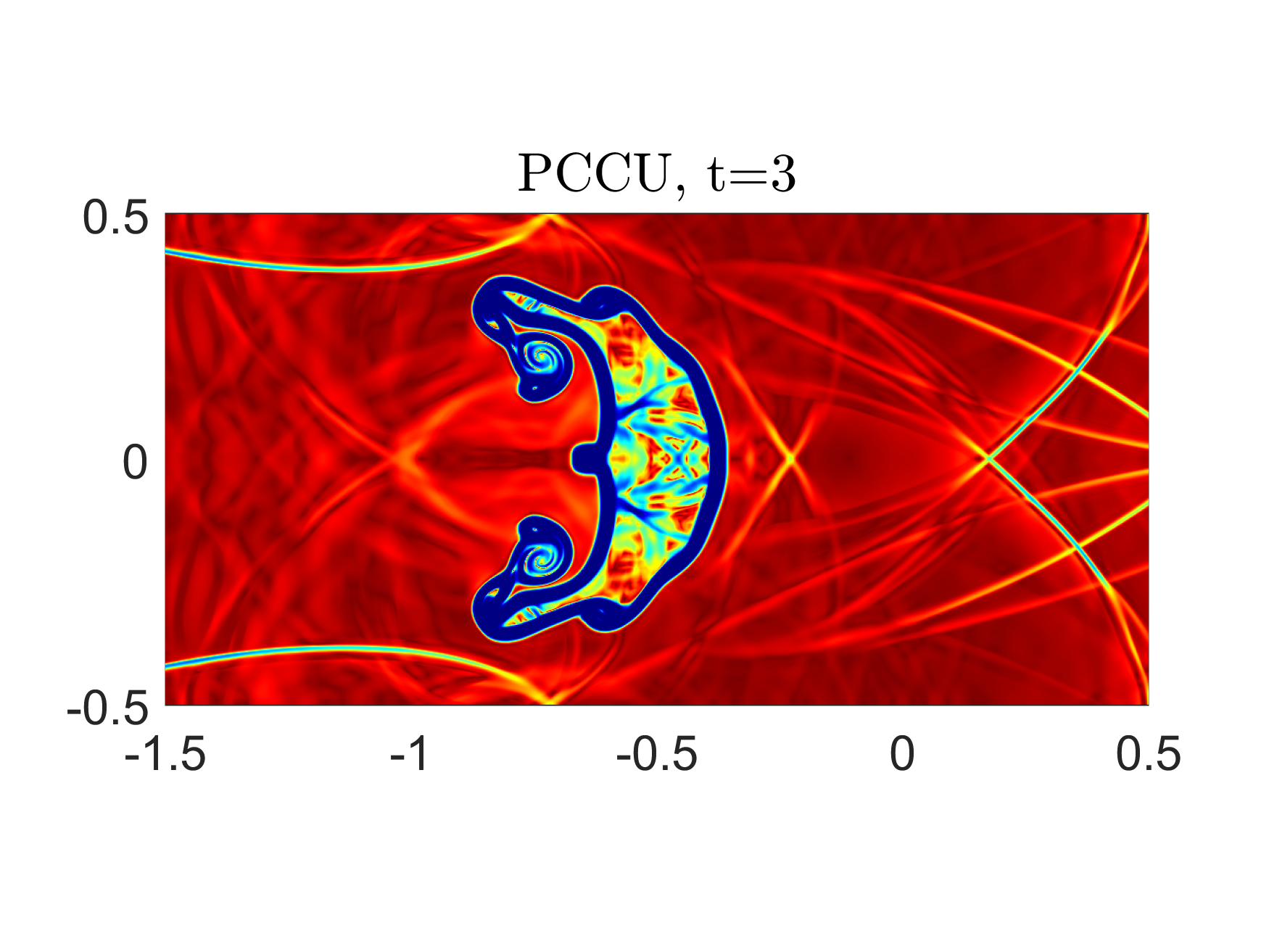}\hspace*{1.cm}
            \includegraphics[trim=0.5cm 1.9cm 1.cm 1.5cm, clip, width=6.cm]{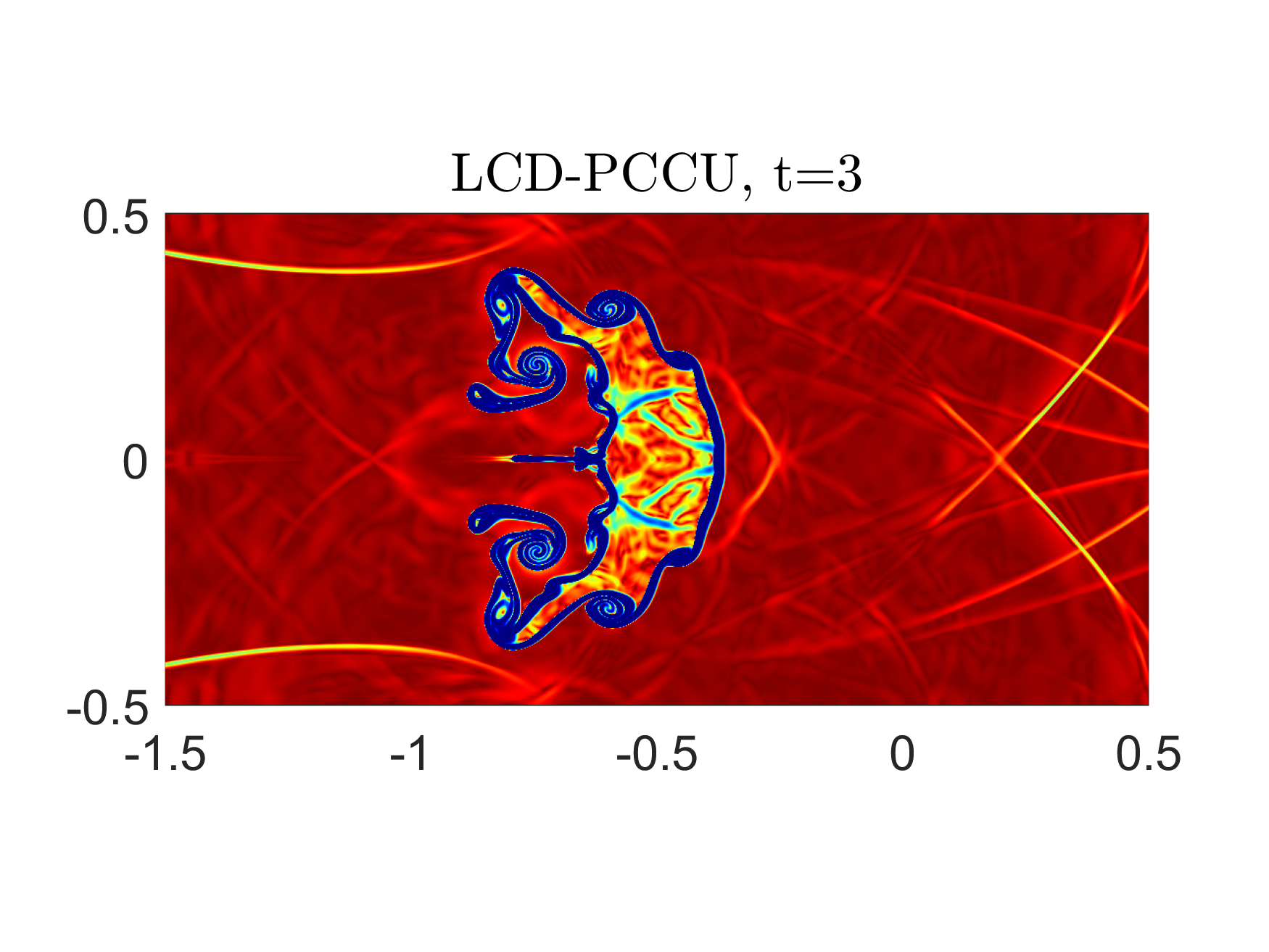}}
\caption{\sf Example 4: Shock-R22 bubble interaction by the PCCU (left column) and LCD-PCCU (right column) schemes at different times.
\label{fig45a}}
\end{figure}

We compute the numerical solutions until the final time $t=3$ on a uniform mesh with $\dx=\dy=1/500$ and present different stages of the
shock-bubble interaction computed by the PCCU and LCD-PCCU schemes in Figure \ref{fig45a}. One can clearly see that upon encountering the
bubble, the shock wave undergoes both refraction and reflection. Under the force of the shock, the bubble undergoes compression and motion.
These findings are in a good agreement with the computational results reported in \cite{Quirk1996,chertock7}. Figure \ref{fig45a}
demonstrates that the LCD-PCCU scheme captures material interfaces more accurately (compared to the PCCU scheme) and also enhances the
resolution of fine details of the computed solution. Once again, this indicates that the LCD-PCCU scheme is less diffusive than its PCCU
counterpart.

\subsubsection*{Example 5---A Cylindrical Explosion Problem}
In the last 2-D example taken from \cite{CKX_23} (see also \cite{XL17}), we consider the case where a cylindrical explosive source is
located between an air-water interface and an impermeable wall. The initial conditions are given by
\begin{equation*}
(\rho,u,v,p;\gamma,\pi_\infty)=\begin{cases}
(1.27,0,0,8290,2,0),&\mbox{$(x-5)^2+(y-2)^2<1$},\\
(0.02,0,0,1;1.4,0),&\mbox{$y>4$},\\
(1,0,0,1;7.15,3309),&\mbox{otherwise},\end{cases}
\end{equation*}
with the solid wall boundary conditions imposed at the bottom, and the free boundary conditions set on the other sides of the computational
domain $[0,10]\times[0,6]$.

\begin{figure}[ht!]
\centerline{\includegraphics[trim=1.2cm 1.5cm 1.1cm 1cm, clip, width=6.cm]{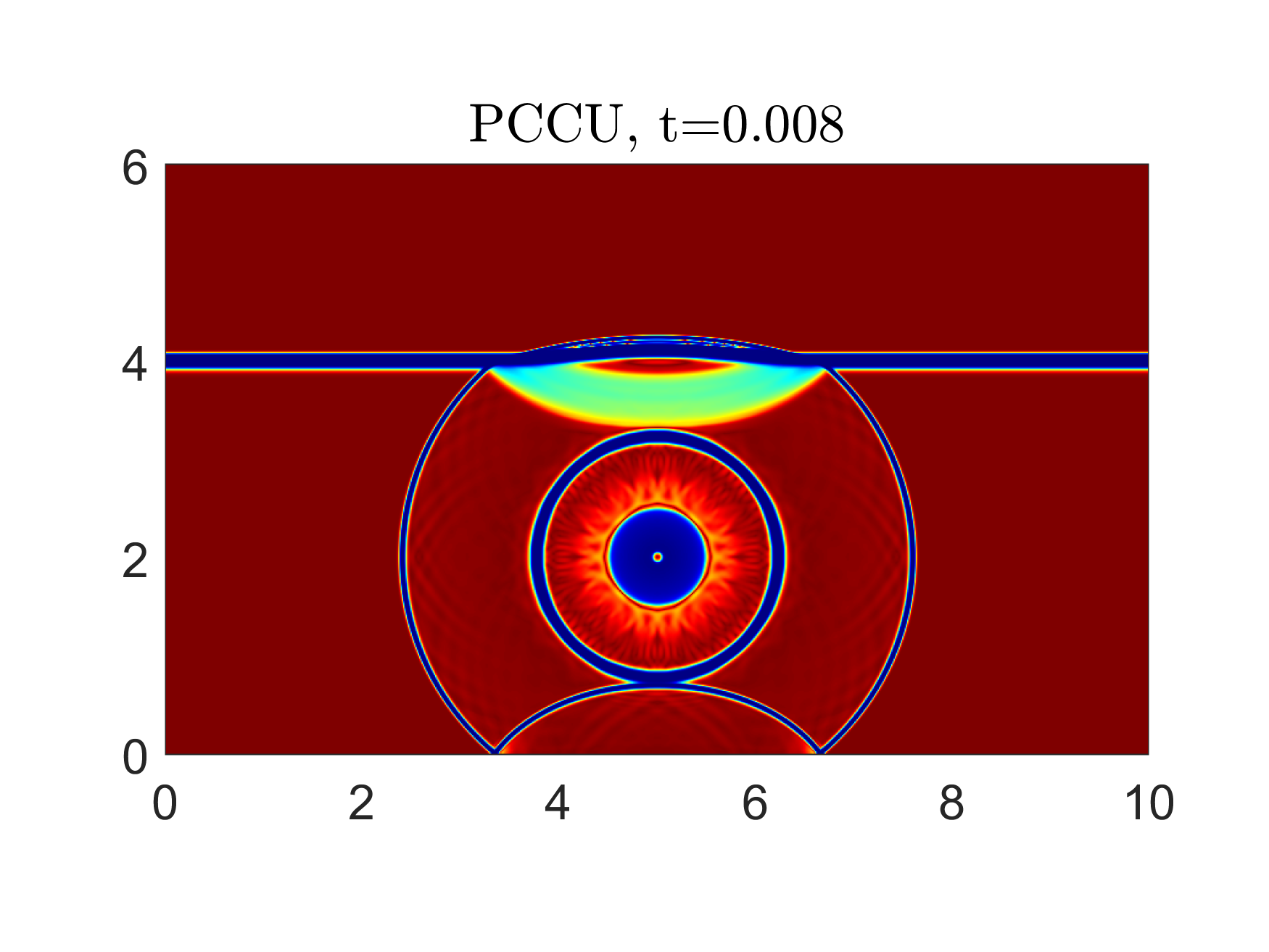}\hspace*{1cm}
            \includegraphics[trim=1.2cm 1.5cm 1.1cm 1cm, clip, width=6.cm]{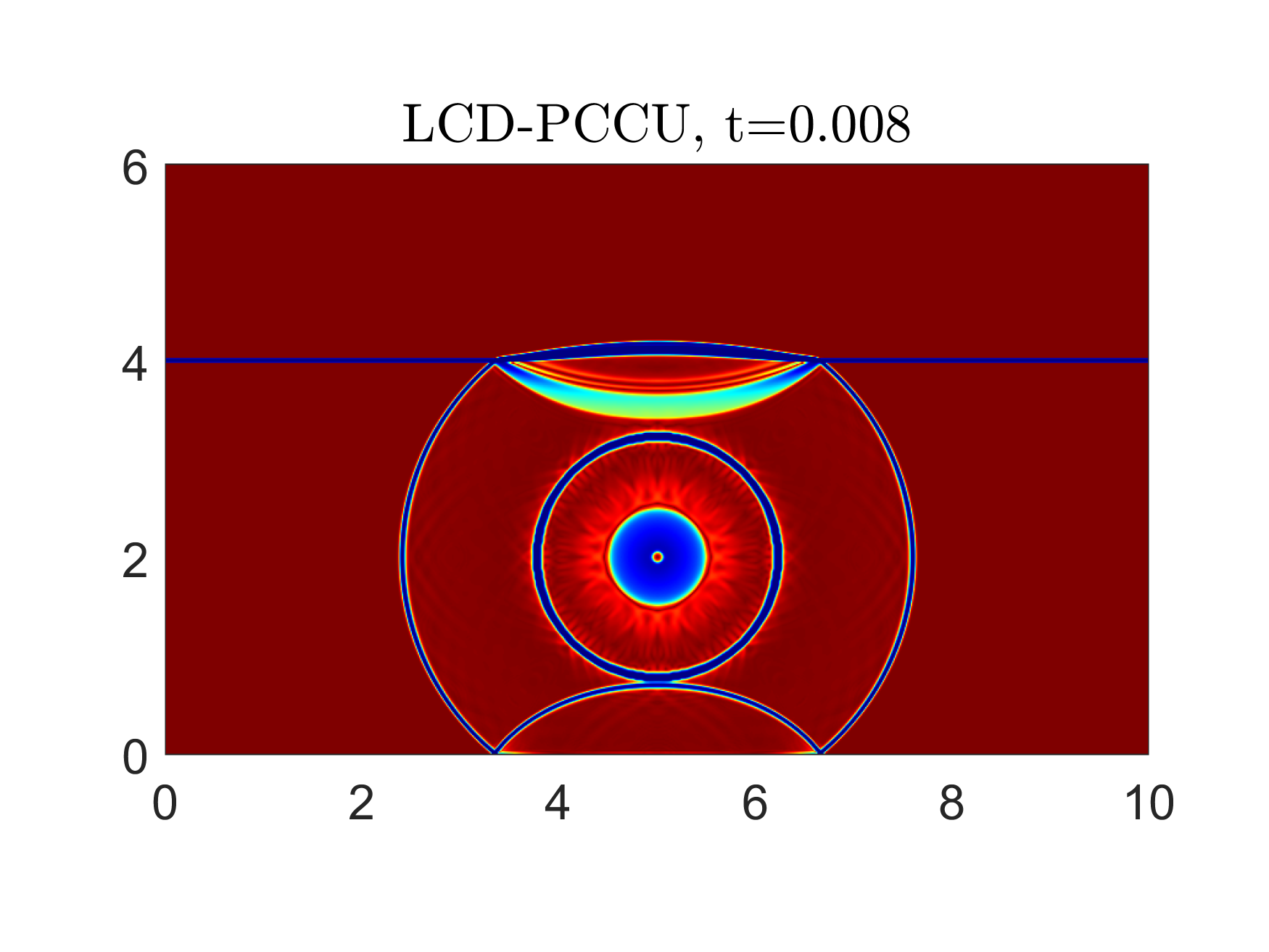}}
\vskip3pt
\centerline{\includegraphics[trim=1.2cm 1.5cm 1.1cm 1cm, clip, width=6.cm]{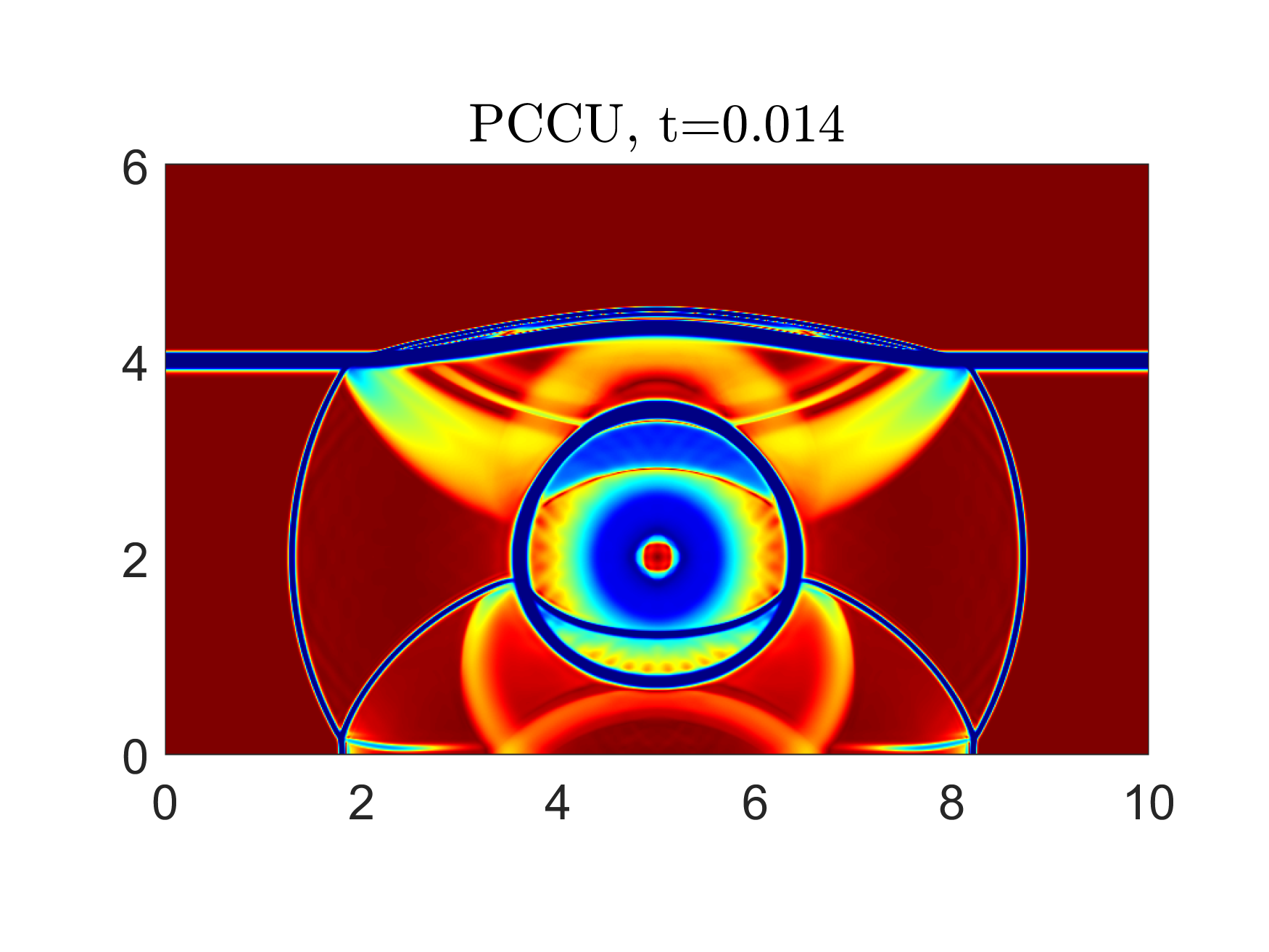}\hspace*{1cm}
            \includegraphics[trim=1.2cm 1.5cm 1.1cm 1cm, clip, width=6.cm]{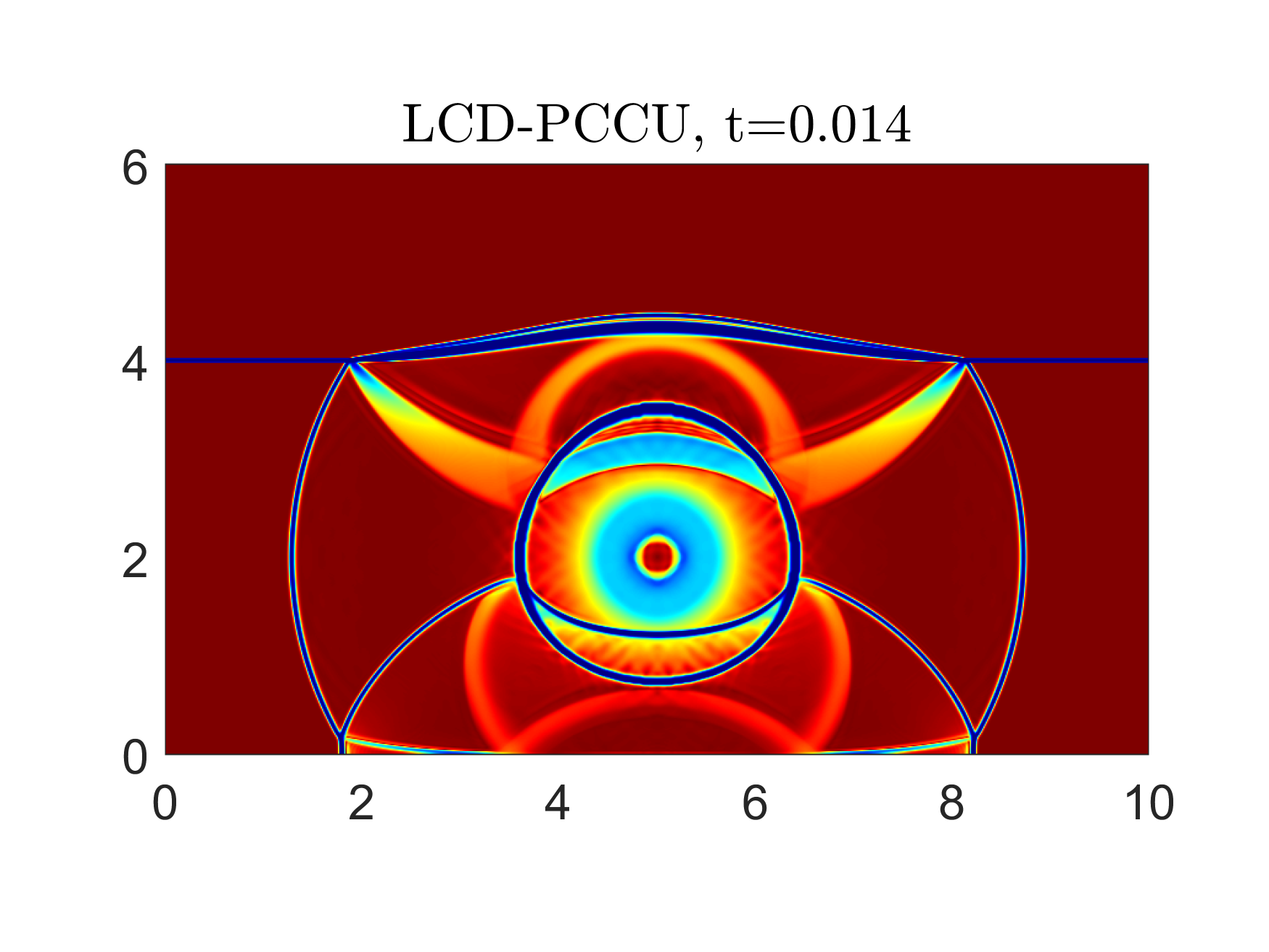}}
\vskip3pt
\centerline{\includegraphics[trim=1.2cm 1.5cm 1.1cm 1cm, clip, width=6.cm]{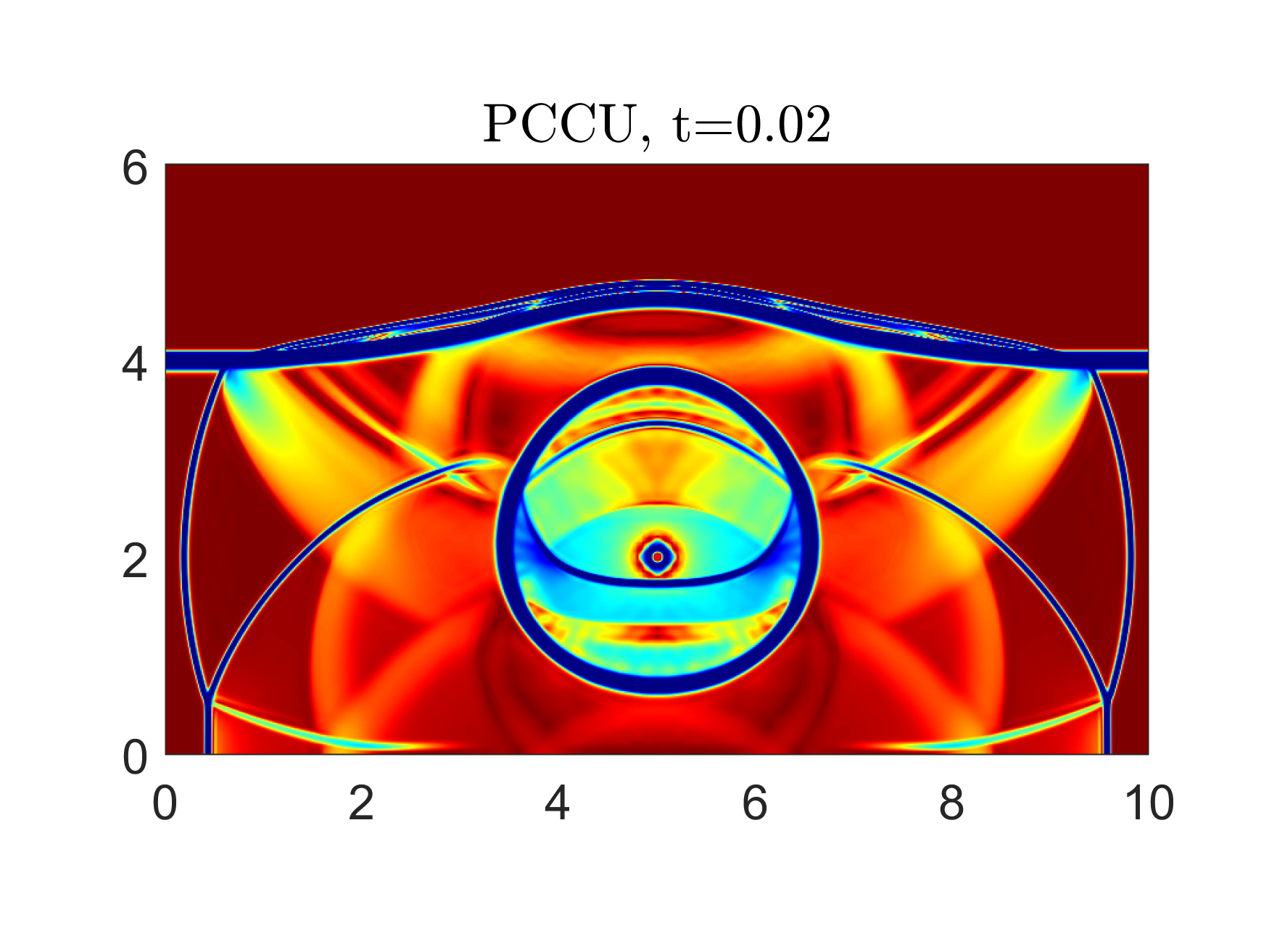}\hspace*{1cm}
            \includegraphics[trim=1.2cm 1.5cm 1.1cm 1cm, clip, width=6.cm]{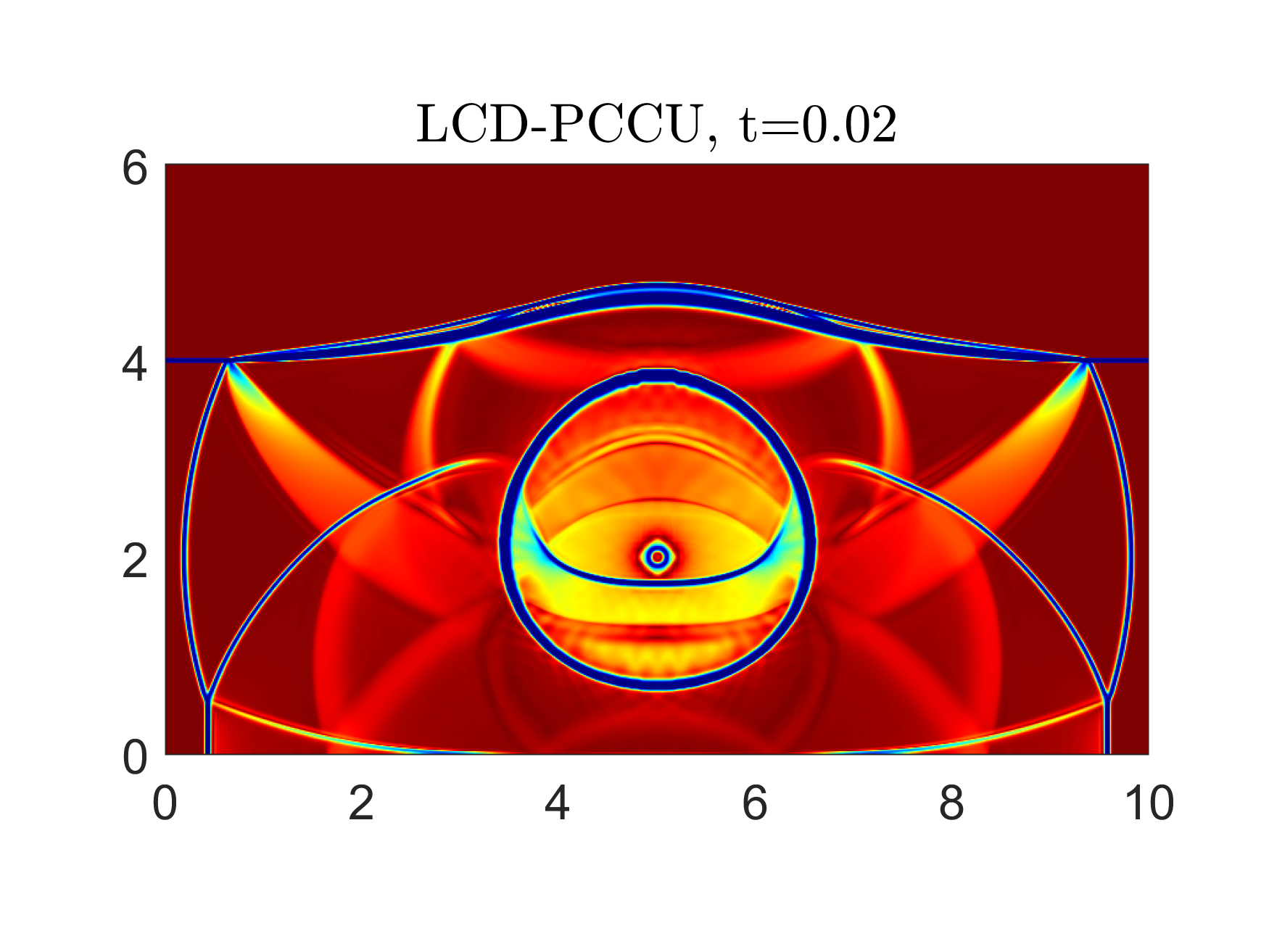}}
\caption{\sf Example 5: Time evolution of the solution computed by the PCCU (left column) and LCD-PCCU (right column) schemes.
\label{fig46a}}
\end{figure}

We compute the numerical solutions until the final time $t=0.02$ on a uniform mesh with $\dx=\dy=1/80$ by the studied PCCU and LCD-PCCU
schemes. In Figure \ref{fig46a}, we present time snapshots of the obtained results, which qualitatively look similar to the numerical
results reported in \cite{CKX_23,XL17}. As one can see, the LCD-PCCU scheme captures both the material interfaces and many of the developed
wave structures substantially sharper than the PCCU scheme.

\begin{figure}[ht!]
\centerline{\includegraphics[trim=1.2cm 1.5cm 1.1cm 1cm, clip, width=6.cm]{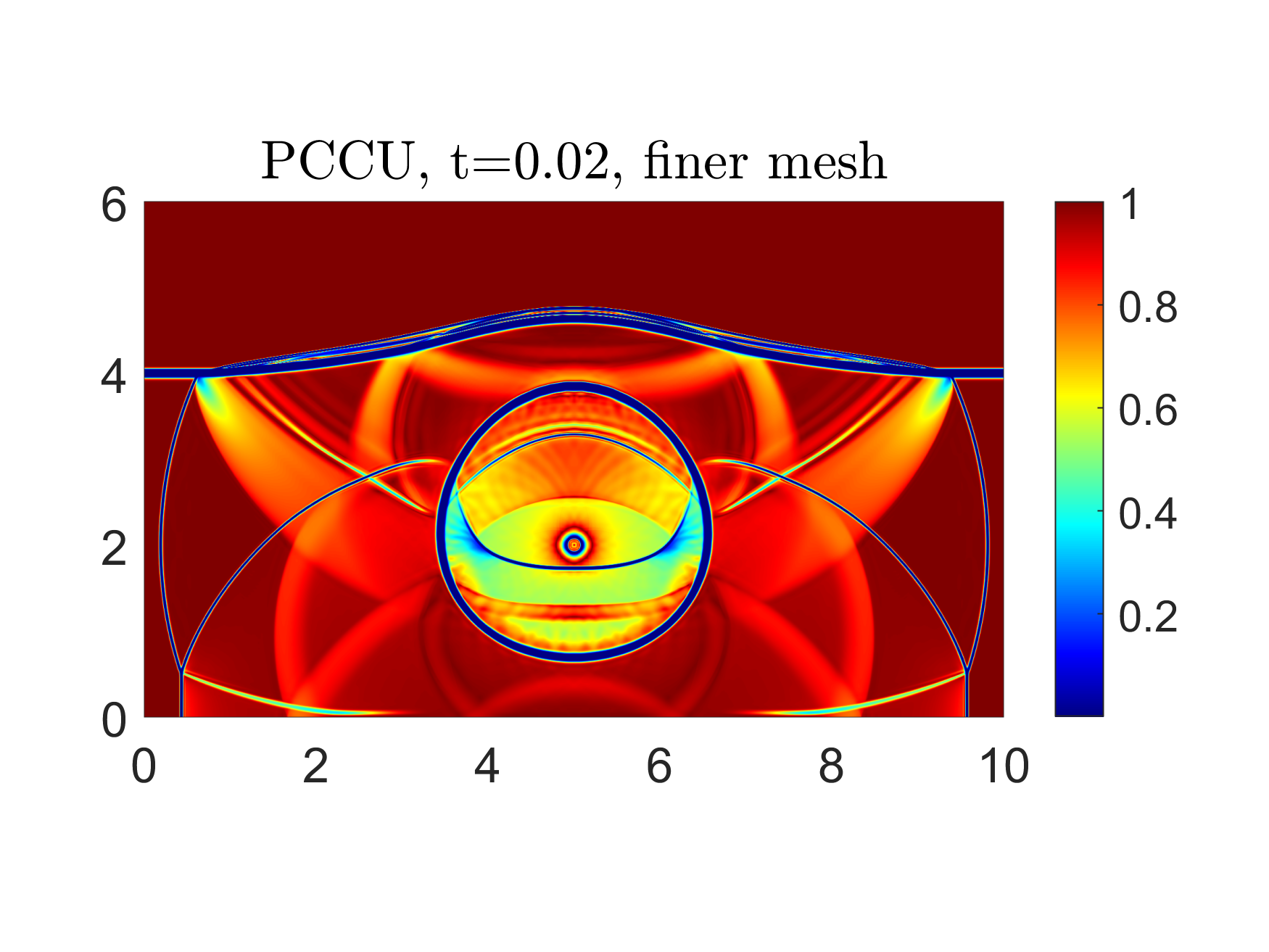}}
\caption{\sf Example 5: A finer mesh ($\dx=\dy=1/160$) PCCU solution at $t=0.02$.\label{fig46b}}
\end{figure}

One can also notice a substantial difference between the PCCU and LCD-PCCU solutions at the final time. In order to clarify this point, we
refine the mesh and compute the numerical solution on a finer mesh with $\dx=\dy=1/160$ by the PCCU scheme and present the obtained result
in Figure \ref{fig46b}, where one can clearly see that the PCCU solution computed on a finer mesh is close to the LCD-PCCU solution computed
with $\dx=\dy=1/80$. Once again, this indicates that the LCD-PCCU solution is more accurate than the PCCU one.

\subsection{Thermal Rotating Shallow Water Equations}
In this section, we apply the developed LCD-PCCU schemes to the 1-D and 2-D TRSW equations \cite{KLZ_2D}, which in the 2-D case can be
written as the system \eref{1.2} with 
\begin{equation}
\resizebox{0.92\hsize}{!}
{$\begin{aligned}
&\mU=(h,q,p,hb)^\top,\quad\mF(\mU)=\Big(hu,hu^2+\frac{b}{2}h^2,huv,hub\Big)^\top,\quad
\mG(\mU)=\Big(hv,huv,hv^2+\frac{b}{2}h^2,hvb\Big)^\top,\\
&B(\mU)=C(\mU)\equiv\mo,\quad\bm S^x(\mU)=(0,fhv-hbZ_x,0,0)^\top,\quad\bm S^y(\mU)=(0,0,-fhu-hbZ_y,0)^\top,
\end{aligned}$}
\label{5.21}
\end{equation}
where $h$ denotes the thickness of the fluid layer, $u$ and $v$ are the zonal and meridional velocities, $q:=hu$ and $p:=hv$ are the
corresponding discharges, $b$ is the layer-averaged buoyancy variable, and $Z(x,y)$ stands for a time-independent bottom topography. In the
oceanic context, $b=g\rho/\rho_0$, where $g$ is the acceleration due to gravity, $\rho$ and $\rho_0$ are the variable and constant parts of
the water density, respectively. In the atmospheric context, the densities $\rho$ and $\rho_0$ should be replaced with the potential
temperatures $\theta$ and $\theta_0$. Finally, $f(y)=f_0+\beta y$ is the Coriolis parameter, where $f_0$ and $\beta$ are positive constants.
\begin{rmk}
For the TRSW equations, we use the equilibrium variables introduced in \cite{CKL23} and compute the point values ($\bm U^\pm_\kph$ and 
$\breve{\bm U}^\pm_\kph$ in the 1-D case and $\bm U^\pm_{\jph,k}$, $\bm U^\pm_{j,\kph}$, $\breve{\bm U}^\pm_{\jph,k}$, and
$\breve{\bm U}^\pm_{j,\kph}$ in the 2-D case) as well as the global fluxes ($\bm L^\pm_\kph$ in the 1-D case and $\bm K^\pm_{\jph,k}$ and
$\bm L^\pm_{j,\kph}$ in the 2-D case) precisely as it was done in \cite{CKL23}. We have also used the same numerical diffusion switch
function, which helps to enforce the WB property, as in \cite{CKL23}.
\end{rmk}
\begin{rmk}
In Appendix \ref{appc}, we provide the expressions for the matrices $\widehat{\cal A}_\kph$ and $\widehat{\cal A}^{\,\mF}_{\jph,k}$ and the
corresponding matrices $R_\kph$ and $R_{\jph,k}$ for the 1-D and 2-D TRSW equations, respectively. The matrices
$\widehat{\cal A}^{\,\mG}_{j,\kph}$ and $R_{j,\kph}$ can be obtained in an analogous way.
\end{rmk}

\subsubsection{One-Dimensional Examples}
In this section, we consider two numerical examples for the so-called 1-D Ripa system, which is the 1-D version of \eref{1.2}, \eref{5.21},
but with $f(y)\equiv0$:
\begin{equation*}
\begin{cases}
h_t+(hv)_y=0,\\[0.3ex]
p_t+\Big(hv^2+\dfrac{b}{2}h^2\Big)_y=-hbZ_y,\\[0.5ex]
(hb)_t+(hvb)_y=0.
\end{cases}
\end{equation*}

\subsubsection*{Example 6---Constant Pressure Equilibrium}
In the first example taken from \cite{CKL23}, we consider the following initial data:
\begin{equation}
(h,v,b)(y,0)=\begin{cases}(2,0,1)&\mbox{if}~y<0,\\(1,0,4)&\mbox{otherwise},\end{cases}
\label{7.2}
\end{equation}
which are prescribed in the computational domain $[-5,5]$ together with the flat bottom topography $Z(y)\equiv0$ subject to the free
boundary conditions.

\begin{figure}[ht!]
\centerline{\includegraphics[trim=1.0cm 0.4cm 1.2cm 0.3cm, clip, width=5.cm]{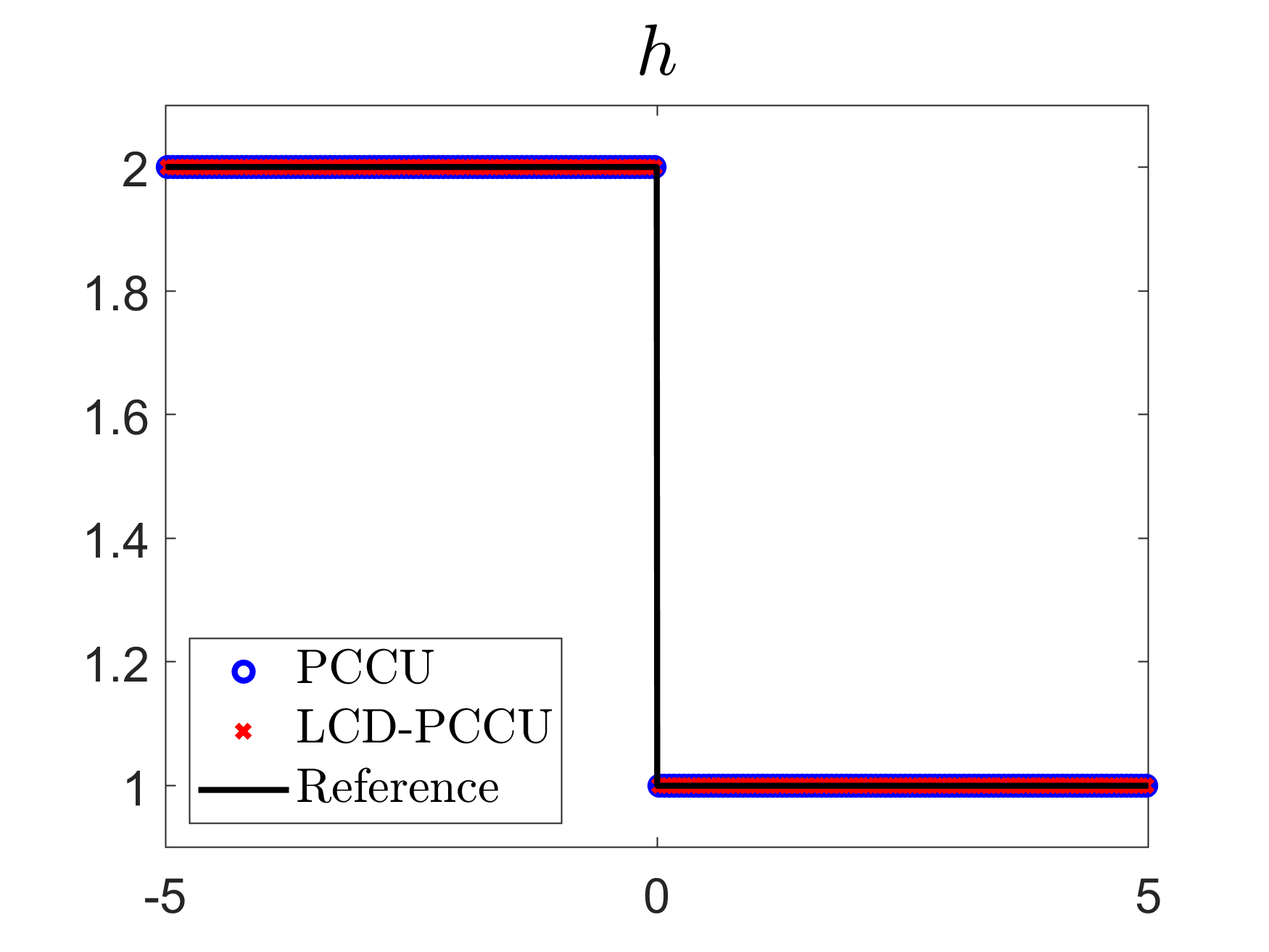}\hspace*{1.0cm}
            \includegraphics[trim=1.0cm 0.4cm 1.2cm 0.3cm, clip, width=5.cm]{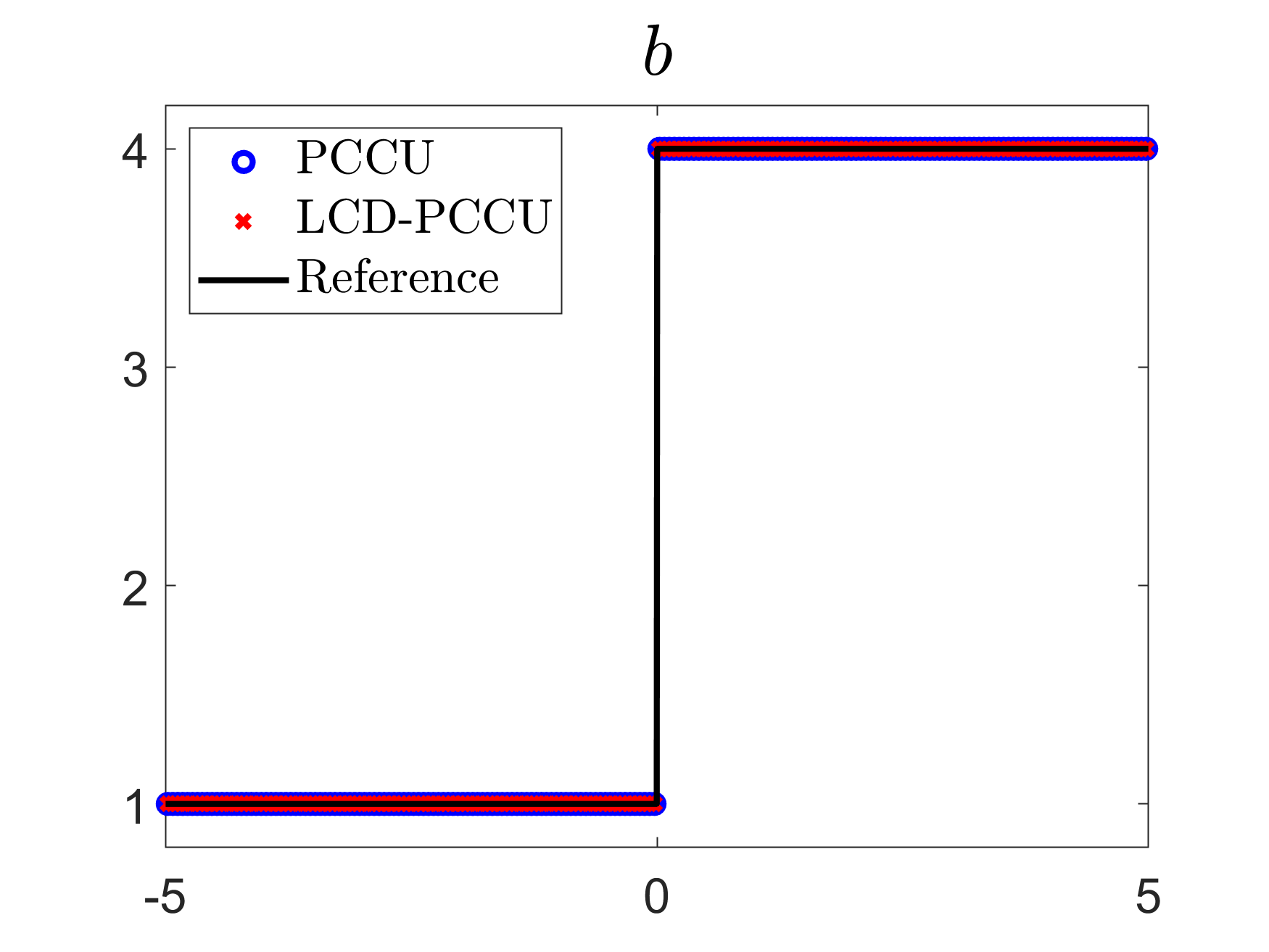}}
\caption{\sf Example 6: $h$ and $b$ computed by the PCCU and LCD-PCCU schemes.\label{fig41}}
\end{figure}

We first compute the numerical solution by the PCCU and LCD-PCCU schemes until the final time $t=10$ on the uniform mesh with $\dx=1/20$ and
plot the obtained fluid thickness and buoyancy variables in Figure \ref{fig41}, where the reference solution is computed by the PCCU scheme
on a much finer mesh with $\dx=1/400$. One can clearly see that both the PCCU and LCD-PCCU schemes can exactly preserve the steady state.

\begin{figure}[ht!]
\centerline{\includegraphics[trim=1.cm 0.4cm 1.0cm 0.3cm, clip, width=5cm]{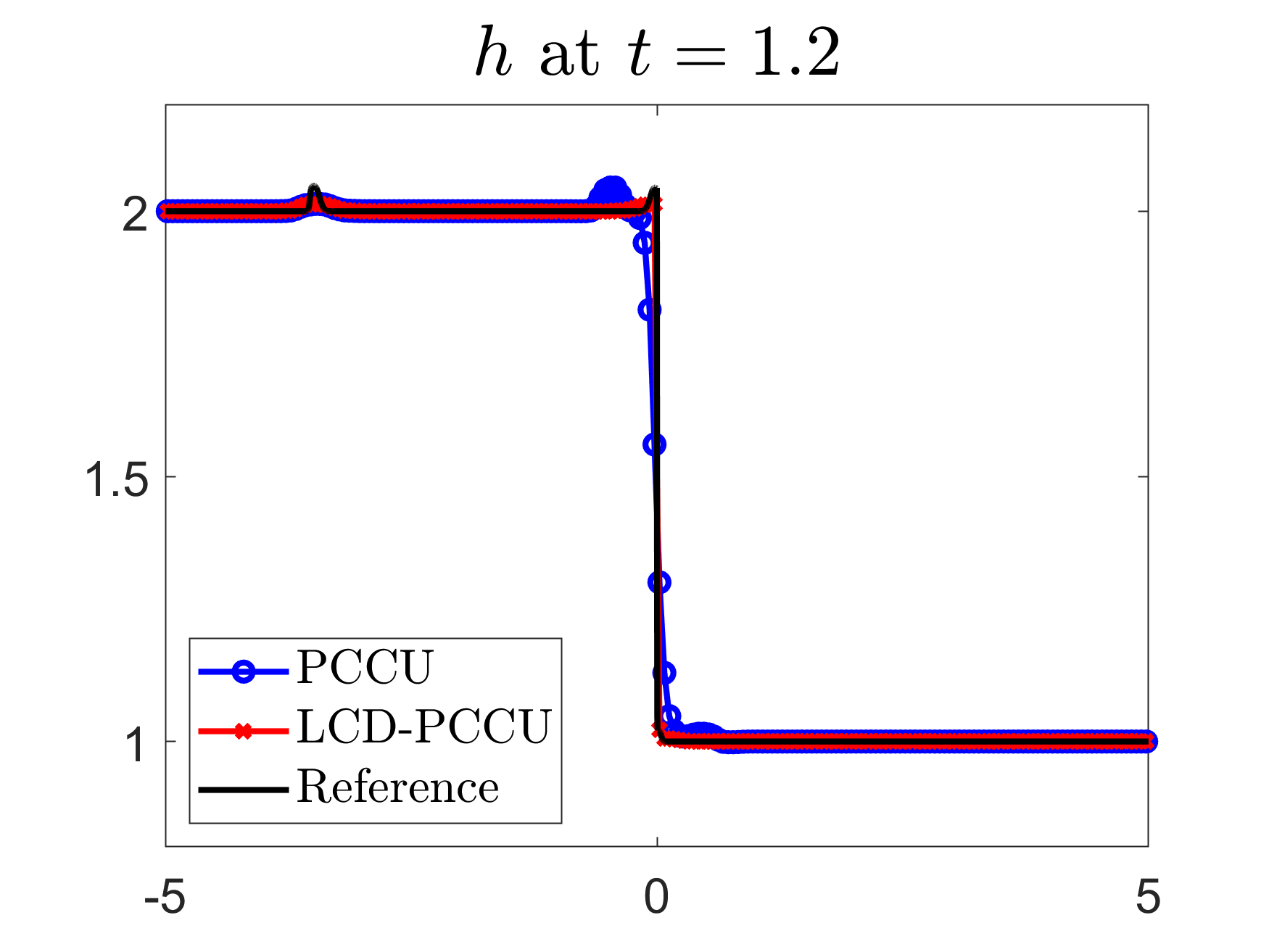}\hspace*{1.0cm}
            \includegraphics[trim=1.cm 0.4cm 1.0cm 0.3cm, clip, width=5cm]{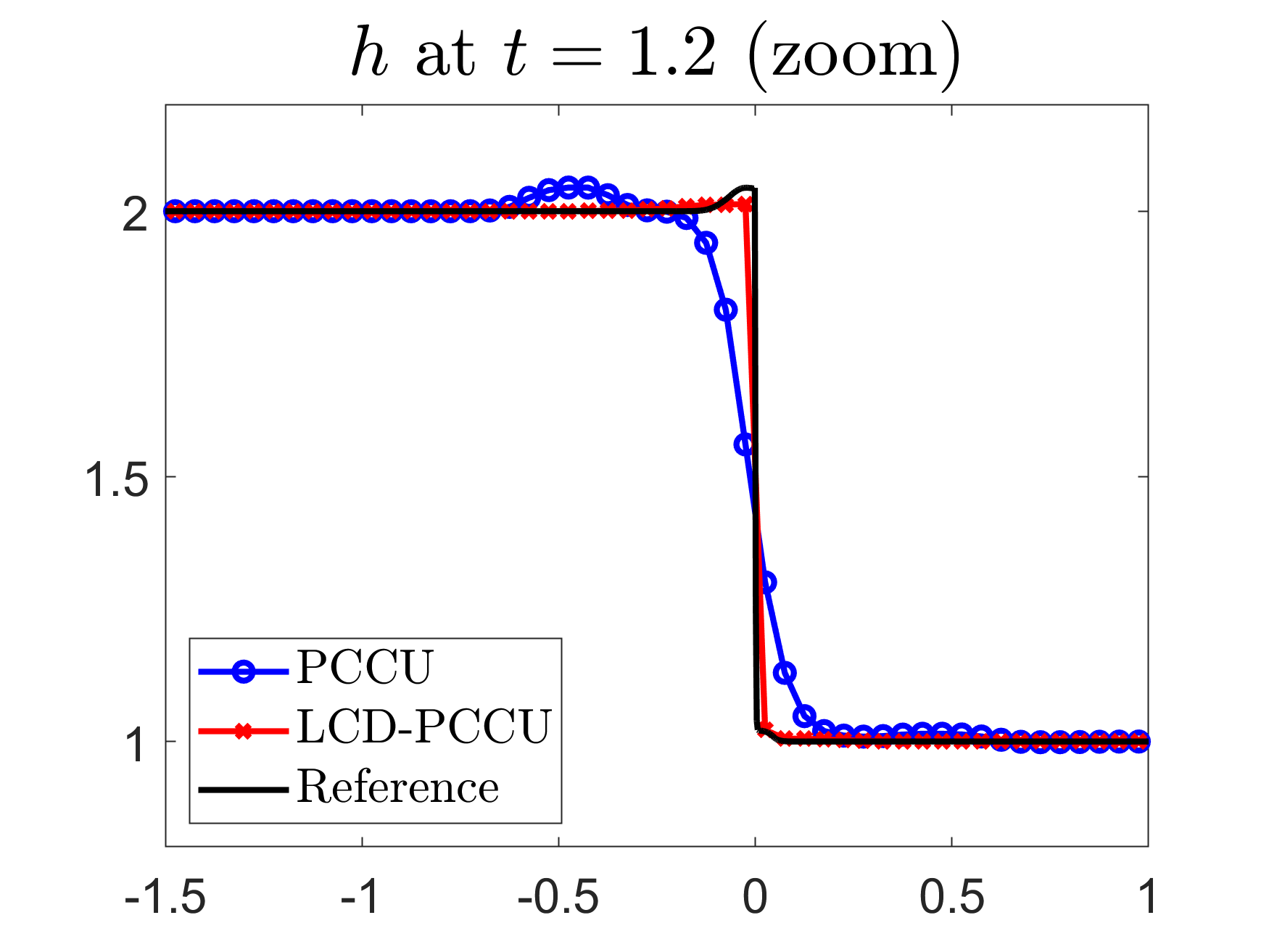}}
\vskip5pt
\centerline{\includegraphics[trim=1.cm 0.4cm 1.0cm 0.3cm, clip, width=5cm]{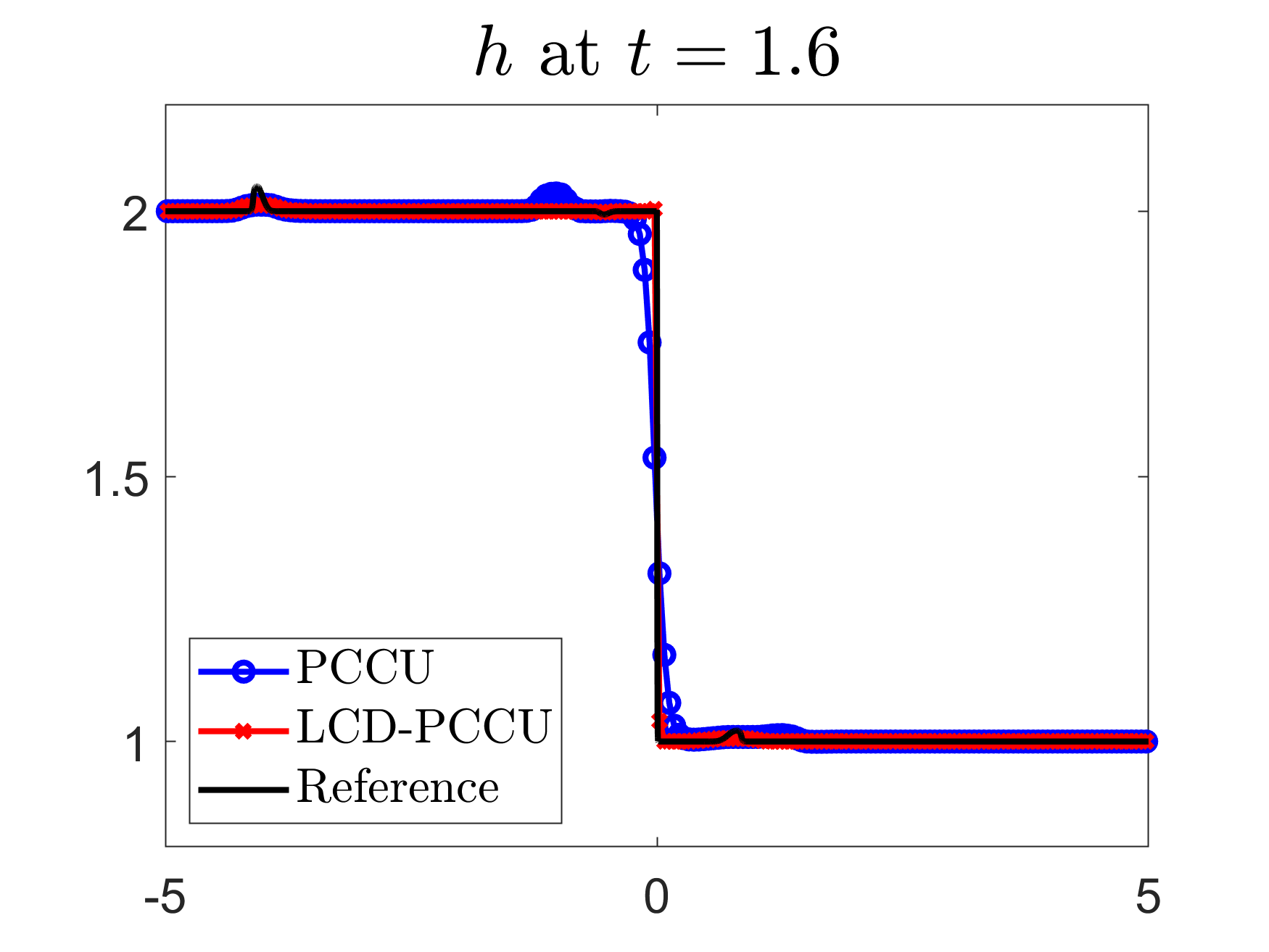}\hspace*{1.0cm}
            \includegraphics[trim=1.cm 0.4cm 1.0cm 0.3cm, clip, width=5cm]{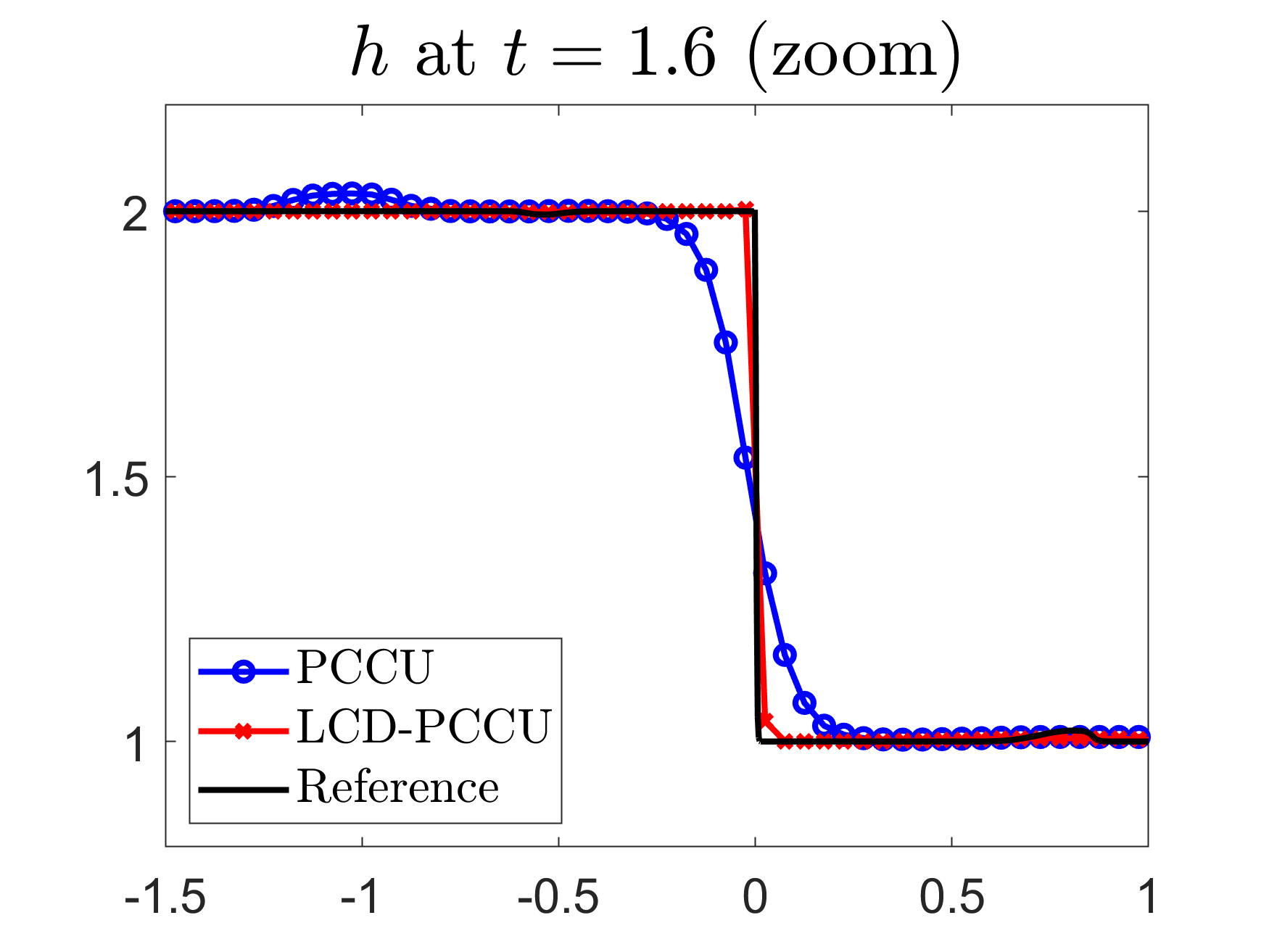}}
\caption{\sf Example 6: $h$ computed by the PCCU and LCD-PCCU schemes at $t=1.2$ (top row) and 1.6 (bottom row) and zoom at $x\in[-1.5,1]$
(right column).\label{fig455a}}
\end{figure}

We then test the ability of the proposed schemes to capture small perturbations of this steady state. To this end, we first denote by
$h_{\rm eq}(y)=h(y,0)$ given by \eref{7.2}, and then modify the initial condition of $h$ in the following manner:
\begin{equation*}
h(y,0)=h_{\rm eq}(y)+\begin{cases}0.1&\mbox{if}~y\in[-1.8,-1.7],\\0&\mbox{otherwise}.\end{cases}
\end{equation*}
We compute the solution until the final time $t=1.6$ on a uniform mesh with $\dx=1/20$ and plot the obtained fluid thickness at different
times $t=1.2$ and 1.6 in Figure \ref{fig455a} together with the reference solution computed by the LCD-PCCU scheme on a much finer mesh with
$\dx=1/400$. One can clearly see that the LCD-PCCU scheme achieves sharper resolution than the PCCU scheme. Moreover, the LCD-PCCU solution
is substantially less oscillatory than the PCCU one.

\subsubsection*{Example 7---Dam Break Over a Nonflat Bottom}
In the second example taken from \cite{CKL14}, we consider the dam break problem subject to the following initial data:
\begin{equation*}
(h,v,b)(y,0)=\begin{cases}(5-Z(y),0,1)&\mbox{if}~y<0,\\(2-Z(y),0,5)&\mbox{otherwise},\end{cases}
\end{equation*}
and the nonflat bottom topography
\begin{equation*}
Z(y)=\begin{cases}
\cos\big(10\pi(x+0.3)\big)+1&\mbox{if}~-0.4\le y\le-0.2,\\
\hf\big[\cos\big(10\pi(x-0.3)\big)+1\big]&0.2\le y\le0.4,\\
0&\mbox{otherwise}.
\end{cases}
\end{equation*}
The initial conditions are prescribed in the computational domain $[-1,1]$ subject to the free boundary conditions.

\begin{figure}[ht!]
\centerline{\includegraphics[trim=1.3cm 0.4cm 1.1cm 0.2cm, clip, width=5cm]{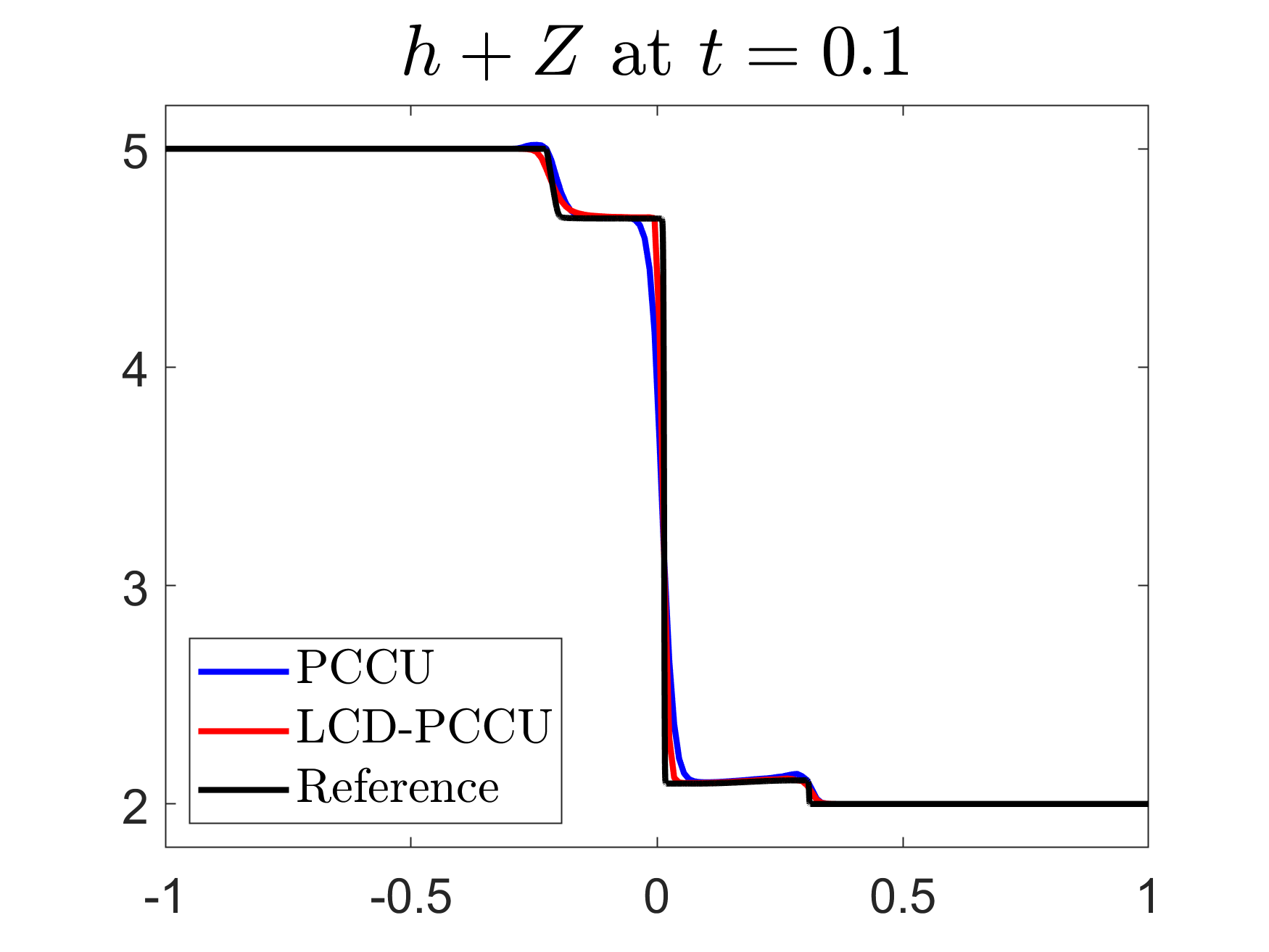}\hspace*{1.0cm}
            \includegraphics[trim=1.3cm 0.4cm 1.1cm 0.2cm, clip, width=5cm]{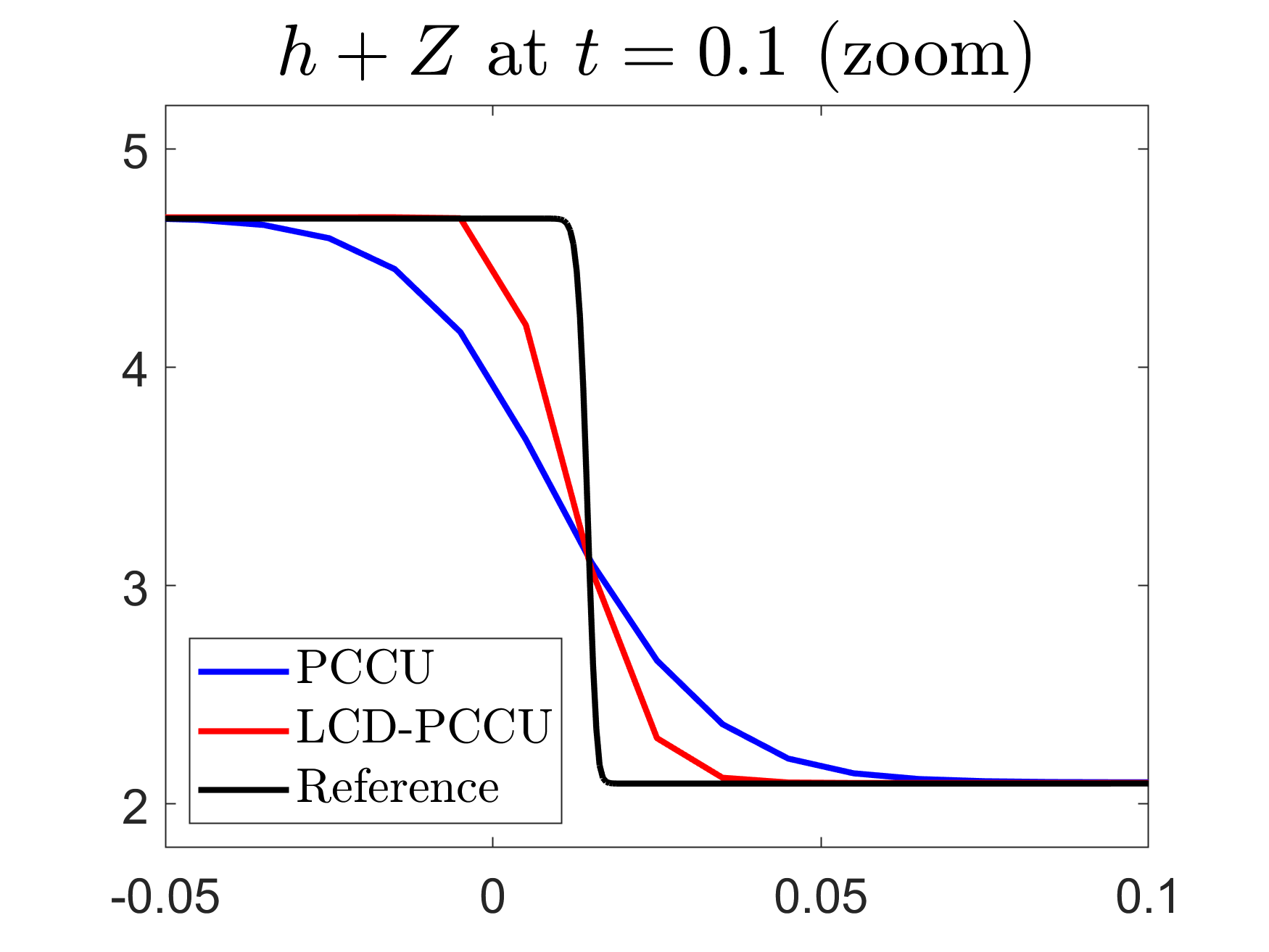}}
\vskip5pt
\centerline{\includegraphics[trim=1.3cm 0.4cm 1.1cm 0.2cm, clip, width=5cm]{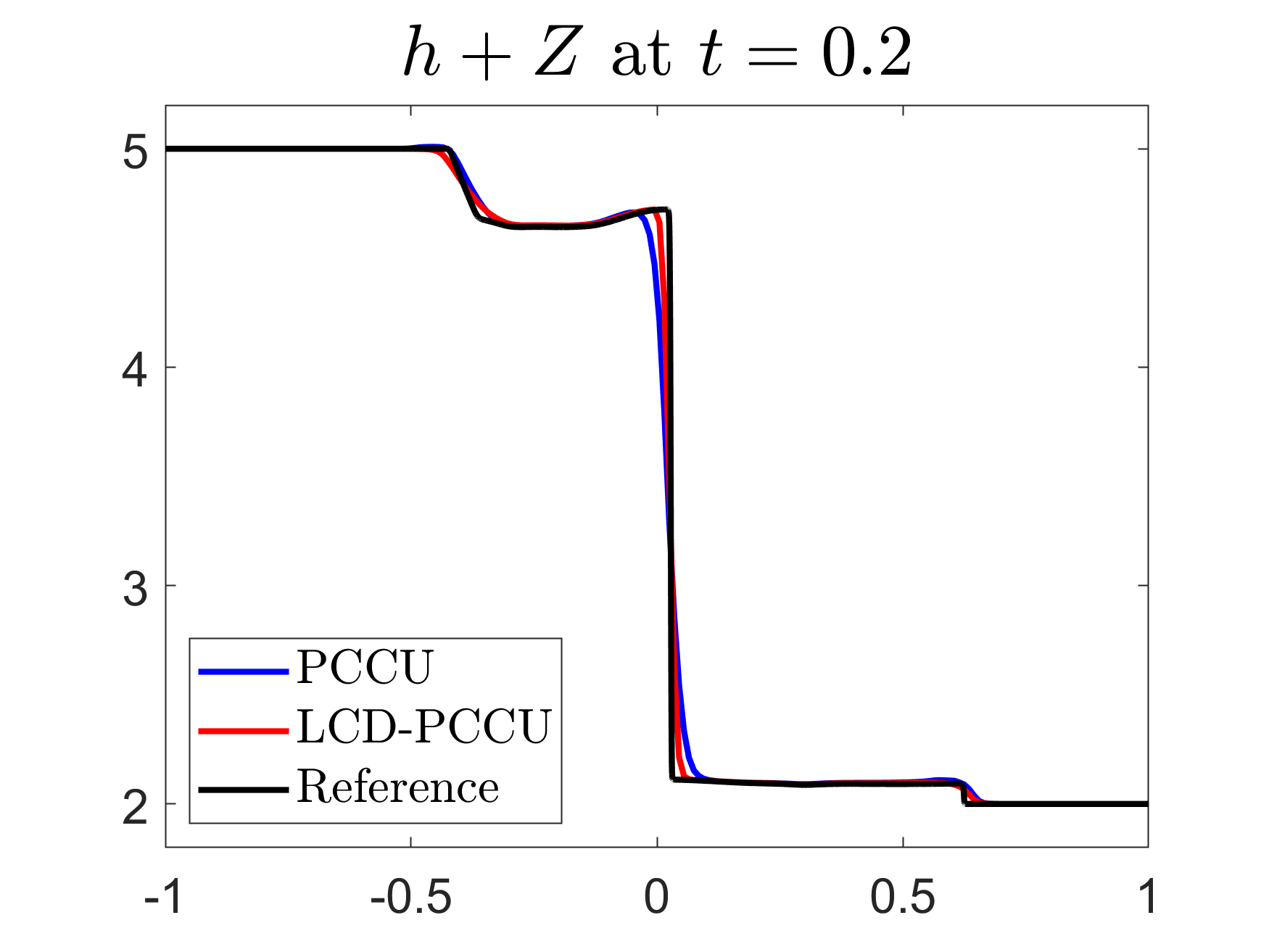}\hspace*{1.0cm}
            \includegraphics[trim=1.3cm 0.4cm 1.1cm 0.2cm, clip, width=5cm]{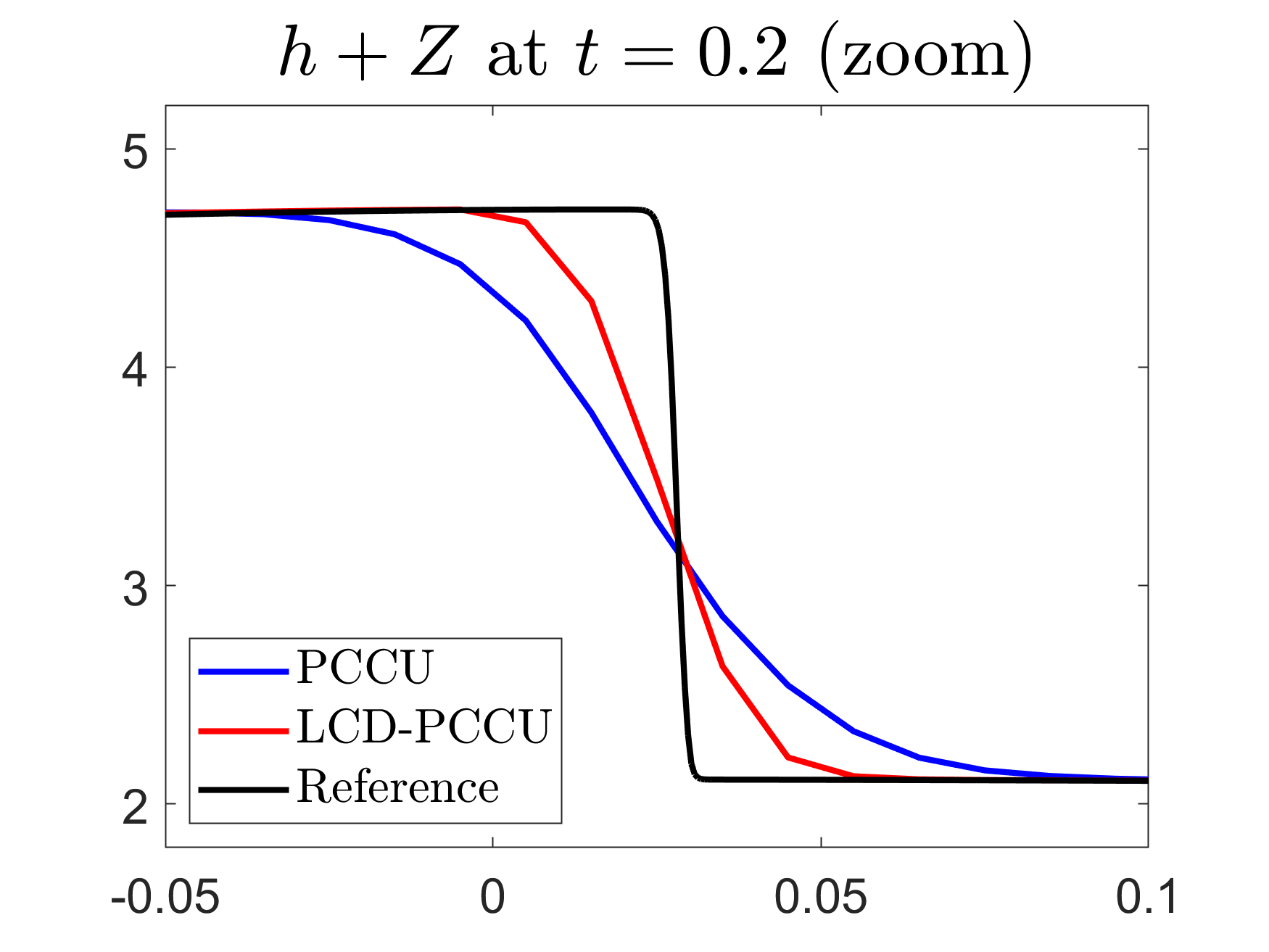}}
\caption{\sf Example 7: $h+Z$ computed by the PCCU and LCD-PCCU schemes at $t=0.1$ (top row) and 0.2 (bottom row) and zoom at
$x\in[-0.05,0.1]$ (right column).\label{fig61}}
\end{figure}

We compute the numerical solution by the PCCU and LCD-PCCU schemes until the final time $t=0.2$ on a uniform mesh with $\dx=1/100$ and plot
the obtained water surface $h+Z$ in Figure \ref{fig61}, where the reference solution is computed by the PCCU scheme on a much finer mesh
with $\dx=1/2000$. One can clearly see that the LCD-PCCU scheme produces higher-resolution results than the PCCU one.

\subsubsection{Two-Dimensional Examples}
In this section, we consider three numerical examples for the 2-D TRSW equations.

\subsubsection*{Example 8---Small Perturbations of a 2-D Steady State}
In this example taken from \cite{CKL14}, we take $f(y)\equiv0$ and consider the following initial conditions: 
\begin{equation*}
(h,u,v,b)(x,y,0)=\begin{cases}
\big(3.1-Z(x,y),0,0,\frac{4}{3}\big)&\mbox{if}~0.01<x^2+y^2<0.09,\\
\big(3-Z(x,y),0,0,\frac{4}{3}\big)&\mbox{if}~0.09<x^2+y^2<0.25~{\rm or}~x^2+y^2<0.01,\\
(2-Z(x,y),0,0,3)&\mbox{otherwise},
\end{cases}
\end{equation*}
with the nonflat bottom topography consisting of two Gaussian shaped humps:
\begin{equation*}
Z(x,y)=\begin{cases}
0.5\,e^{-100\left[(x+0.5)^2+(y+0.5)^2\right]}&\mbox{if}~x<0,\\
0.6\,e^{-100\left[(x-0.5)^2+(y-0.5)^2\right]}&\mbox{otherwise}.
\end{cases}
\end{equation*}
The initial data are prescribed in the computational domain $[-1,1]\times[-1,1]$ subject to the free boundary conditions.

\begin{figure}[ht!]
\centerline{\includegraphics[trim=1.2cm 0.4cm 1.7cm 0.2cm, clip, width=5cm]{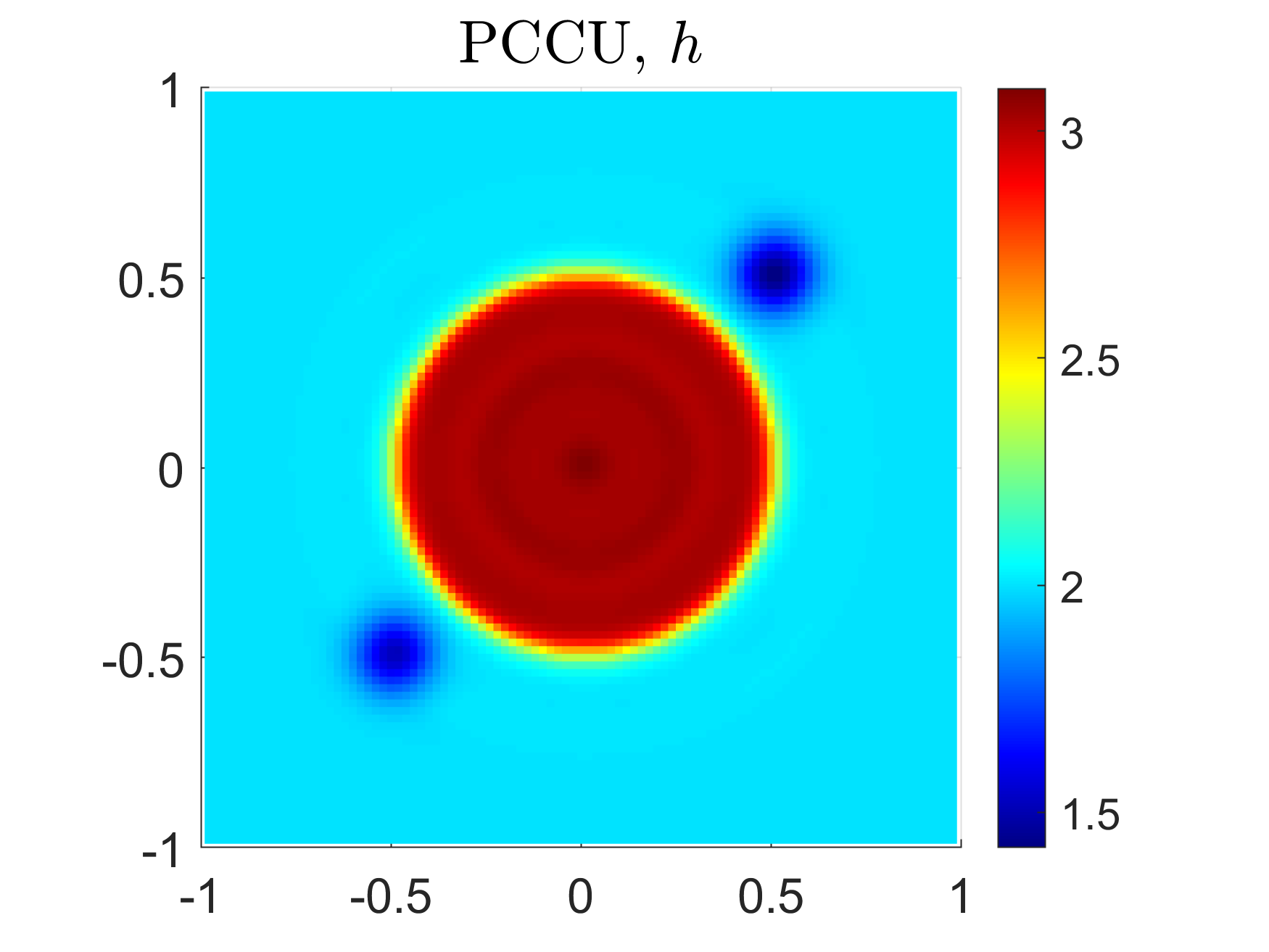}\hspace*{1.0cm}
            \includegraphics[trim=1.2cm 0.4cm 1.7cm 0.2cm, clip, width=5cm]{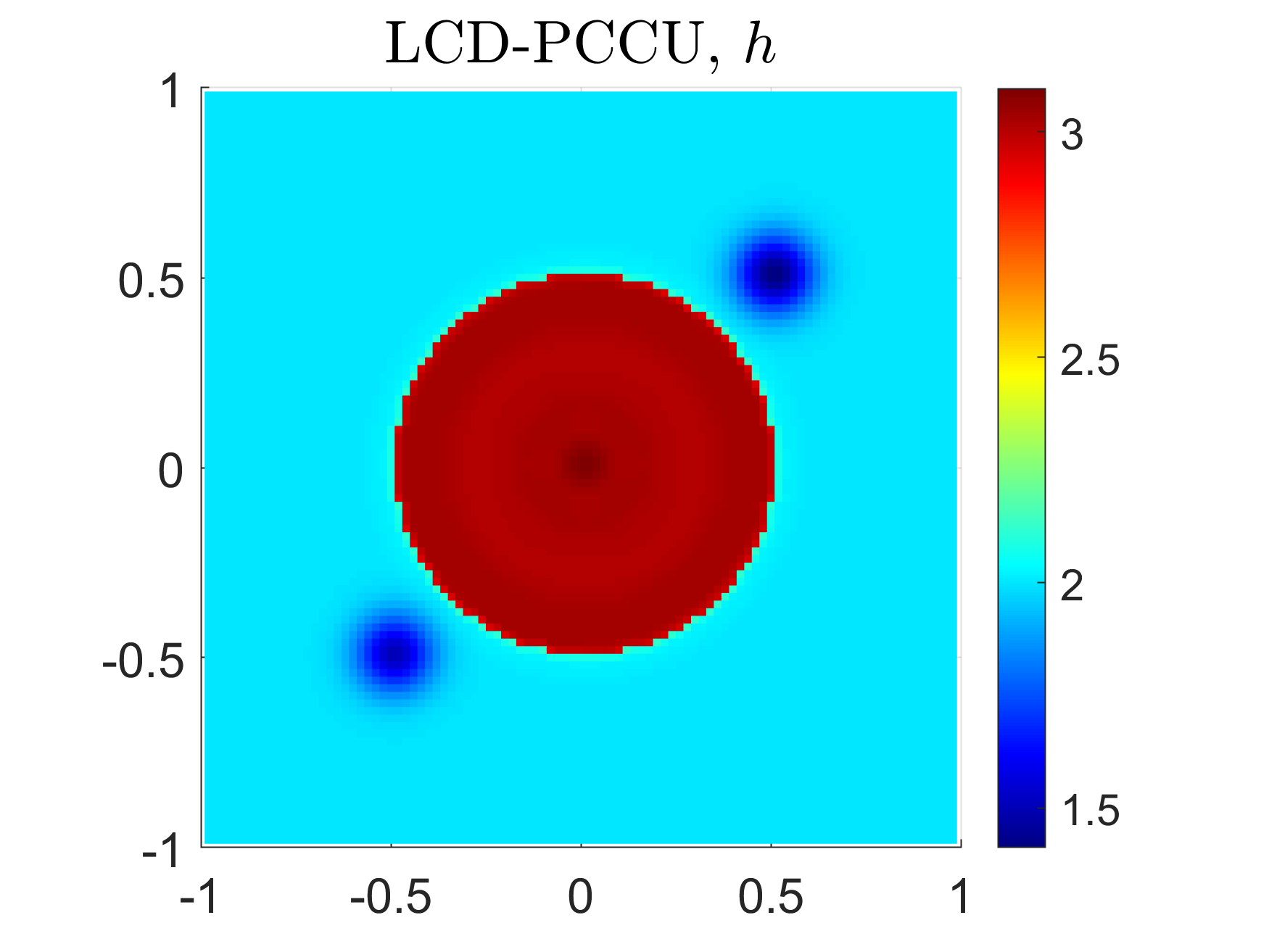}}
\vskip7pt
\centerline{\includegraphics[trim=1.2cm 0.4cm 1.7cm 0.2cm, clip, width=5cm]{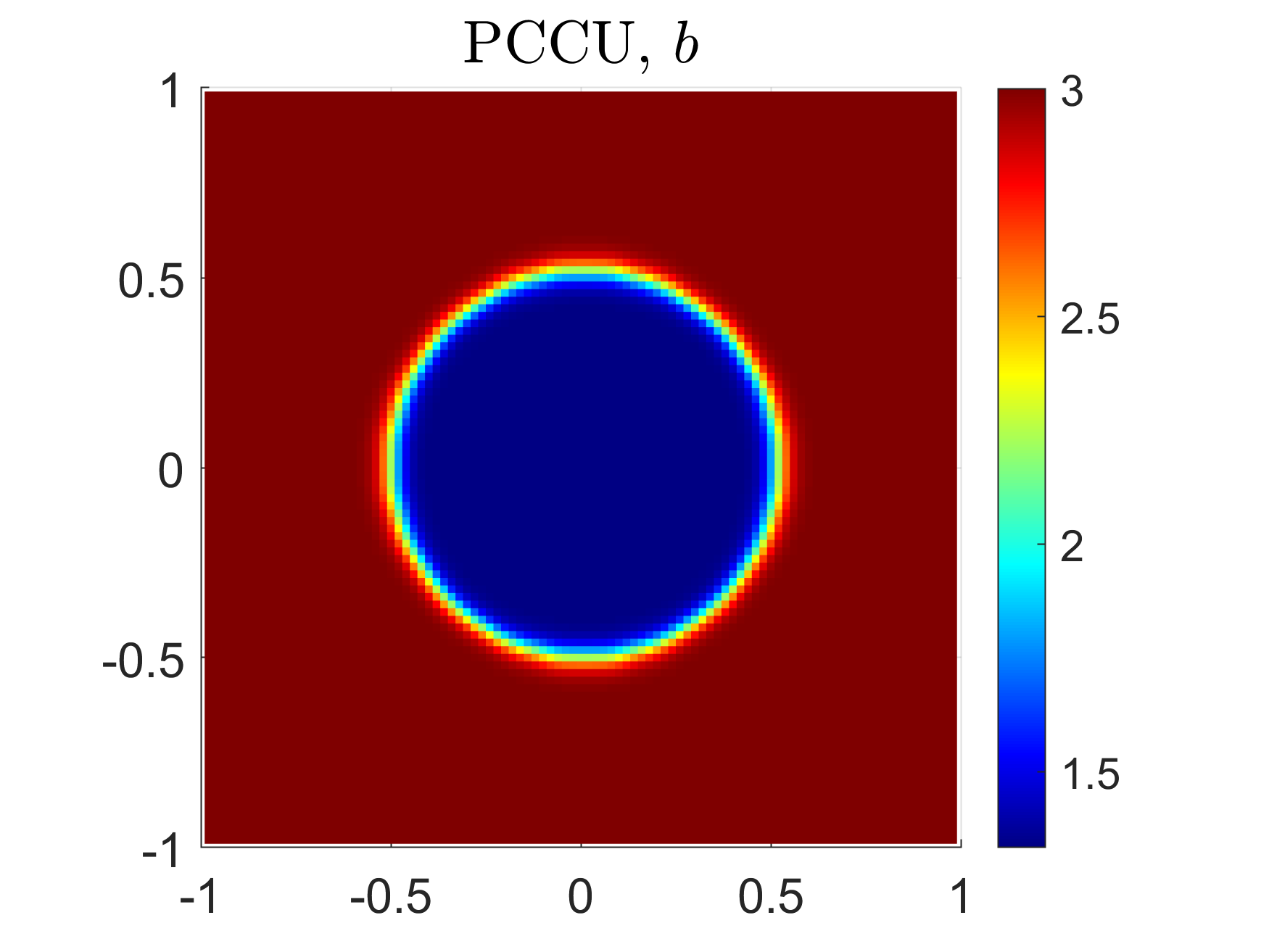}\hspace*{1.0cm}
            \includegraphics[trim=1.2cm 0.4cm 1.7cm 0.2cm, clip, width=5cm]{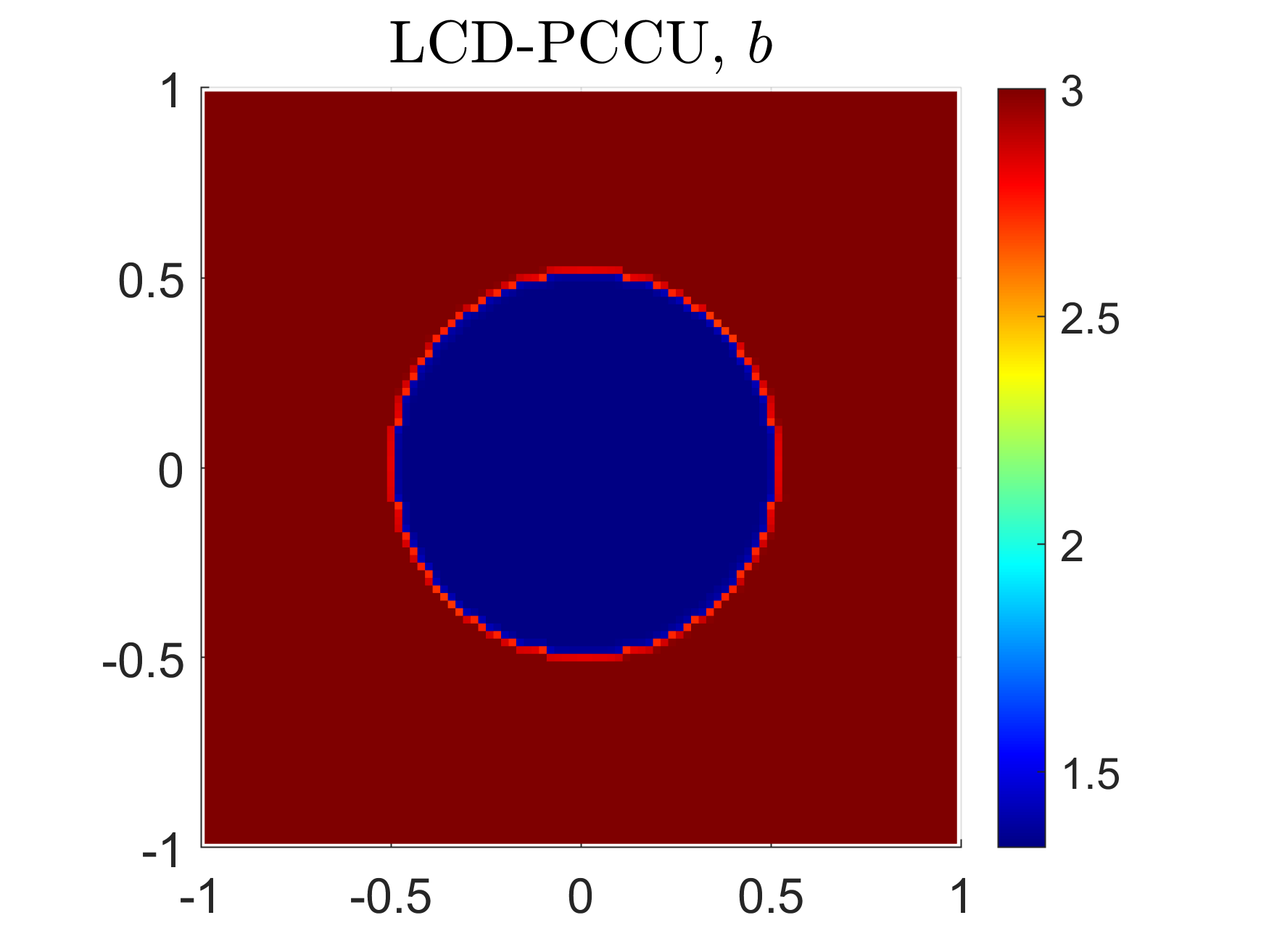}}
\caption{\sf Example 8: $h$ and $b$ computed by the PCCU (left column) and LCD-PCCU (right column) schemes.\label{fig101}}
\end{figure}

\begin{figure}[ht!]
\centerline{\includegraphics[trim=0.8cm 0.3cm 1.3cm 0.2cm, clip, width=5cm]{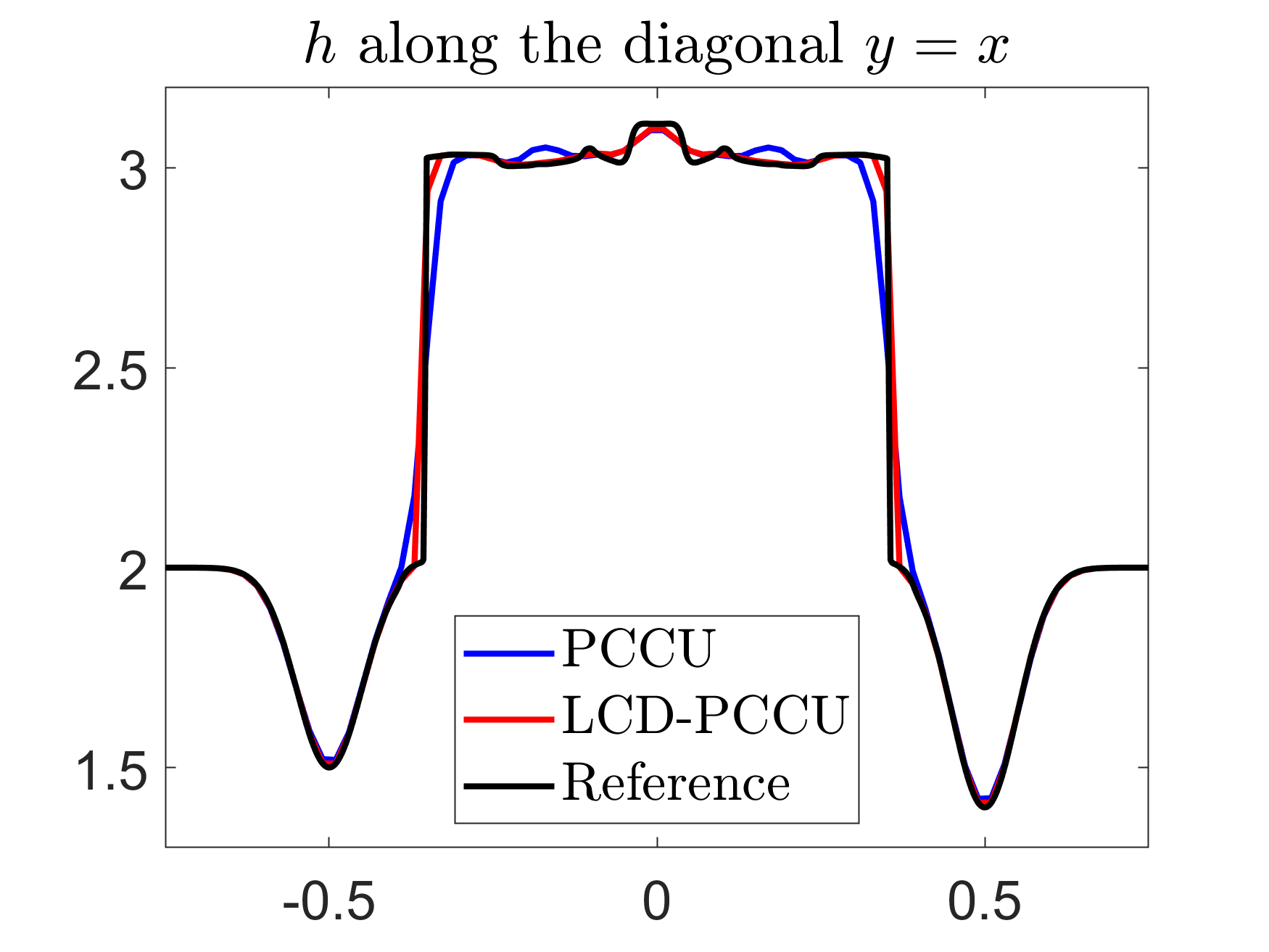}\hspace*{1.0cm}
            \includegraphics[trim=0.8cm 0.3cm 1.3cm 0.2cm, clip, width=5cm]{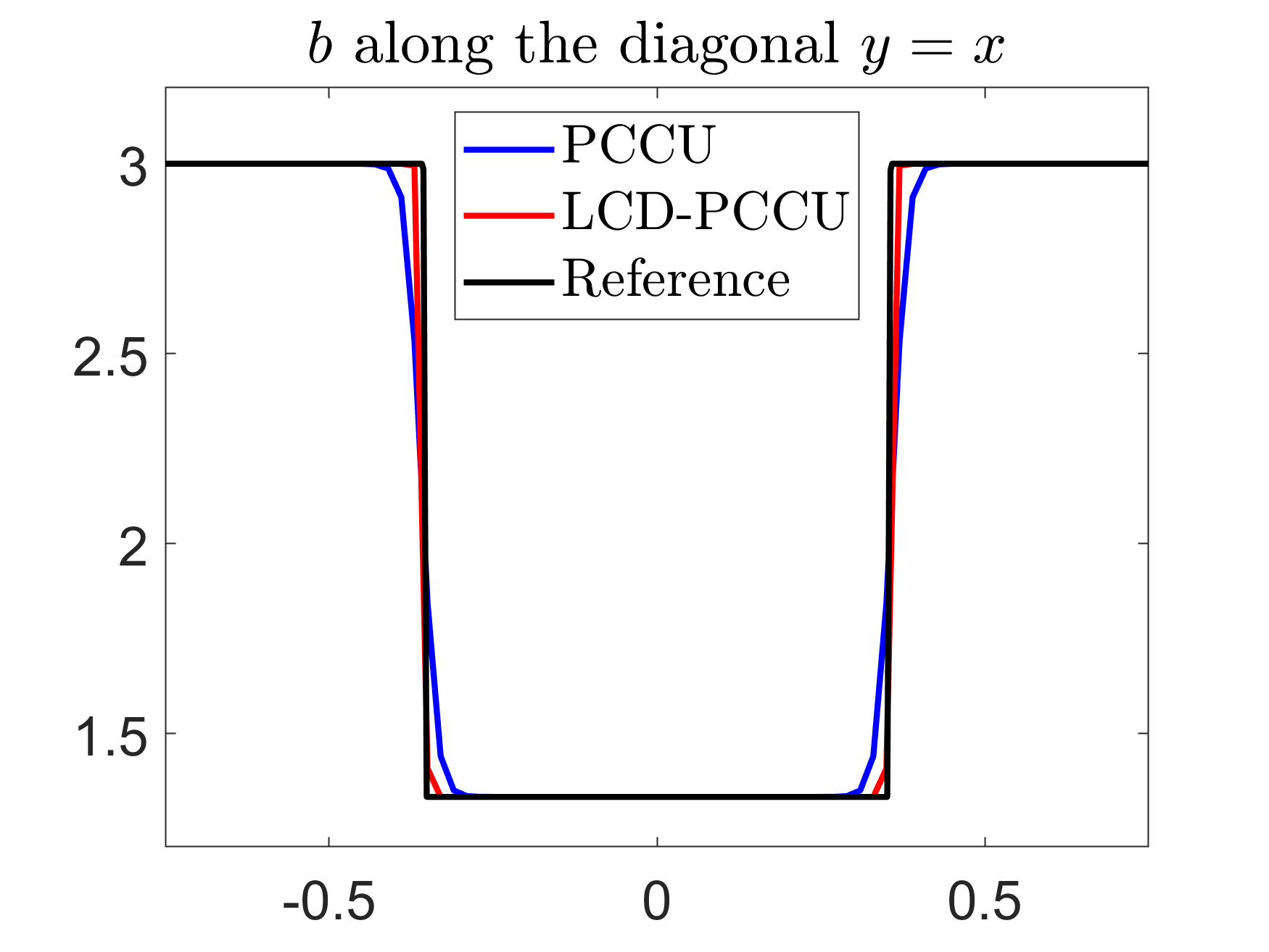}}
\caption{\sf Example 8: Diagonal slices of $h$ (left) and $b$ (right) computed by the PCCU and LCD-PCCU schemes; zoom at $[-0.75,0.75]$.
\label{fig102}}
\end{figure}

We compute the numerical solution by the PCCU and LCD-PCCU schemes until the final time $t=0.12$ on a uniform mesh with $\dx=\dy=1/50$. The
obtained results are presented in Figure \ref{fig101}, where one can clearly see that the LCD-PCCU scheme outperforms the PCCU one, since
the $h$- and $b$-components computed by the PCCU scheme are smeared; see Figure \ref{fig101}. We also plot the 1-D slices of $h$ and $b$
along the diagonal $y=x$; see Figure \ref{fig102}, where the reference solution is computed by the PCCU scheme on a much finer mesh with
$\dx=\dy=1/400$. As one can see, the LCD-PCCU scheme achieves higher resolution than the PCCU one.

\subsubsection*{Example 9---Radial Dam Break Over the Flat Bottom}
In this example taken from \cite{CKL14}, we compare the PCCU and LCD-PCCU schemes on a circular dam-break problem with a flat bottom
topography ($Z(x,y)\equiv0$) and no Coriolis forces $(f(y)\equiv0)$. The initial conditions,
\begin{equation*}
(h,u,v,b)(x,y,0)=\begin{cases}(1.5,0,0,1)&\mbox{if}~x^2+y^2<0.25,\\(1.2,0,0,1.5)&\mbox{otherwise},\end{cases}
\end{equation*}
are prescribed in the computational domain $[-1,1]\times[-1,1]$ subject to the free boundary conditions.

\begin{figure}[ht!]
\centerline{\includegraphics[trim=1.2cm 0.4cm 1.5cm 0.2cm, clip, width=5cm]{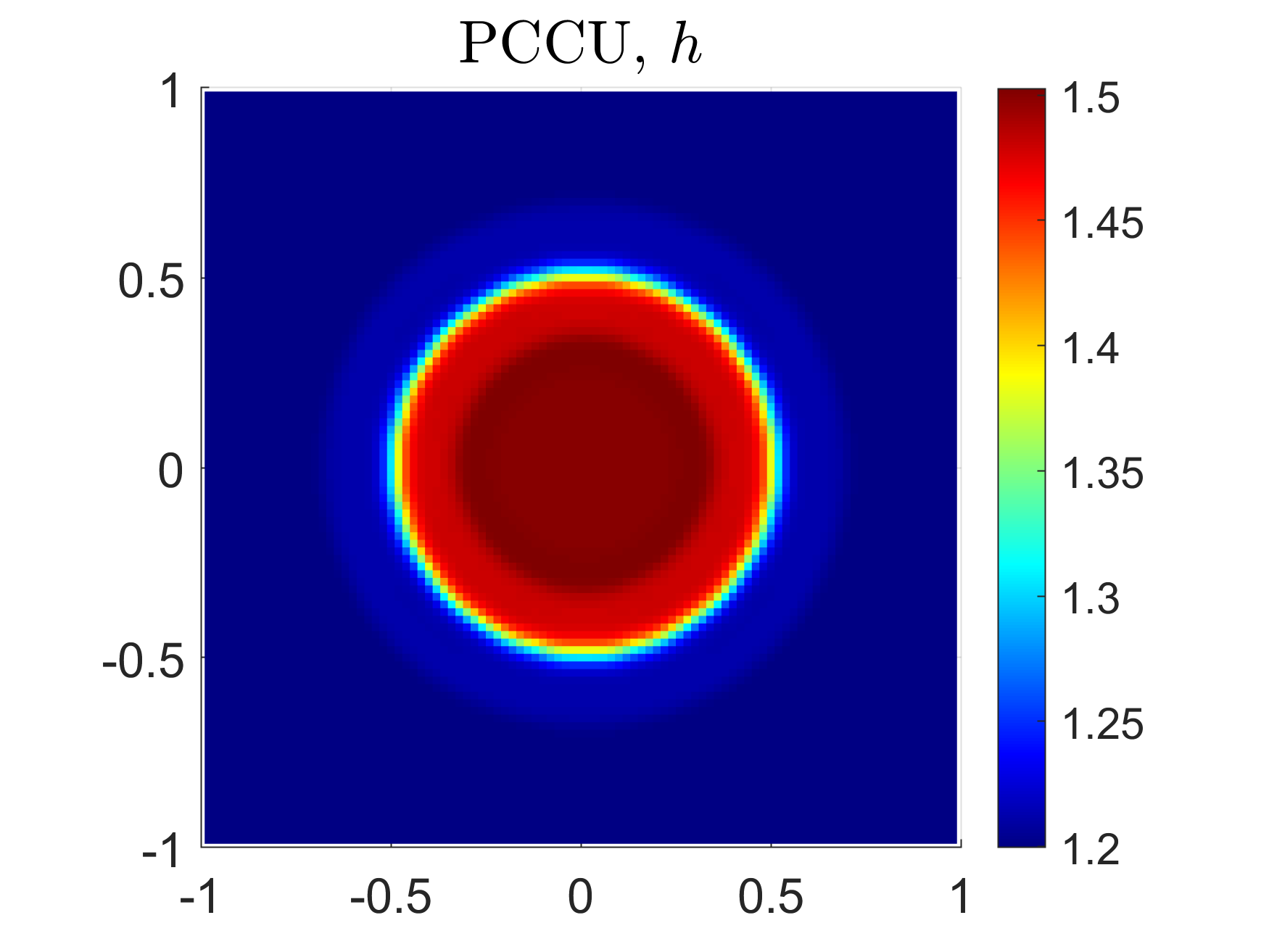}\hspace*{1.0cm}
            \includegraphics[trim=1.2cm 0.4cm 1.5cm 0.2cm, clip, width=5cm]{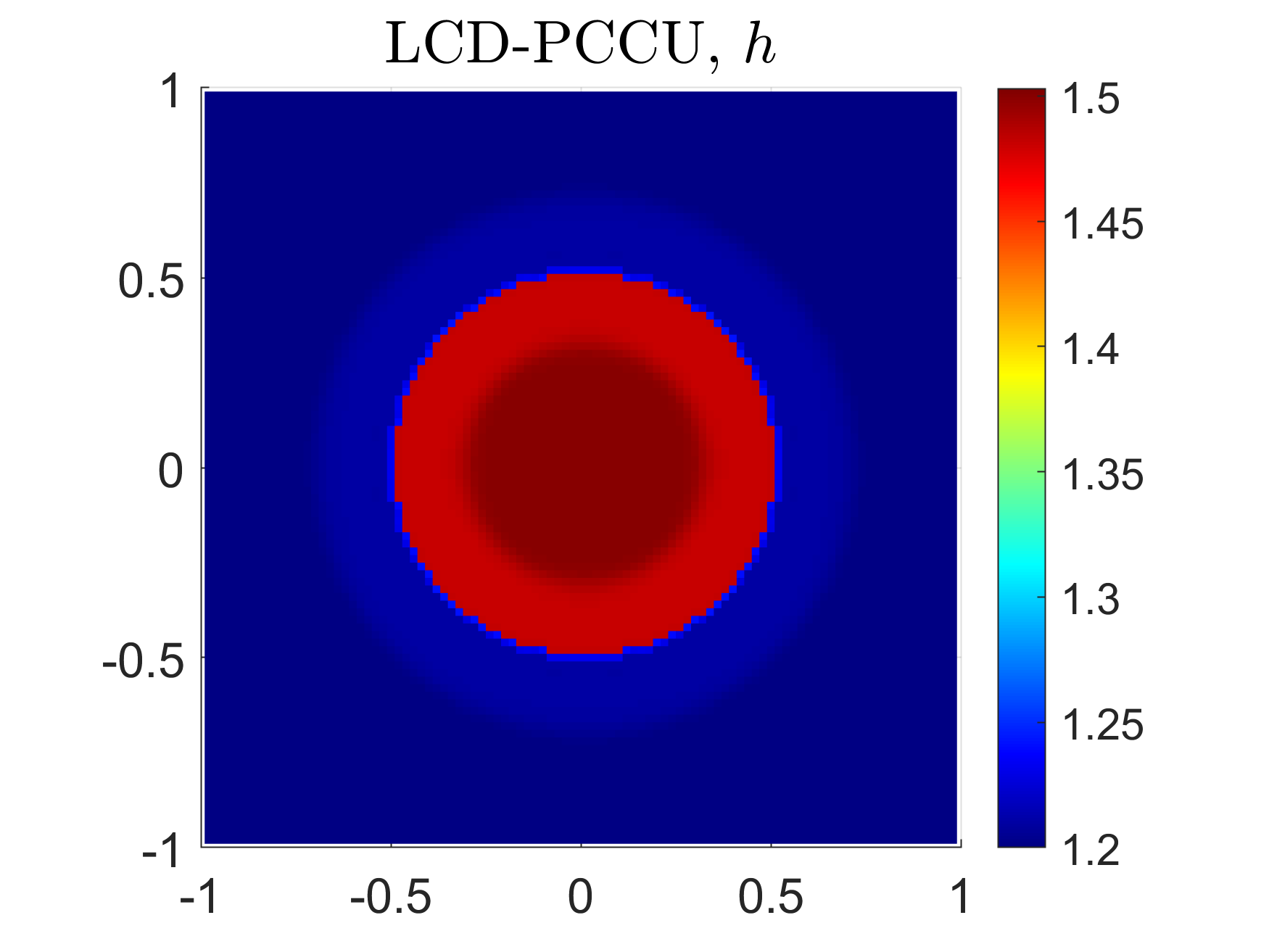}}
\vskip7pt
\centerline{\includegraphics[trim=1.2cm 0.4cm 1.5cm 0.2cm, clip, width=5cm]{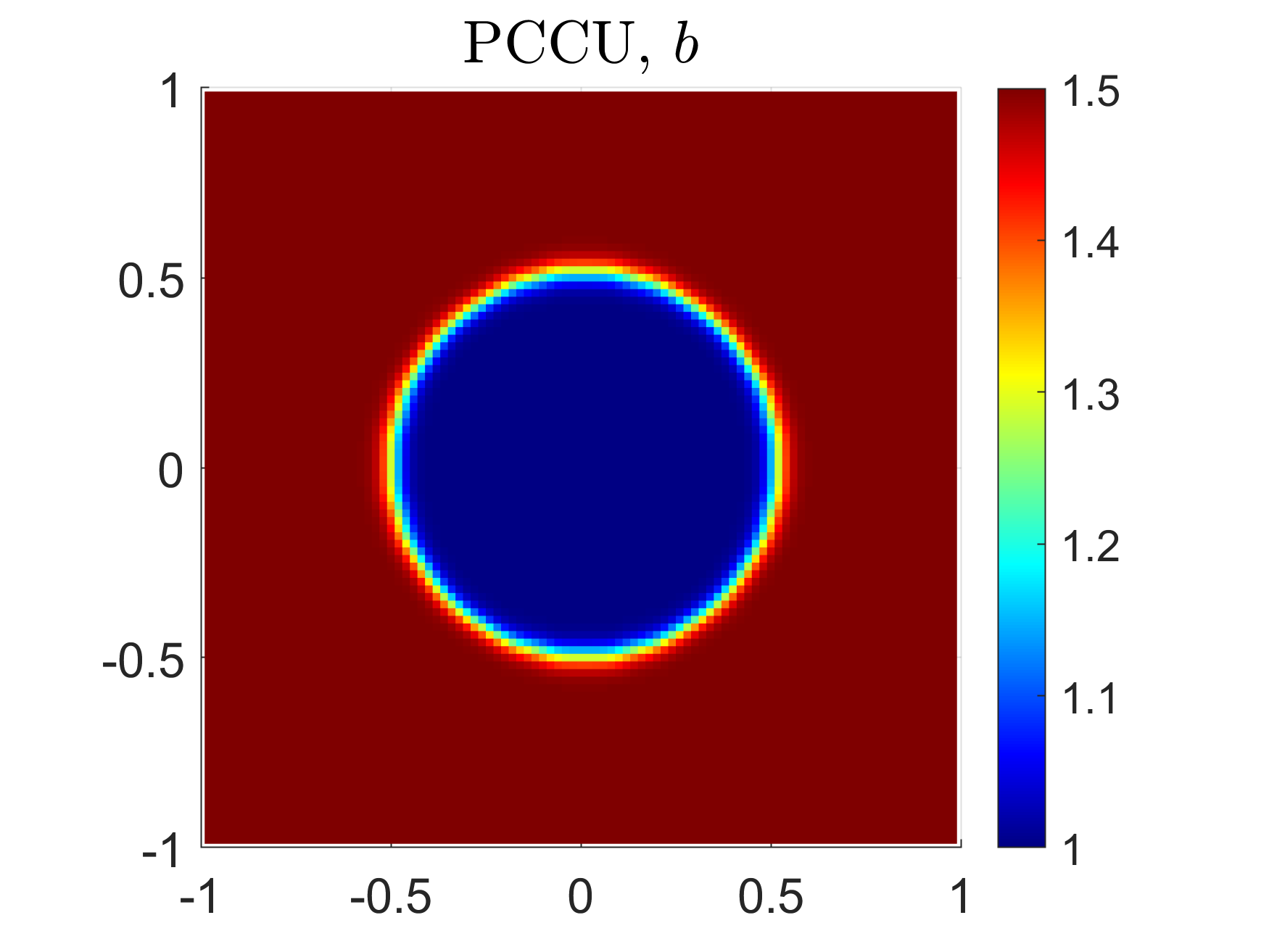}\hspace*{1.0cm}
            \includegraphics[trim=1.2cm 0.4cm 1.5cm 0.2cm, clip, width=5cm]{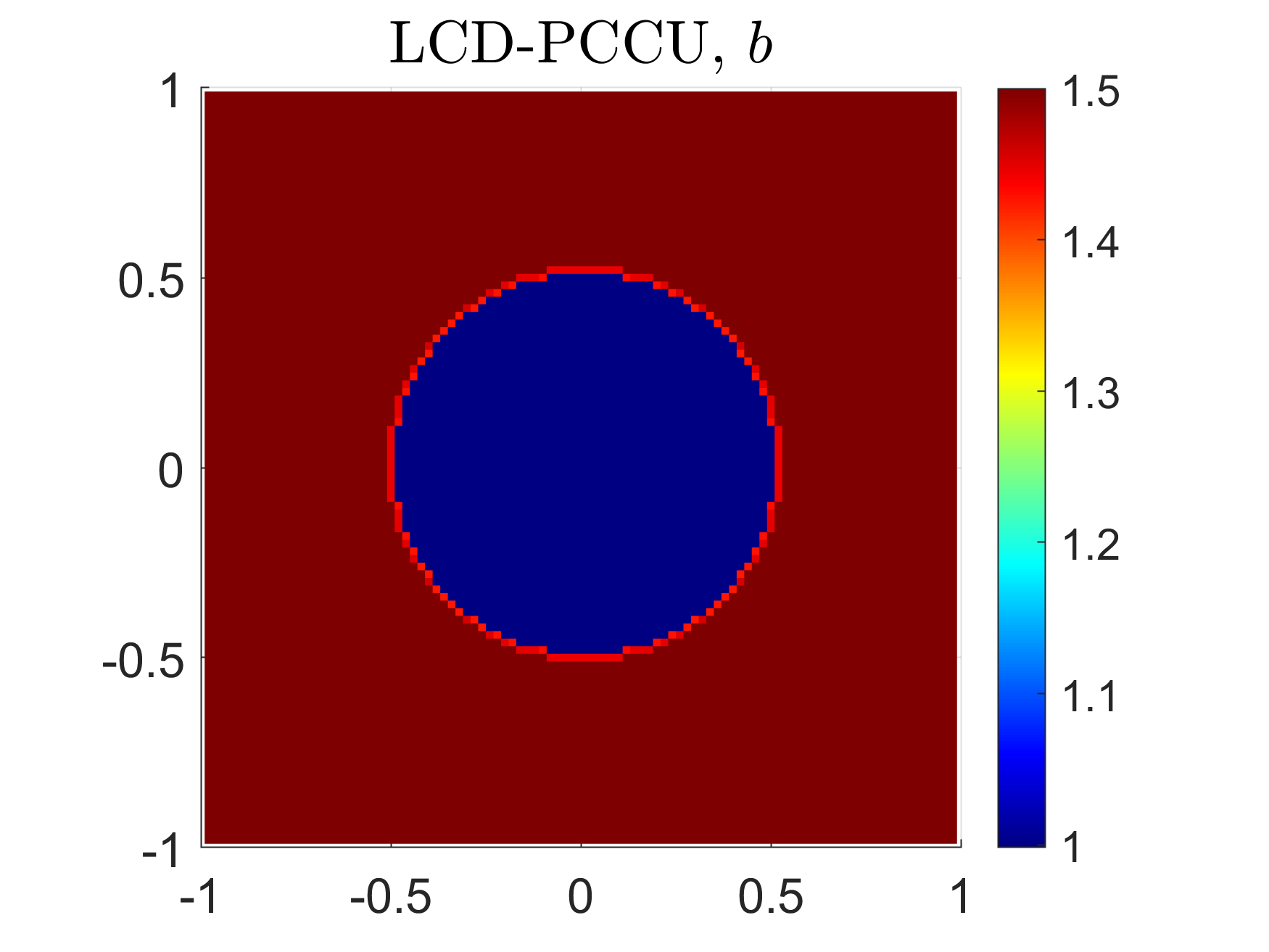}}
\caption{\sf Example 9: $h$ and $b$ computed by the PCCU (left column) and LCD-PCCU (right column) schemes.\label{fig111}}
\end{figure}

\begin{figure}[ht!]
\centerline{\includegraphics[trim=0.8cm 0.3cm 1.0cm 0.2cm, clip, width=5cm]{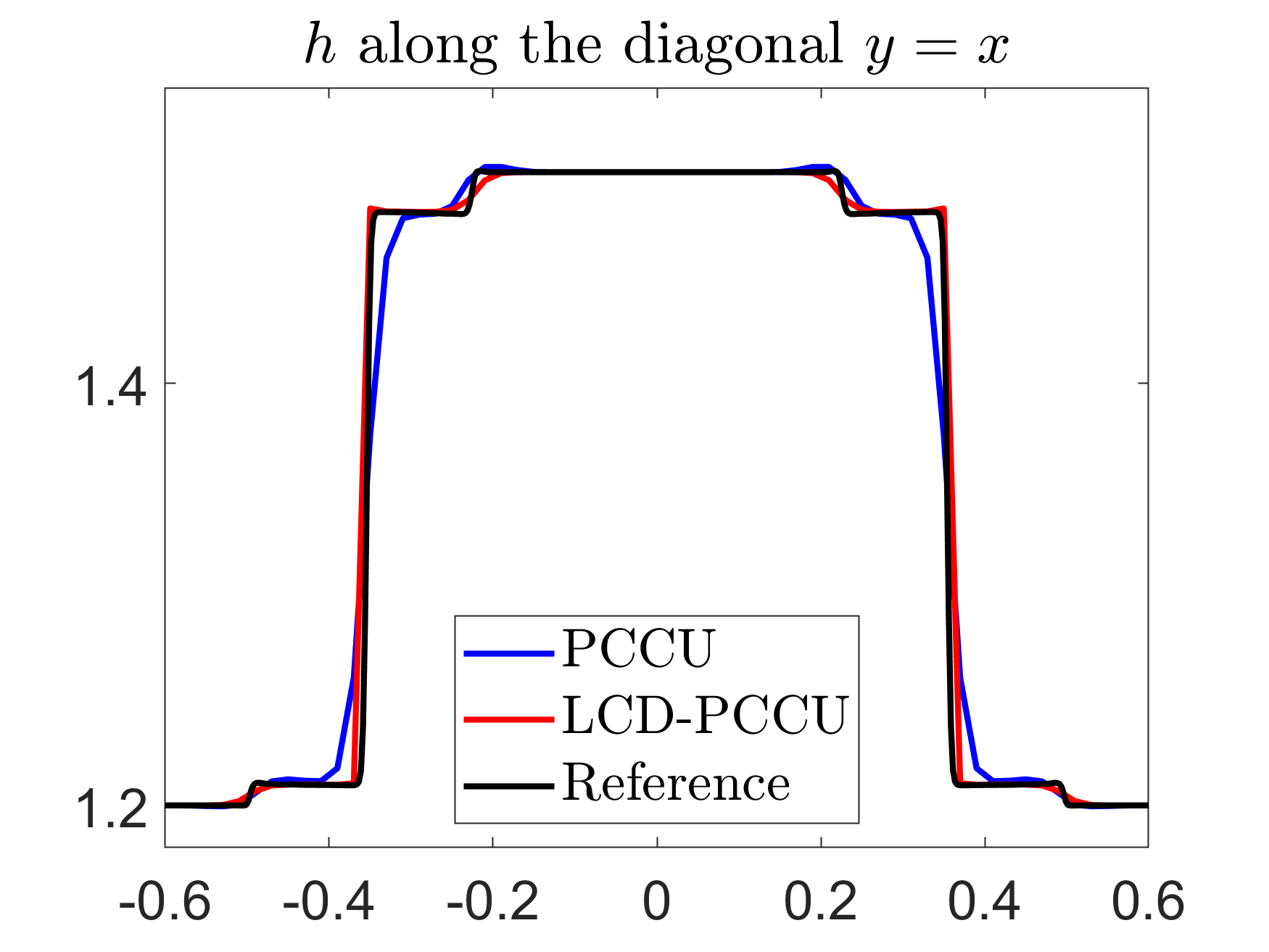}\hspace*{1.0cm}
            \includegraphics[trim=0.8cm 0.3cm 1.0cm 0.2cm, clip, width=5cm]{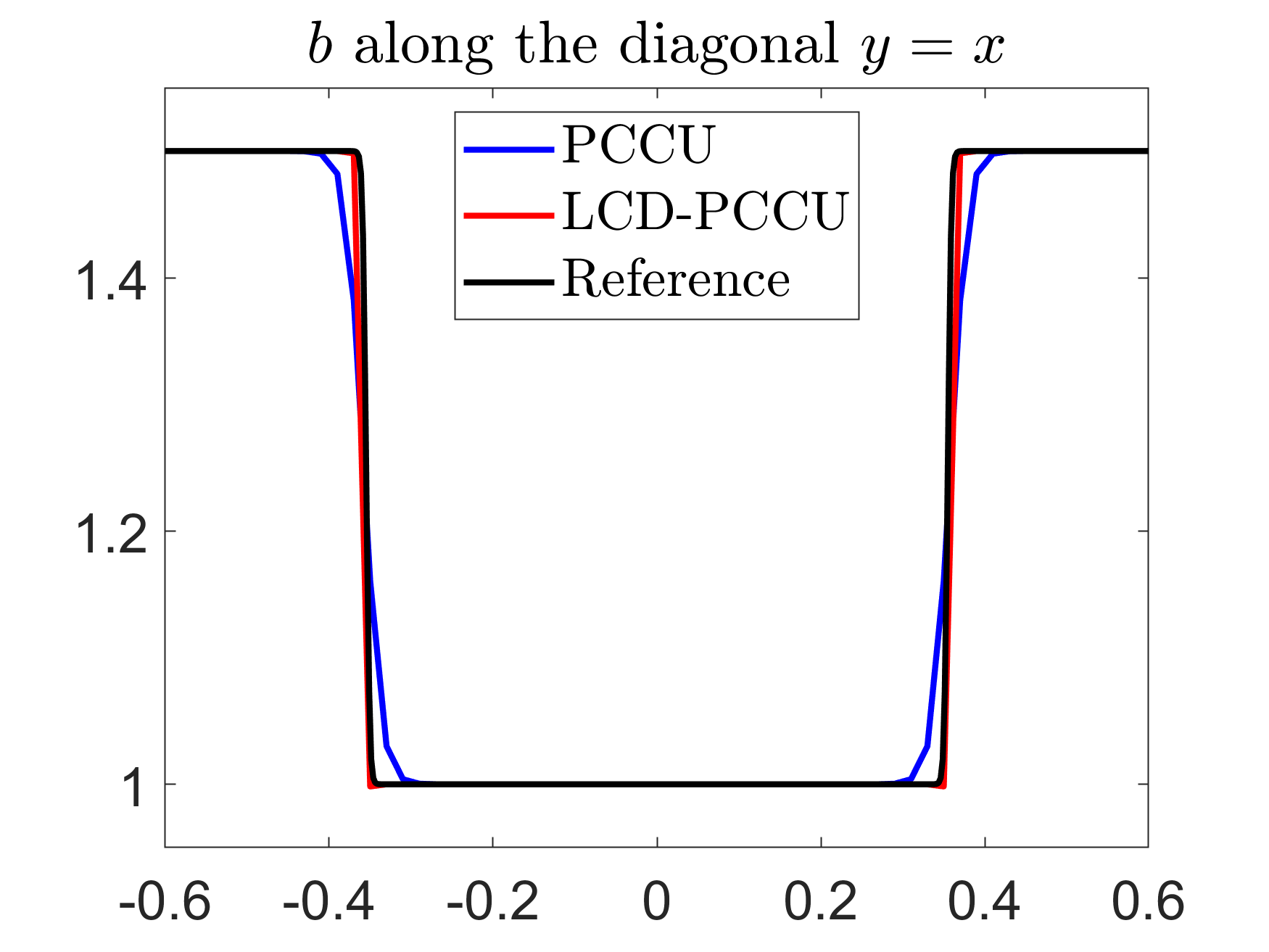}}
\caption{\sf Example 9: Diagonal slices of $h$ (left) and $b$ (right) computed by the PCCU and LCD-PCCU schemes; zoom at $[-0.6,0.6]$.
\label{fig112}}
\end{figure}

We compute the numerical solution by the PCCU and LCD-PCCU schemes until the final time $t=0.15$ on a uniform mesh with $\dx=\dy=1/50$. The
computed results are presented in Figure \ref{fig111}, where one can clearly see that the LCD-PCCU scheme outperforms the PCCU scheme. We
also plot the 1-D slices of $h$ and $b$ along the diagonal $y=x$; see Figure \ref{fig112}, where the reference solution is computed by the
PCCU scheme on a much finer mesh with $\dx=\dy=1/400$. As one can see, the LCD-PCCU scheme outperforms its PCCU counterpart.

\subsubsection*{Example 10---Relaxation of Localized Anomalies on the Equatorial $\beta$-Plane}
In the last example taken from \cite{KLZ_2D}, we consider the relaxation of localized anomalies on the $\beta$-plane, where the Coriolis
parameter is $f(y)=y$, with the initial data
\begin{equation*}
h(x,y,0)\equiv1,\quad q(x,y,0)=p(x,y,0)\equiv0,\quad b(x,y,0)=1+\hf\,e^{-\left(\frac{1}{50}x^2+\frac{1}{2}y^2\right)},
\end{equation*}
prescribed in the computational domain $[-40,80]\times[-10,10]$ subject to the free boundary conditions.

\begin{figure}[ht!]
\centerline{\includegraphics[trim=0.7cm 0.4cm 0.7cm 0.2cm, clip, width=5cm]{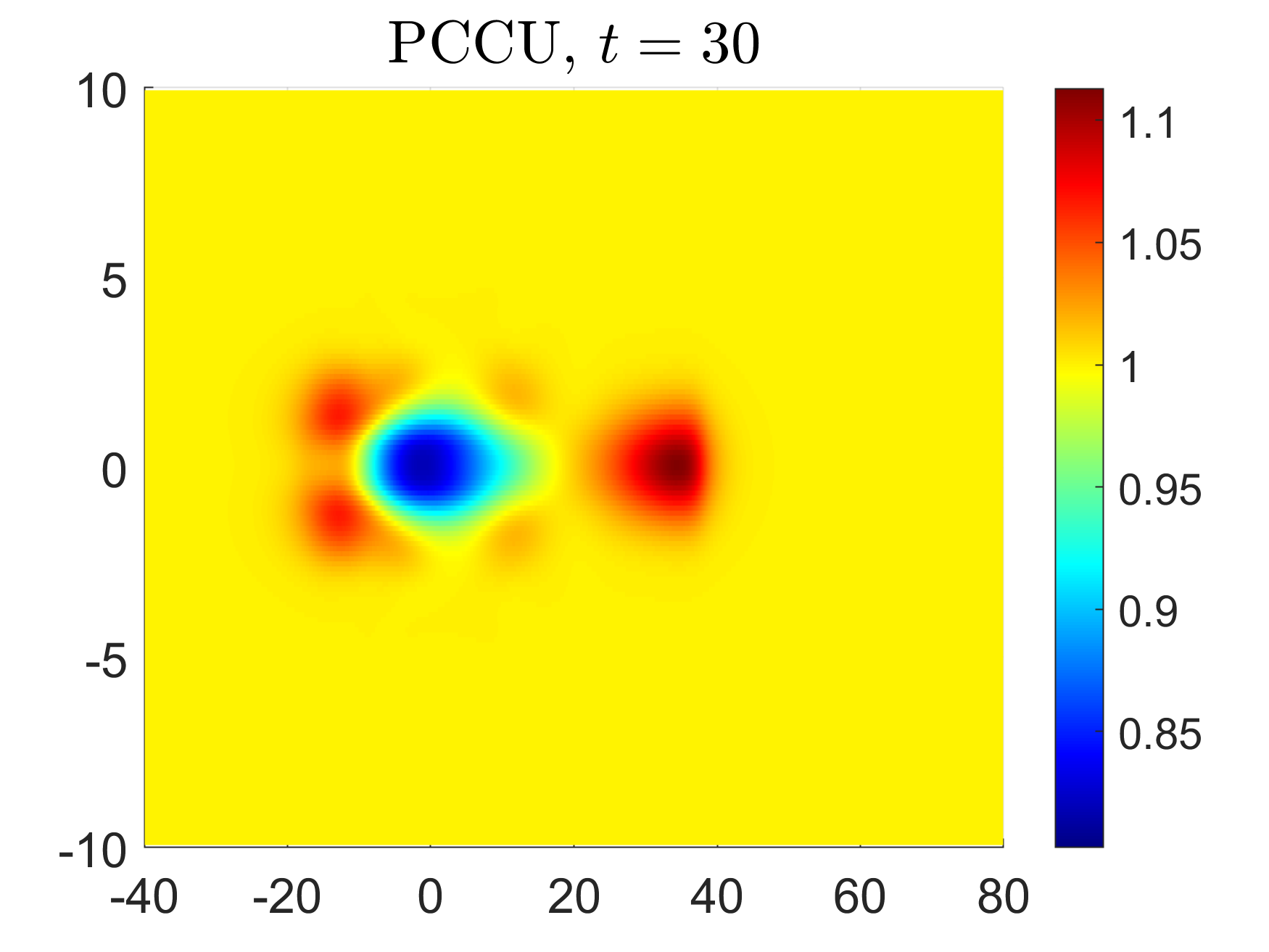}\hspace*{1.0cm}
            \includegraphics[trim=0.7cm 0.4cm 0.7cm 0.2cm, clip, width=5cm]{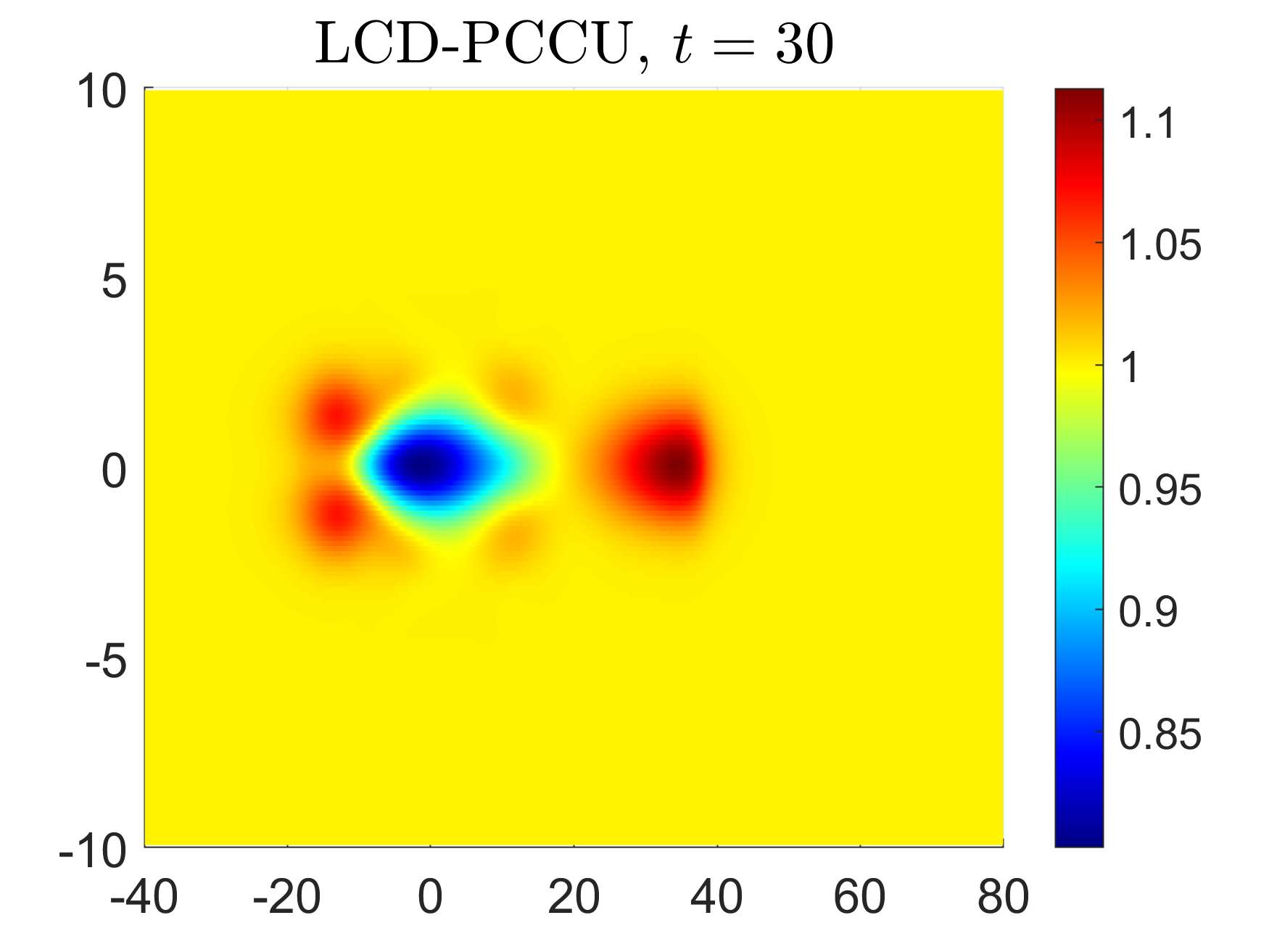}}
\vskip8pt
\centerline{\includegraphics[trim=0.7cm 0.4cm 0.7cm 0.2cm, clip, width=5cm]{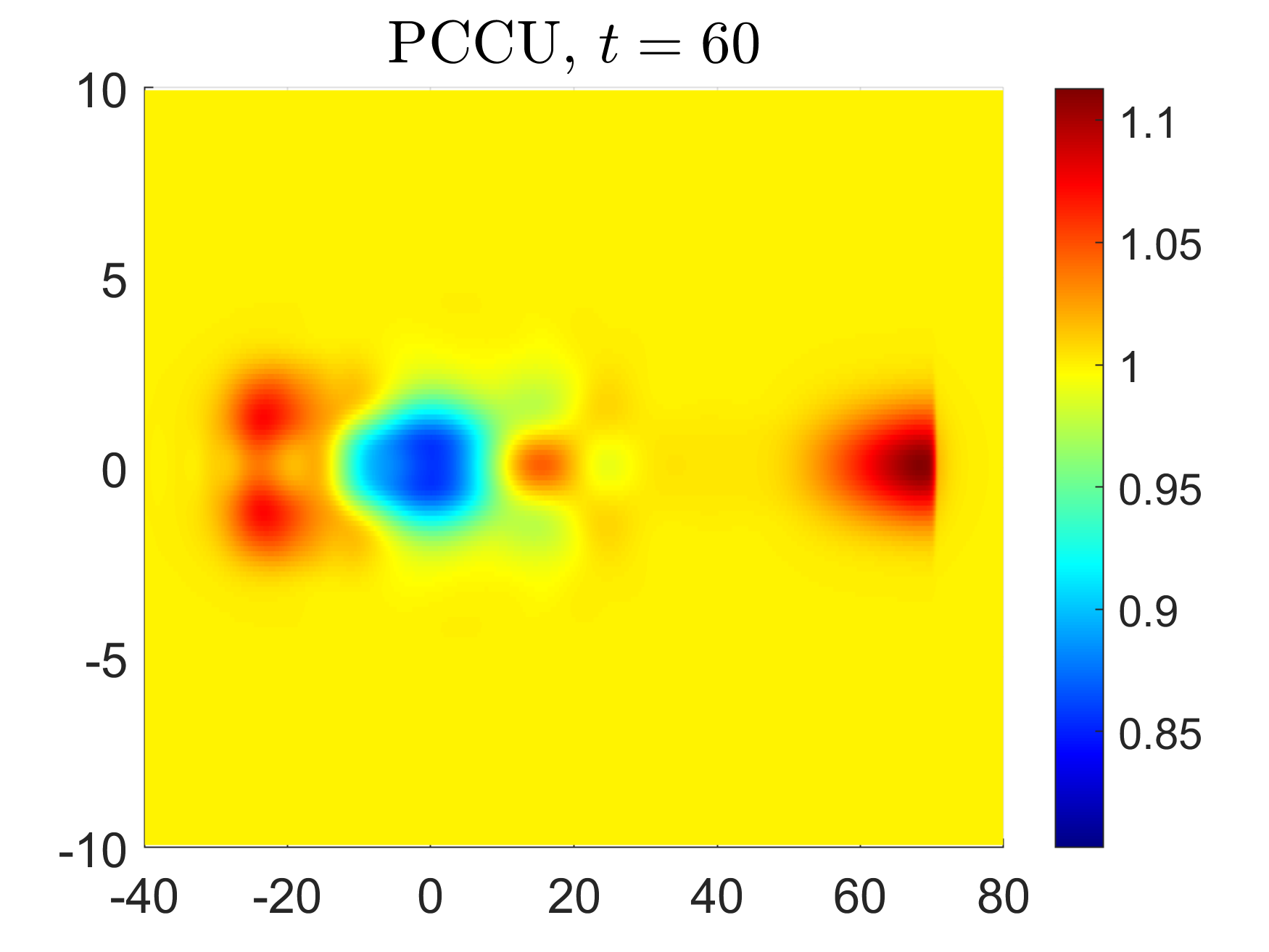}\hspace*{1.0cm}
            \includegraphics[trim=0.7cm 0.4cm 0.7cm 0.2cm, clip, width=5cm]{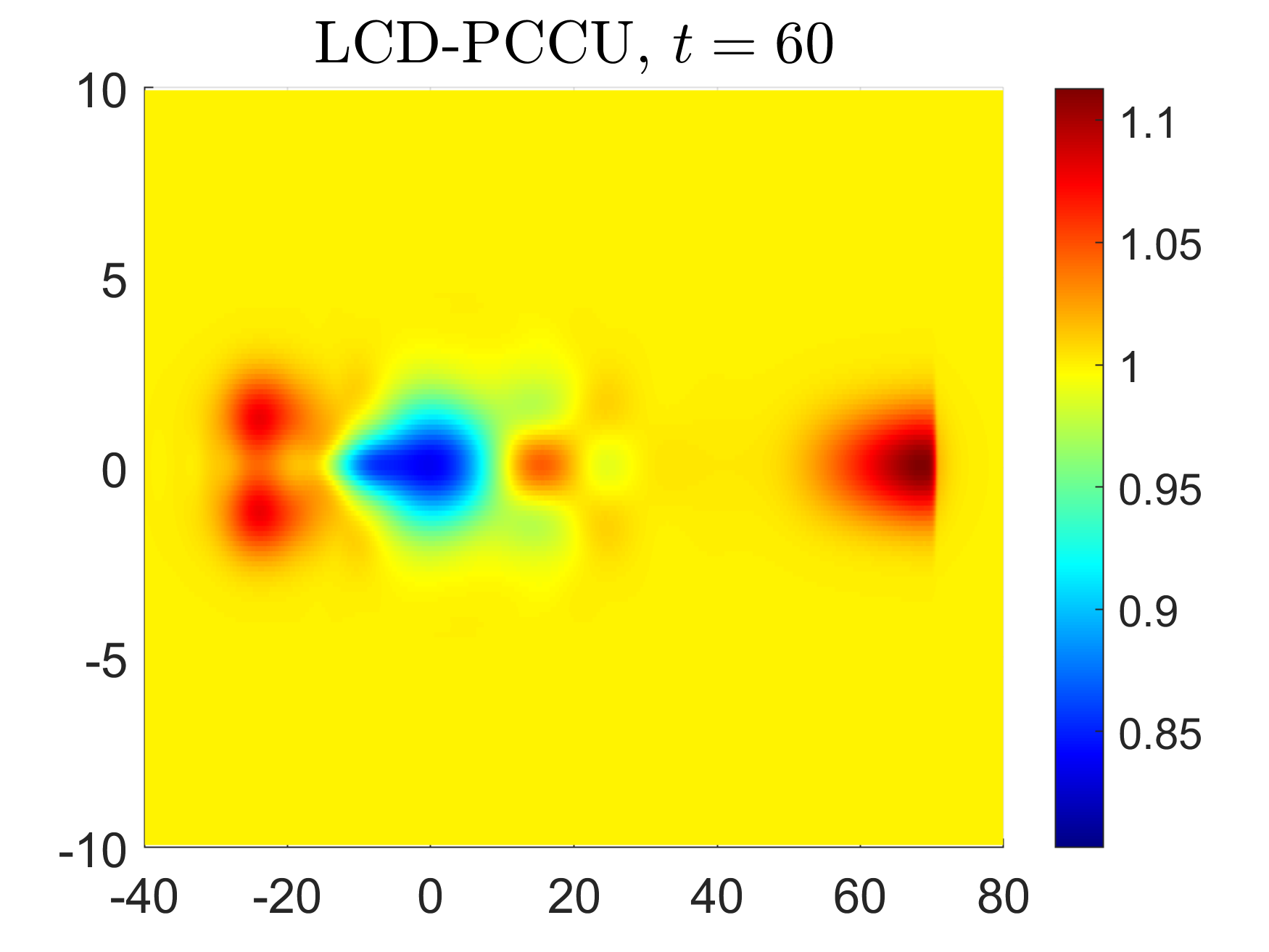}}
\vskip8pt
\centerline{\includegraphics[trim=0.7cm 0.4cm 0.7cm 0.2cm, clip, width=5cm]{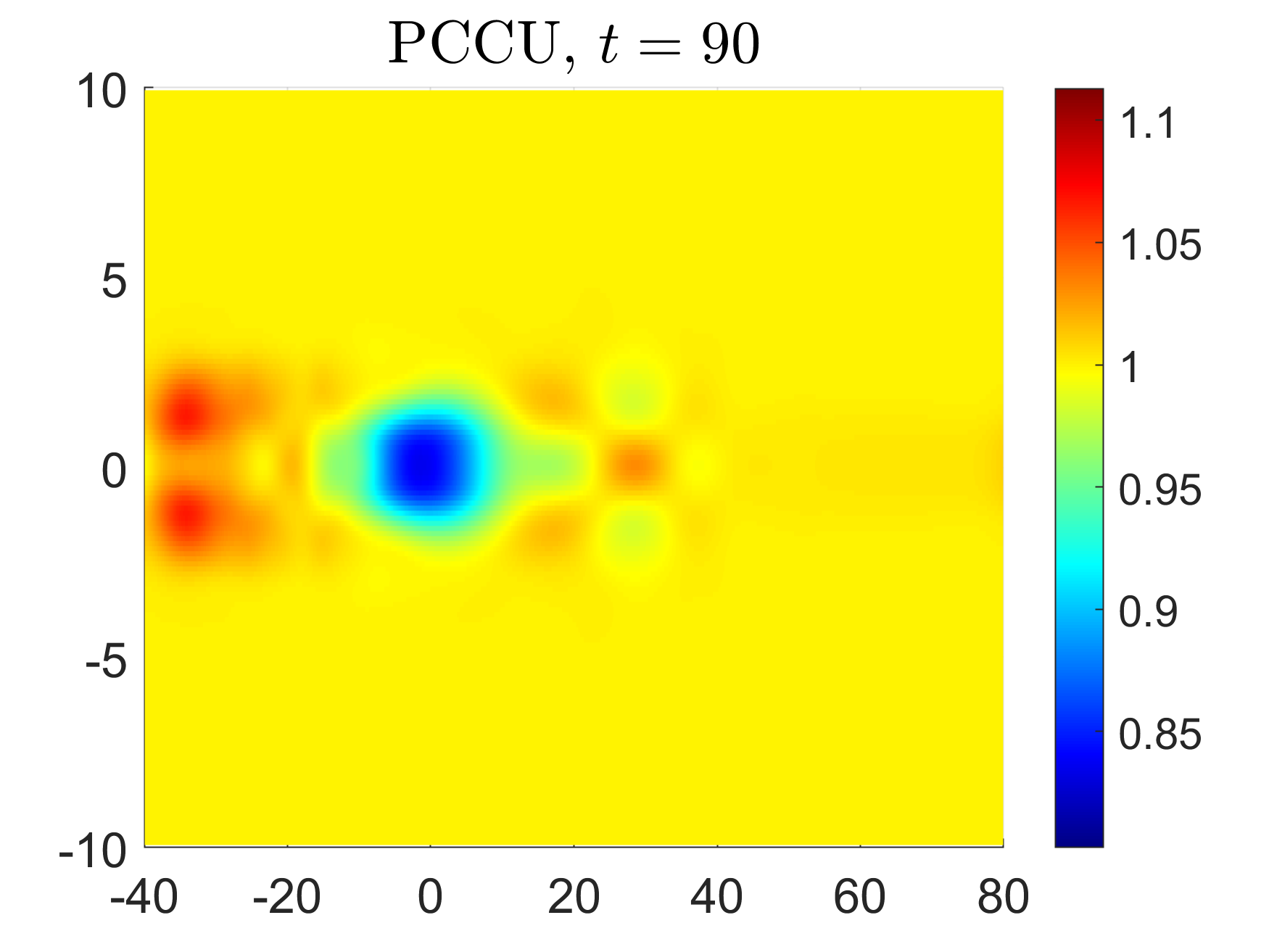}\hspace*{1.0cm}
            \includegraphics[trim=0.7cm 0.4cm 0.7cm 0.2cm, clip, width=5cm]{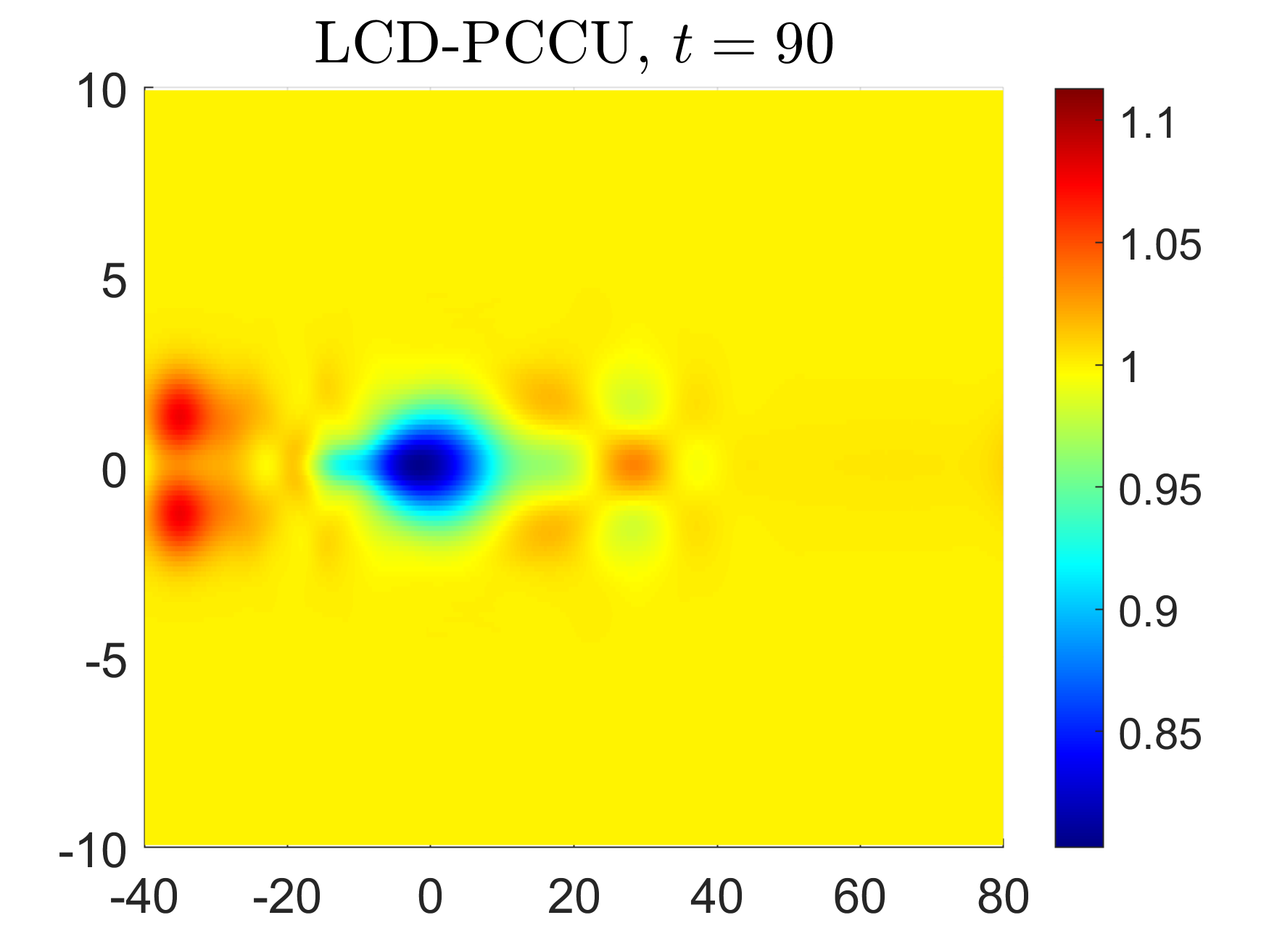}}
\vskip8pt
\centerline{\includegraphics[trim=0.7cm 0.4cm 0.7cm 0.2cm, clip, width=5cm]{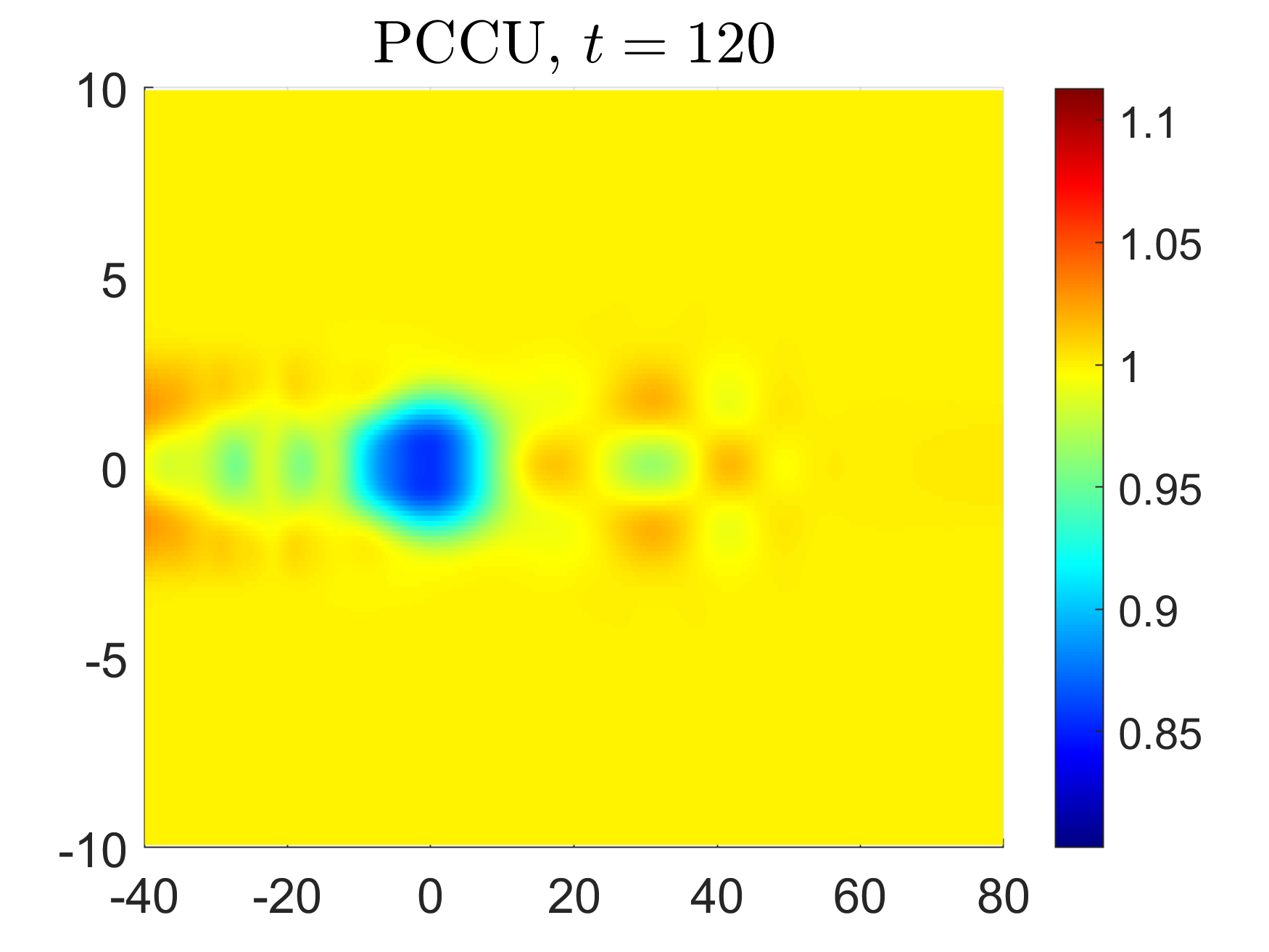}\hspace*{1.0cm}
            \includegraphics[trim=0.7cm 0.4cm 0.7cm 0.2cm, clip, width=5cm]{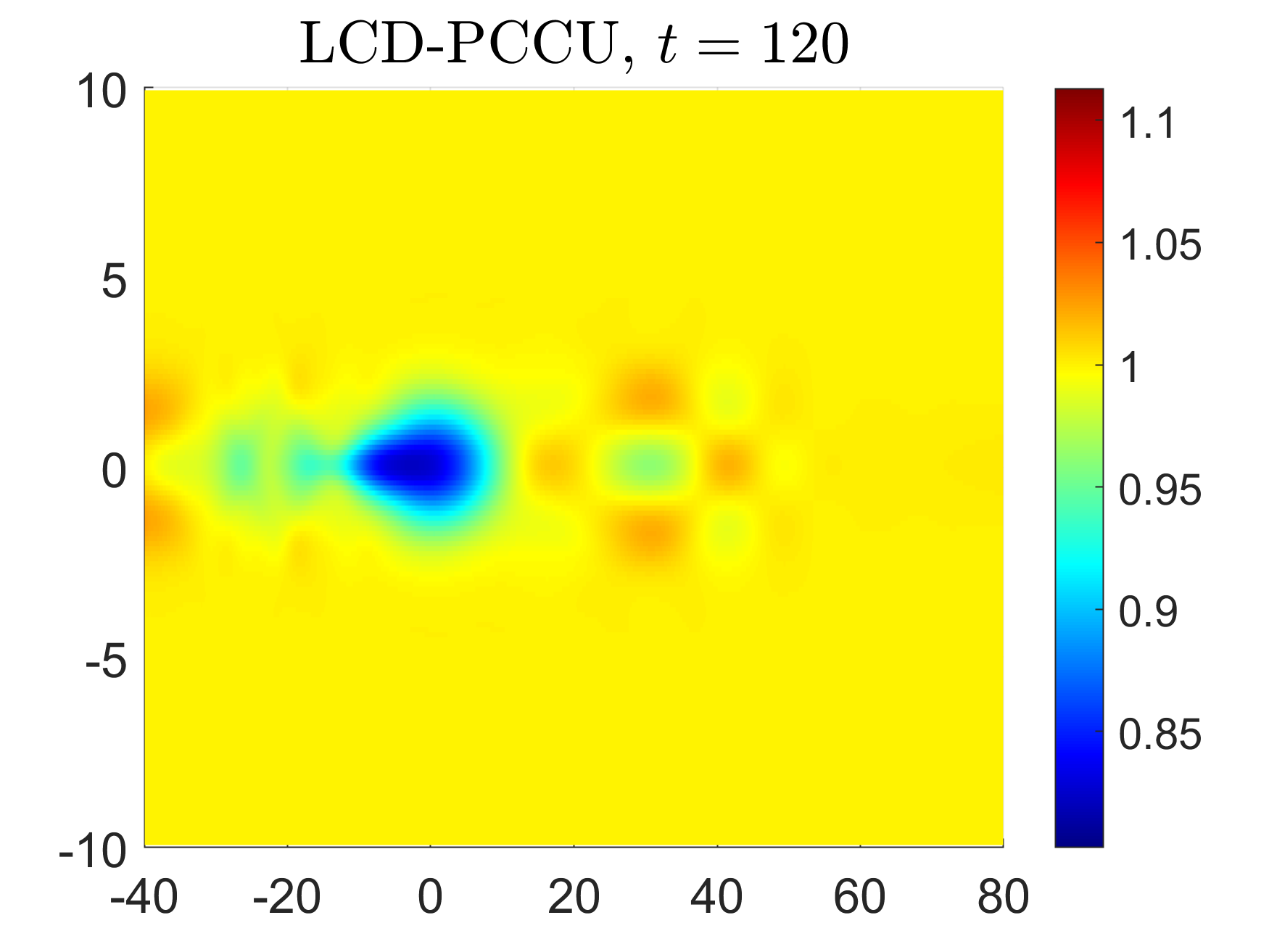}}
\caption{\sf Example 10: Time snapshots of $h$ computed by the PCCU (left column) and LCD-PCCU (right column) schemes.\label{fig131}}
\end{figure}

\begin{figure}[ht!]
\centerline{\includegraphics[trim=0.8cm 0.3cm 1.1cm 0.1cm, clip, width=5.cm]{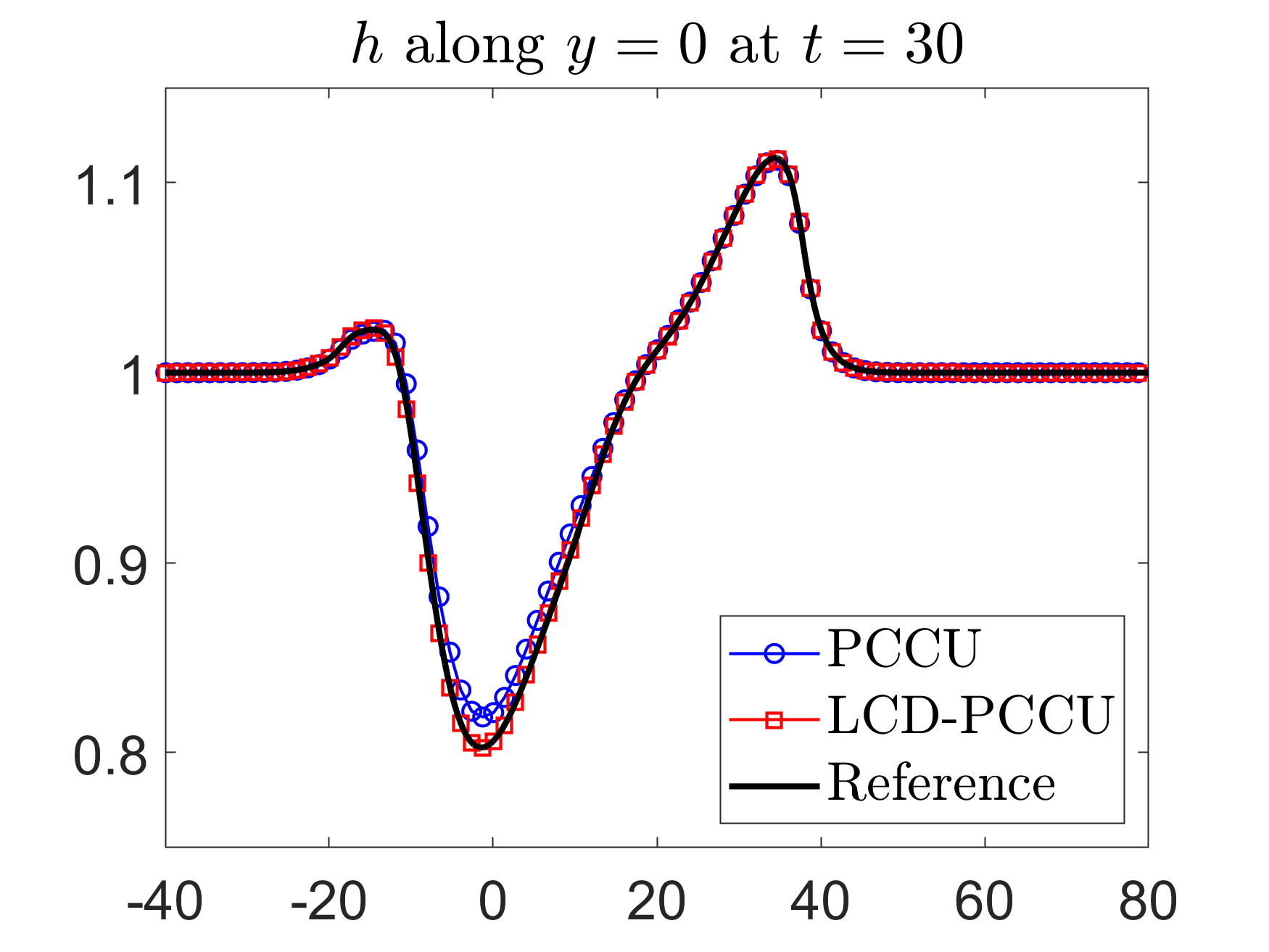}\hspace*{1.0cm}
            \includegraphics[trim=0.8cm 0.3cm 1.1cm 0.1cm, clip, width=5.cm]{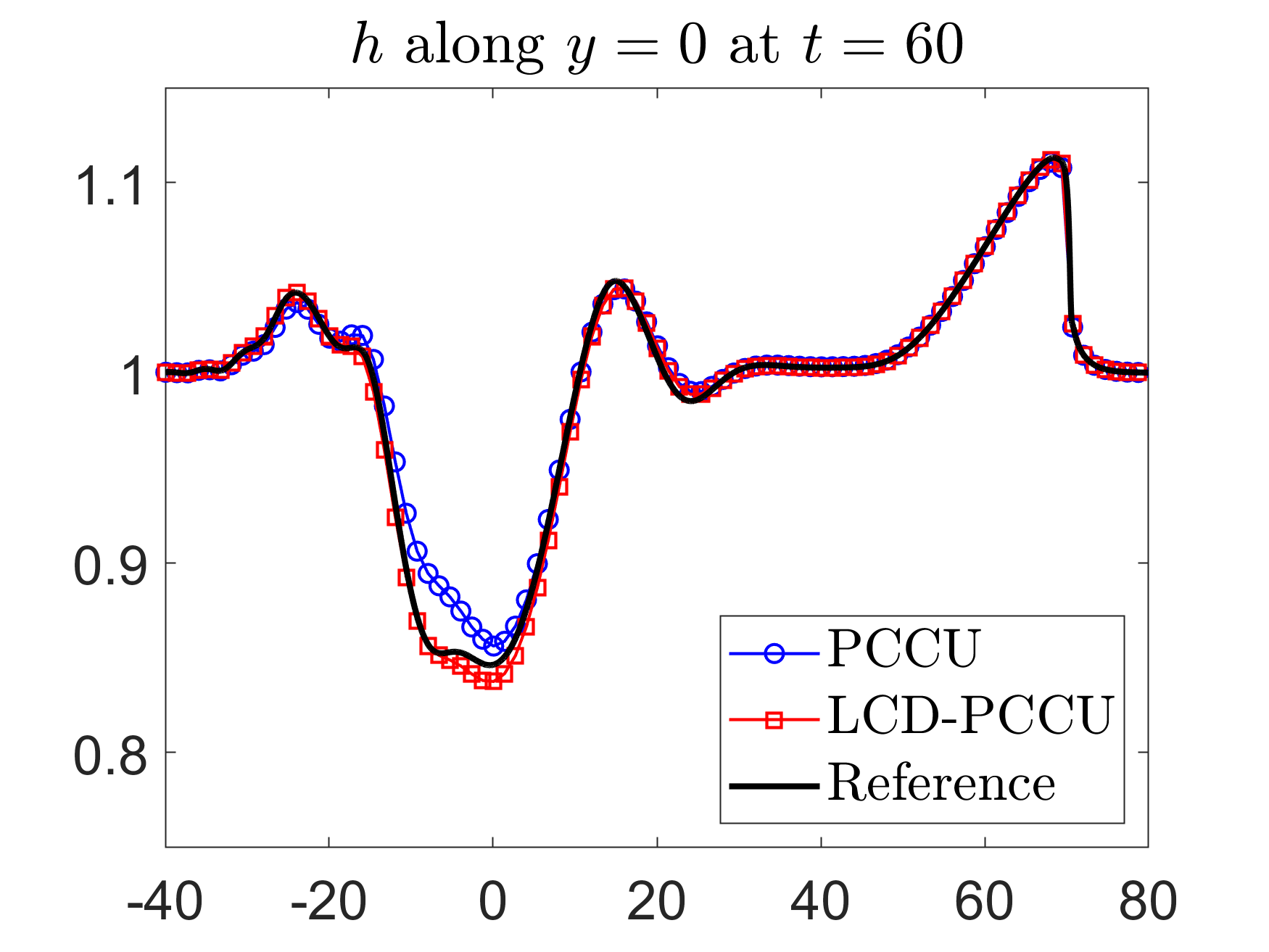}}
\vskip7pt
\centerline{\includegraphics[trim=0.8cm 0.3cm 1.1cm 0.1cm, clip, width=5.cm]{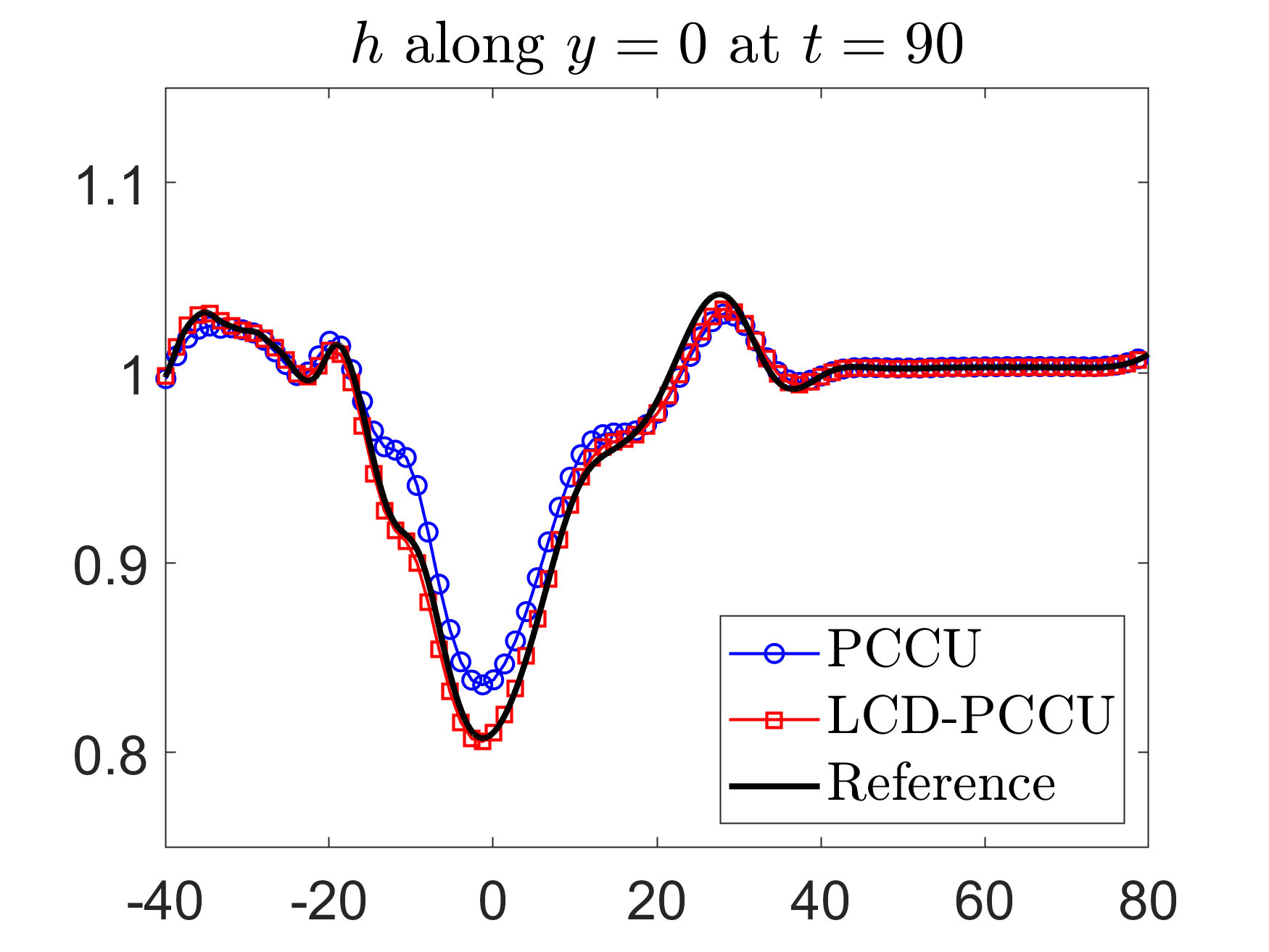}\hspace*{1.0cm}
            \includegraphics[trim=0.8cm 0.3cm 1.1cm 0.1cm, clip, width=5.cm]{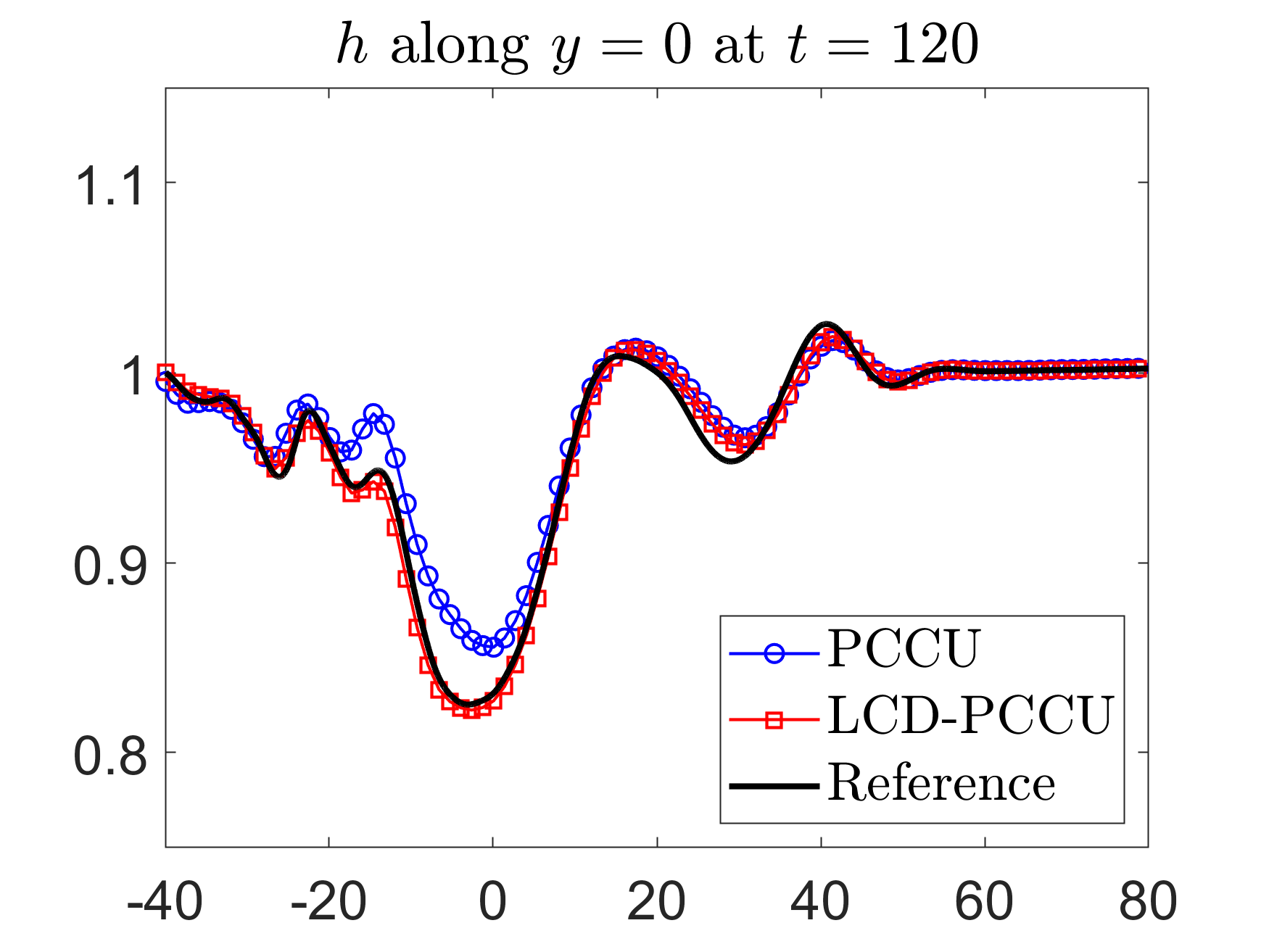}}
\caption{\sf Example 10: Slices of $h$ and $b$ along the line $y=0$ computed by the PCCU and LCD-PCCU schemes using $\dx=\dy=2/15$ at
different times.\label{fig132}}
\end{figure}

We compute the numerical solution by the PCCU and LCD-PCCU schemes until the final time $t=120$ on a uniform mesh with $\dx=\dy=2/15$ and
plot the obtained numerical results in Figure \ref{fig131}, where one can see that our results are consistent with those reported in
\cite{KLZ_2D}, and at the same time, the LCD-PCCU scheme outperforms the PCCU one. In order to better see the difference between the
computed solutions, we also plot the 1-D slices of $h$ along $y=0$ in Figure \ref{fig132}, where the reference solution is computed by the
PCCU scheme on a finer mesh with $\dx=\dy=1/25$. One can see that the LCD-PCCU scheme achieves significantly better resolution than the PCCU
one.

\begin{figure}[ht!]
\centerline{\includegraphics[trim=0.8cm 0.3cm 1.1cm 0.1cm, clip, width=5.cm]{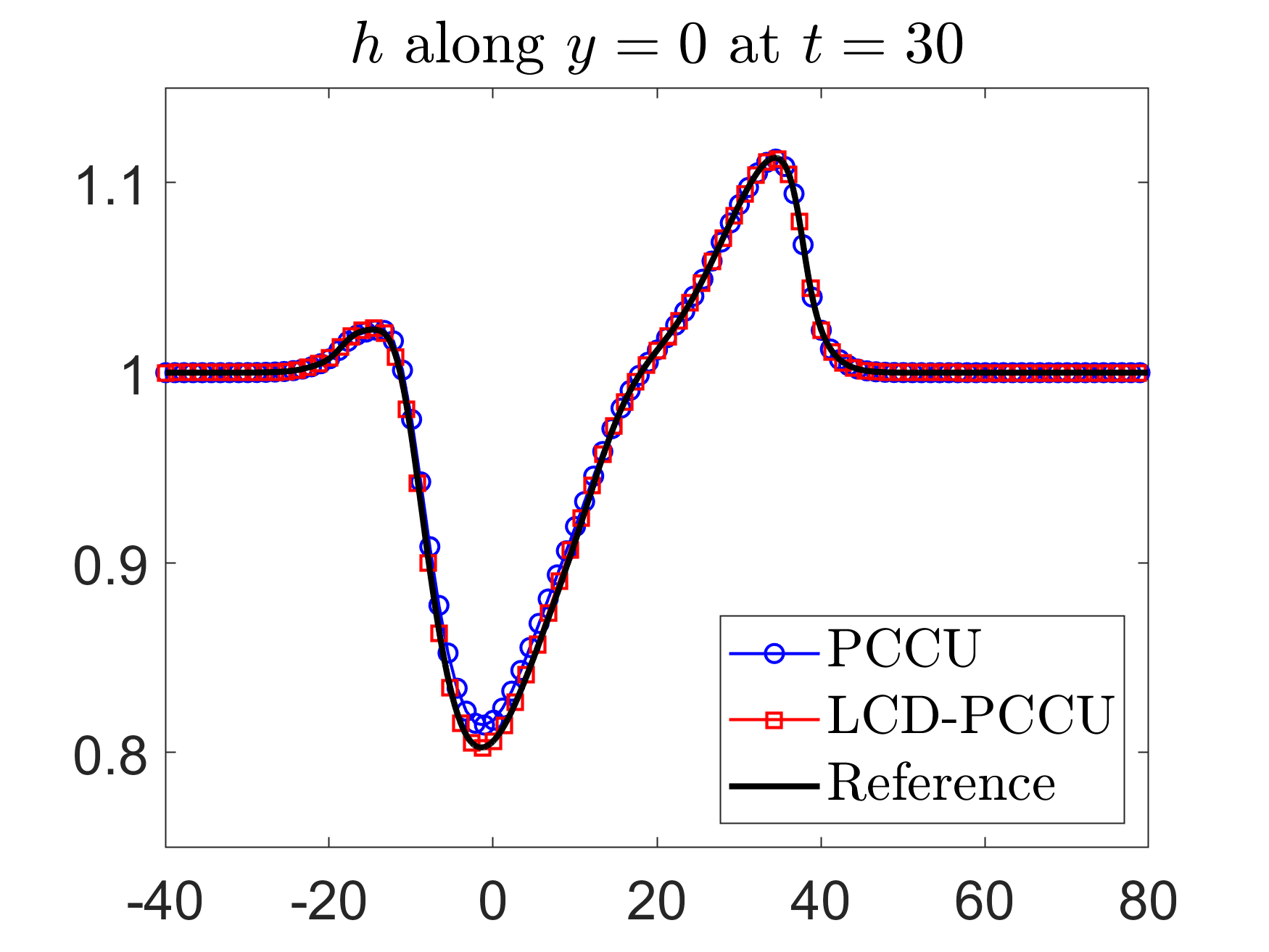}\hspace*{1.0cm}
            \includegraphics[trim=0.8cm 0.3cm 1.1cm 0.1cm, clip, width=5.cm]{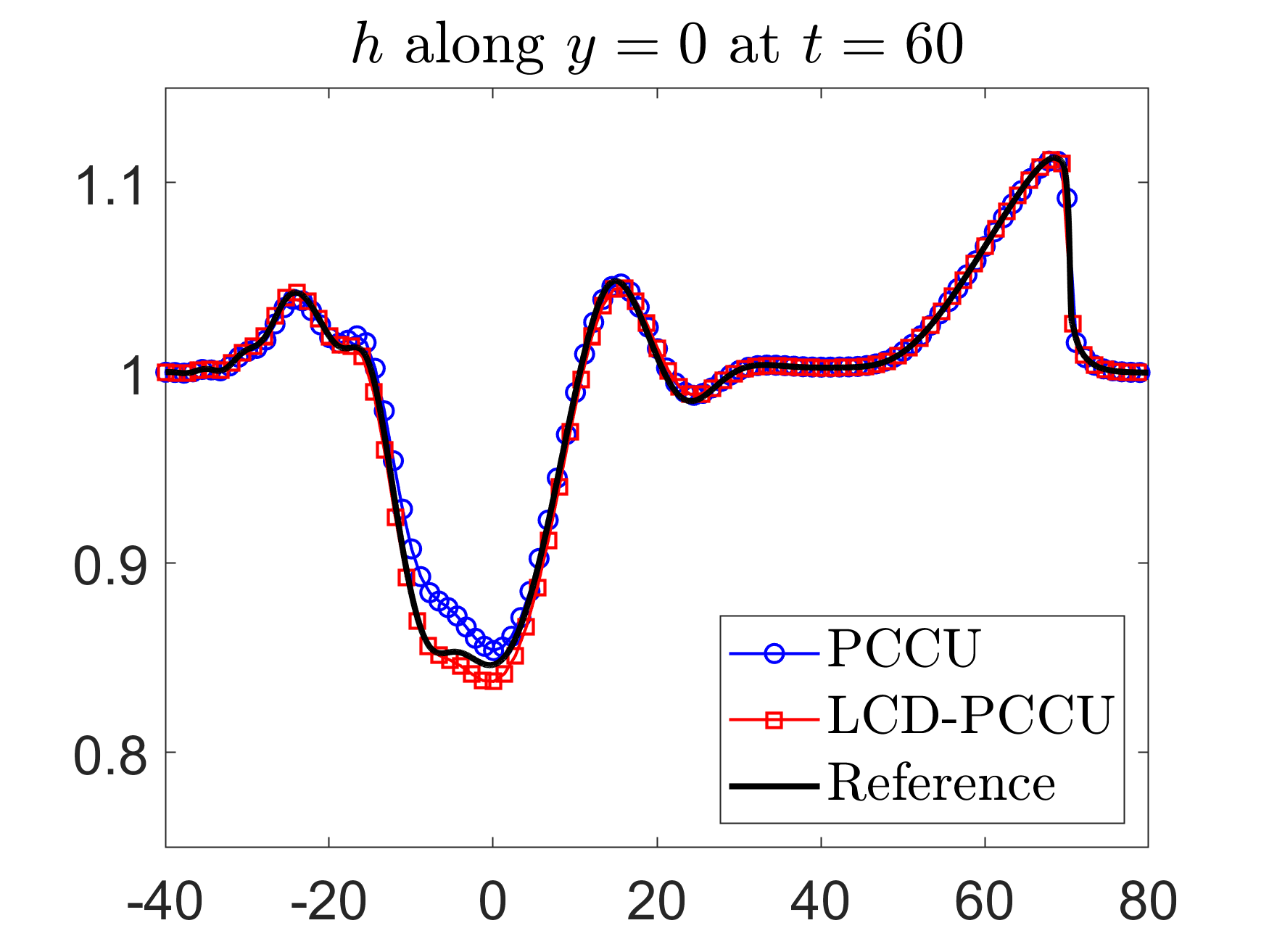}}
\vskip7pt
\centerline{\includegraphics[trim=0.8cm 0.3cm 1.1cm 0.1cm, clip, width=5.cm]{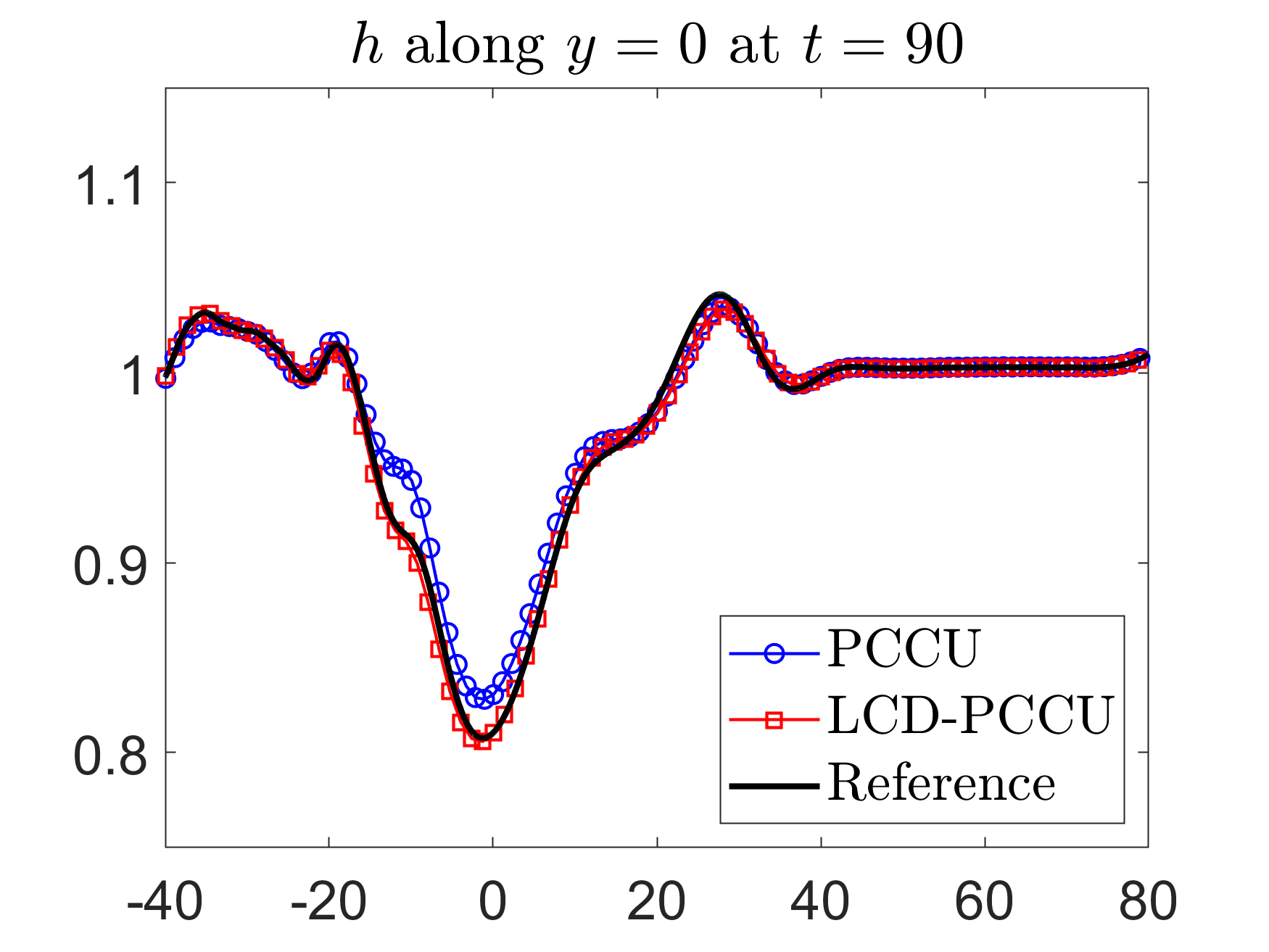}\hspace*{1.0cm}
            \includegraphics[trim=0.8cm 0.3cm 1.1cm 0.1cm, clip, width=5.cm]{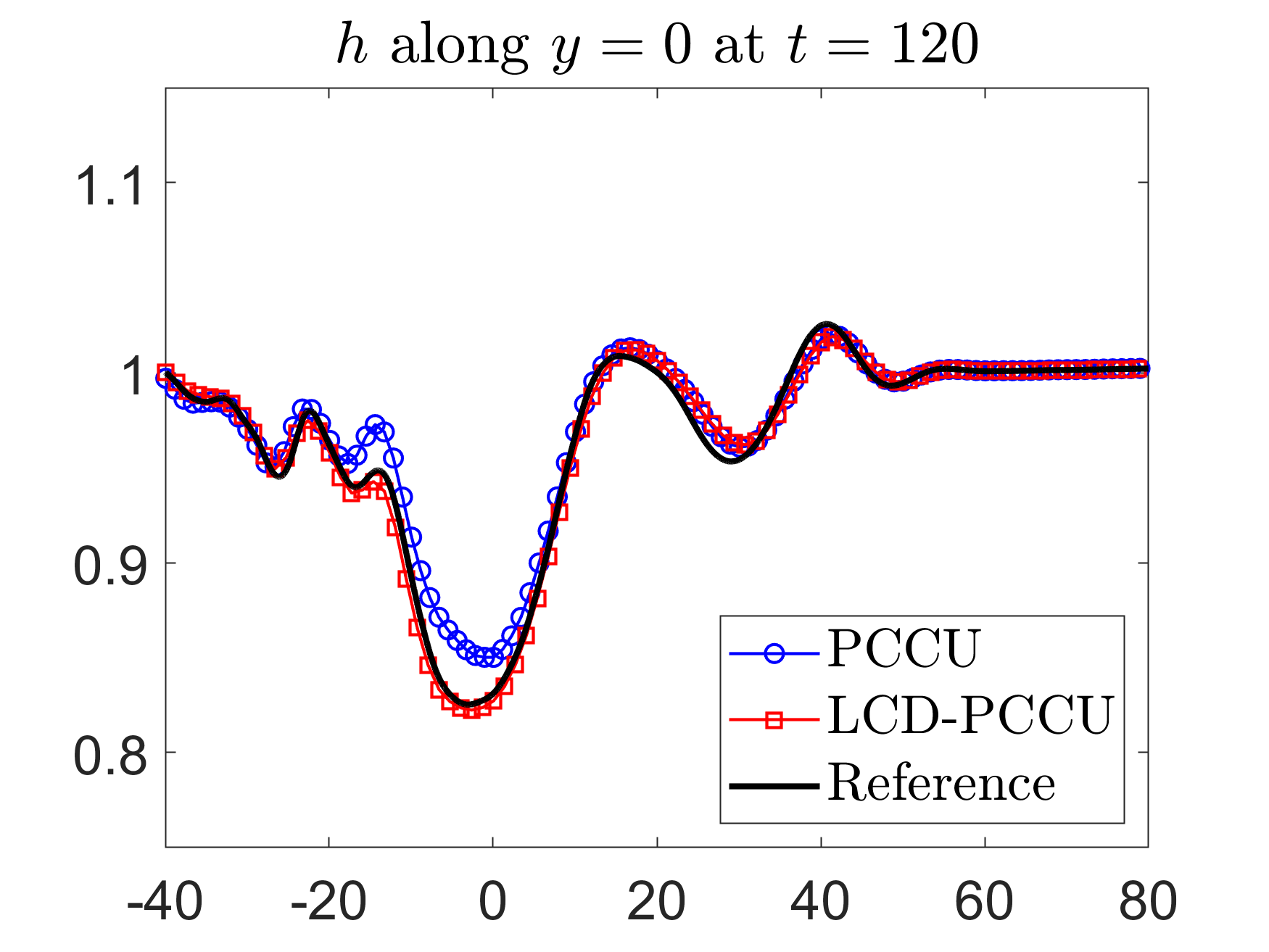}}
\caption{\sf Example 10: Slices of $h$ and $b$ along the line $y=0$ computed by the PCCU scheme using $\dx=\dy=1/9$ and LCD-PCCU scheme using $\dx=\dy=2/15$ at different times.}\label{fig1132}
\end{figure}

In this example, we also check the efficiency of the proposed LCD-PCCU scheme by performing a more thorough comparison between the studied
schemes. To this end, we measure the CPU time consumed during the above calculations by the LCD-PCCU scheme and refine the mesh used by the
PCCU scheme to the level that exactly the same CPU time is consumed to compute both of the numerical solutions: This way, we take into
account additional computational cost of the LCD-PCCU scheme. The corresponding grids are $\dx=\dy=2/15$ for the LCD-PCCU scheme, and
$\dx=\dy=1/9$ for the PCCU scheme. The obtained numerical results are plotted in Figure \ref{fig1132}, where one can clearly see that the
LCD-PCCU solution is still more accurate than the PCCU one.

\subsection*{Acknowledgments}
The work of Michael Herty was funded by the Deutsche Forschungsgemeinschaft (DFG, German Research Foundation)--SPP 2183:
Eigenschaftsgeregelte Umformprozesse with the Project(s) HE5386/19-2,19-3 Entwicklung eines flexiblen isothermen Reckschmiedeprozesses
f\"ur die eigenschaftsgeregelte Herstellung von Turbinenschaufeln aus Hochtemperaturwerkstoffen (424334423) and by the Deutsche
Forschungsgemeinschaft (DFG, German Research Foundation)--SPP 2410 Hyperbolic Balance Laws in Fluid Mechanics: Complexity, Scales,
Randomness (CoScaRa) within the Project(s) HE5386/26-1 (Numerische Verfahren f\"ur gekoppelte Mehrskalenprobleme,525842915) and
(Zuf\"allige kompressible Euler Gleichungen: Numerik und ihre Analysis, 525853336) HE5386/27-1. The work of A. Kurganov was supported in
part by NSFC grant 12171226 and by the fund of the Guangdong Provincial Key Laboratory of Computational Science and Material Design (No.
2019B030301001).

\appendix
\section{LCD Matrices for the $\gamma$-Based Compressible Multifluid Equations}\label{appa}
\subsection{One-Dimensional Case}
In order to compute the LCD matrices for the 1-D $\gamma$-based compressible multifluid system, we first compute the matrix
${\cal A}=\frac{\partial\mF}{\partial\mU}-B$:

\begin{equation}
{\cal A}(\mU)=\begin{pmatrix}0&1&0&0&0\\\dfrac{\gamma-3}{2}u^2&(3-\gamma)u&\gamma-1&(1-\gamma)p&1-\gamma\\[0.8ex]
\dfrac{uc^2}{1-\gamma}+\Big(\dfrac{\gamma}{2}-1\Big)u^3&\dfrac{c^2}{\gamma-1}+\Big(\dfrac{3}{2}-\gamma\Big)u^2&\gamma u&(1-\gamma)pu&
(1-\gamma)u\\0&0&0&u&0\\0&0&0&0&u\end{pmatrix},
\label{A1}
\end{equation}
where $c=\sqrt{\frac{\gamma(p+\pi_\infty)}{\rho}}$ is the speed of sound, and then introduce
$\widehat{\cal A}_\jph={\cal A}(\widehat\mU_\jph)$, which is given by \eref{A1} with $\gamma$, $u$, $p$, and $c$ replaced with
$$
\hat\gamma=\frac{\gamma_j+\gamma_{j+1}}{2},\quad\hat u=\frac{u_j+u_{j+1}}{2},\quad\hat p=\frac{p_j+p_{j+1}}{2},\quad{\rm and}\quad
\hat c=\sqrt{\frac{\hat\gamma(\hat p+\hat\pi_\infty)}{\hat\rho}},
$$
respectively. Here,
$$
\begin{aligned}
&\gamma_j=\frac{1}{\xbar\Gamma_j}+1,\quad(\pi_\infty)_j=\frac{\xbar\Pi_j}{\xbar\Gamma_j+1},\quad\hat
\rho=\frac{\xbar\rho_j+\xbar\rho_{j+1}}{2},\quad u_j=\frac{(\xbar{\rho u})_j}{\xbar\rho_j},\\ 
&p_j=(\gamma_j-1)\Big[\xbar E_j-\hf\,\xbar\rho_j u_j^2\Big]-\gamma_j(\pi_\infty)_j,\quad
\hat\pi_\infty=\frac{(\pi_\infty)_j+(\pi_\infty)_{j+1}}{2}.
 \end{aligned}
$$

Notice that all of the $\hat{(\cdot)}$ quantities have to have a subscript index, that is, $\hat{(\cdot)}=\hat{(\cdot)}_\jph$, but we omit
it for the sake of brevity for all of the quantities except for $\widehat{\cal A}_\jph$. We then compute the matrix $R_\jph$ composed of the
right eigenvectors of $\widehat{\cal A}_\jph$ and obtain
\begin{equation*}
\begin{aligned}
R_\jph=&\begin{pmatrix}1&1&0&0&1\\\hat u-\hat c&\hat u&0&0&\hat u+\hat c\\[0.3ex]
\dfrac{\hat c^2}{\hat\gamma-1}+\dfrac{1}{2}\hat u^2-\hat u\hat c&\dfrac{\hat u^2}{2}&\hat p&0&
\dfrac{\hat c^2}{\hat\gamma-1}+\dfrac{1}{2}\hat u^2+\hat u\hat c\\[0.5ex]
0&0&1&-1&0\\0&0&0&\hat p&0\end{pmatrix},\\[1.0ex]
R^{-1}_\jph=\dfrac{1}{2\hat c^2}&\begin{pmatrix}
\dfrac{(\hat\gamma-1)\hat u^2}{2}+\hat u\hat c&-\hat c-(\hat\gamma-1)\hat u&\hat\gamma-1&(1-\hat\gamma)\hat p&1-\hat\gamma\\[0.5ex]
2\hat c^2-(\hat\gamma-1)\hat u^2&2(\hat\gamma-1)\hat u&2(1-\hat\gamma)&2(\hat\gamma-1)\hat p&2(\hat\gamma-1)\\[0.3ex]
0&0&0&2\hat c^2&\dfrac{2\hat c^2}{\hat p}\\[1.5ex]
0&0&0&0&-\dfrac{2\hat c^2}{\hat p}\\
\dfrac{(\hat\gamma-1)\hat u^2}{2}-\hat u\hat c&\hat c-(\hat\gamma-1)\hat u&\hat\gamma-1&(1-\hat\gamma)\hat p&1-\hat\gamma\end{pmatrix}.
\end{aligned}
\end{equation*}

\subsection{Two-Dimensional Case}
We now compute the matrix $\widehat{\cal A}^{\,\mF}_{\jph,k}$ and the corresponding matrices $R_{\jph,k}$ and $R^{-1}_{\jph,k}$ for the 2-D
$\gamma$-based compressible multifluid system.

We first compute ${\cal A}^{\mF}=\frac{\partial\mF}{\partial\mU}-B$:
\begin{equation}
\resizebox{0.92\hsize}{!}{$
{\cal A}^{\mF}(\mU)=\begin{pmatrix}0&1&0&0&0&0\\
\dfrac{\gamma-3}{2}u^2+\dfrac{\gamma-1}{2}v^2&(3-\gamma)u&(1-\gamma)v&\gamma-1&(1-\gamma)p&1-\gamma\\[0.5ex]
-uv&v&u&0&0&0\\
\dfrac{uc^2}{1-\gamma}+\Big(\dfrac{\gamma}{2}-1\Big)u(u^2+v^2)&\dfrac{c^2}{\gamma-1}+\Big(\dfrac{3}{2}-\gamma\Big)u^2+\dfrac{1}{2}v^2&
(1-\gamma)uv&\gamma u&(1-\gamma)pu&(1-\gamma)u\\0&0&0&0&u&0\\0&0&0&0&0&u
\end{pmatrix},$}
\label{A21}
\end{equation}
and then introduce the matrix $\widehat{\cal A}^{\,\mF}_{\jph,k}={\cal A}(\widehat\mU_{\jph,k})$, which is given by \eref{A21} with
$\gamma$, $u$, $v$, $p$, and $c$ replaced with 
\begin{equation*}
\hat\gamma=\frac{\gamma_{j,k}+\gamma_{j+1,k}}{2},\quad\hat u=\frac{u_{j,k}+u_{j+1,k}}{2},\quad\hat v=\frac{v_{j,k}+v_{j+1,k}}{2},\quad
\hat p=\frac{p_{j,k}+p_{j+1,k}}{2},\quad\hat c=\sqrt{\frac{\hat\gamma(\hat p+\hat\pi_\infty)}{\hat\rho}},
\end{equation*}
respectively. Here, 
\begin{equation*}
\begin{aligned}
&\gamma_{j,k}=\frac{1}{\xbar\Gamma_{j,k}}+1,\quad(\pi_\infty)_{j,k}=\frac{\xbar\Pi_{j,k}}{\xbar\Gamma_{j,k}+1},\quad
\hat\rho=\frac{\xbar\rho_{j,k}+\xbar\rho_{j+1,k}}{2},\quad u_{j,k}=\frac{(\xbar{\rho u})_{j,k}}{\xbar\rho_{j,k}},\quad
v_{j,k}=\frac{(\xbar{\rho u})_{j,k}}{\xbar\rho_{j,k}},\\
&p_{j,k}=(\gamma_{j,k}-1)\Big[\xbar E_{j,k}-\hf\,\xbar\rho_{j,k}u_{j,k}^2-\hf\,\xbar\rho_{j,k}v_{j,k}^2\Big]-\gamma_{j,k}(\pi_\infty)_{j,k},
\quad\hat\pi_\infty=\frac{(\pi_\infty)_{j,k}+(\pi_\infty)_{j+1,k}}{2}.
\end{aligned}
\end{equation*}

Notice that all of the $\hat{(\cdot)}$ quantities have to have a subscript index, that is, $\hat{(\cdot)}=\hat{(\cdot)}_{\jph,k}$, but we
omit it for the sake of brevity for all of the quantities except for $\widehat{\cal A}^{\,\mF}_{\jph,k}$. We then compute the matrix
$R_{\jph,k}$ composed of the right eigenvectors of $\widehat{\cal A}^{\,\mF}_{\jph,k}$ and obtain

\begin{equation*}
\begin{aligned}
R_{\jph,k}=&\begin{pmatrix}1&1&0&0&0&1\\\hat u-\hat c&\hat u&0&0&0&\hat u+\hat c\\\hat v&\hat v&1&0&0&\hat v\\
\dfrac{\hat c^2}{\hat\gamma-1}+\dfrac{1}{2}\hat u^2+\dfrac{1}{2}\hat v^2-\hat u\hat c&\dfrac{\hat u^2+\hat v^2}{2}&\hat v&\hat p&0&
\dfrac{\hat c^2}{\hat \gamma-1}+\dfrac{1}{2}\hat u^2+\dfrac{1}{2}\hat v^2+\hat u\hat c\\0&0&0&1&1&0\\0&0&0&0&-\hat p&0\end{pmatrix}\\[1.0ex]
R^{-1}_{\jph,k}=\dfrac{1}{2\hat c^2}&\begin{pmatrix}
\dfrac{(\hat\gamma-1)(\hat u^2+\hat v^2)}{2}+\hat u\hat c&-\hat c-(\hat\gamma-1)\hat u&(1-\hat\gamma)\hat v&\hat\gamma-1&
(1-\hat\gamma)\hat p&1-\hat\gamma\\[0.5ex]
2\hat c^2-(\hat\gamma-1)(\hat u^2+\hat v^2)&2(\hat\gamma-1)\hat u&2(\hat\gamma-1)\hat v&2(1-\hat\gamma)&2(\hat\gamma-1)\hat p&
2(\hat\gamma-1)\\-2\hat c^2\hat v&0&2\hat c^2&0&0&0\\[0.3ex]
0&0&0&0&2\hat c^2&\dfrac{2\hat c^2}{\hat p}\\\dfrac{(\hat\gamma-1)(\hat u^2+\hat v^2)}{2}-\hat u\hat c&\hat c-(\hat\gamma-1) hat u&
(1-\hat\gamma)\hat v&\hat\gamma-1&(1-\hat\gamma)\hat p&1-\hat\gamma\end{pmatrix}.
\end{aligned}
\end{equation*}

\section{LCD Matrices for the TRSW Equations}\label{appc}
\subsection{One-dimensional Case}
In order to compute the LCD matrices for the 1-D TRSW equations, we first compute the matrix ${\cal A}=\frac{\partial\mG}{\partial\mU}-C$:
\begin{equation}
{\cal A}(\mU)=\begin{pmatrix}0&0&1&0\\-uv&v&u&0\\[0.3ex]
\dfrac{bh}{2}-v^2&0&2v&\dfrac{h}{2}\\[0.3ex]
-bv&0&b&v
\end{pmatrix},
\label{C1}
\end{equation}
and then introduce the matrices $\widehat{\cal A}_\kph={\cal A}(\widehat\mU_\kph)$, which are given by \eref{C1} with $h$, $u$, $v$, and $b$
replaced with 
\begin{equation*}
\hat h=\frac{\xbar h_k+\xbar h_{k+1}}{2},\quad\hat u=\frac{u_k+u_{k+1}}{2},\quad\hat v=\frac{v_k+v_{k+1}}{2},\quad
\hat b=\frac{b_k+b_{k+1}}{2},
\end{equation*}
respectively. Here, $u_k=\nicefrac{\xbar q_k}{\,\xbar h_k}$, $v_k=\nicefrac{\xbar p_k}{\,\xbar h_k}$, and
$b_k=\nicefrac{\xbar{(hb)}_k}{\,\xbar h_k}$.

Notice that all of the $\hat{(\cdot)}$ quantities have to have a subscript index, that is, $\hat{(\cdot)}=\hat{(\cdot)}_\kph$, but we omit
it for the sake of brevity for all of the quantities except for $\widehat{\cal A}_\kph$. We then compute the matrix $R_\kph$ composed of the
right eigenvectors of $\widehat{\cal A}_\kph$ and obtain

\begin{equation*}
R_\kph=\frac{1}{\hat b}\begin{pmatrix}1&-1&0&1\\\hat u&0&\hat b&\hat u\\\hat v-\sqrt{\hat b\hat h}&-\hat v&0&\hat v+\sqrt{\hat b\hat h}\\
\hat b&\hat b&0&\hat b\end{pmatrix},\quad 
R^{-1}_\kph=\frac{1}{4}\begin{pmatrix}\hat b+2\hat v\sqrt{\dfrac{\hat b}{\hat h}}&0&-2\sqrt{\dfrac{\hat b}{\hat h}}&1\\[1.8ex]
-2\hat b&0&0&2\\-2\hat u&4&0&-2\dfrac{\hat u}{\hat b}\\[1.0ex]
\hat b-2\hat v\sqrt{\dfrac{\hat b}{\hat h}}&0&2\sqrt{\dfrac{\hat b}{\hat h}}&1\end{pmatrix}.
\end{equation*}

\subsection{Two-Dimensional Case}
We now compute the matrix $\widehat{\cal A}^{\,\mF}_{\jph,k}$ and the corresponding matrices $R_{\jph,k}$ and $R^{-1}_{\jph,k}$ for the 2-D
TRSW equations.

We first compute ${\cal A}^{\,\mF}=\frac{\partial\mF}{\partial\mU}-B$:
\begin{equation}
{\cal A}^{\,\mF}(\mU)=\begin{pmatrix}0&1&0&0\\\dfrac{bh}{2}-u^2&2u&0&\dfrac{h}{2}\\[0.8ex]
-uv&v&u&0\\-bu&b&0&u\end{pmatrix},
\label{D2}
\end{equation}
and then introduce the matrices $\widehat{\cal A}^{\,\mF}_{\jph,k}={\cal A}(\widehat\mU_{\jph,k})$, which are given by \eref{D2} with $h$,
$u$, $v$, and $b$ replaced with 
\begin{equation*}
\hat h=\frac{\xbar h_{j,k}+\xbar h_{j+1,k}}{2},\quad\hat u=\frac{u_{j,k}+u_{j+1,k}}{2},\quad\hat v=\frac{v_{j,k}+v_{j+1,k}}{2},\quad
\hat b=\frac{b_{j,k}+b_{j+1,k}}{2},
\end{equation*}
respectively. Here $u_k=\nicefrac{\xbar q_{j,k}}{\,\xbar h_{j,k}}$, $v_k=\nicefrac{\xbar p_{j,k}}{\,\xbar h_{j,k}}$, and
$b_k=\nicefrac{\xbar{(hb)}_{j,k}}{\,\xbar h_{j,k}}$.

Notice that all of the $\hat{(\cdot)}$ quantities have to have a subscript index, that is, $\hat{(\cdot)}=\hat{(\cdot)}_{\jph,k}$, but we
omit it for the sake of brevity for all of the quantities except for $\widehat{\cal A}^{\,\mF}_{\jph,k}$. We then compute the matrix
$R_\jph$ composed of the right eigenvectors of $\widehat{\cal A}^{\,\mF}_{\jph,k}$ and obtain
\begin{equation*}
R_{\jph,k}=\frac{1}{\hat b}\begin{pmatrix}1&-1&0&1\\\hat u-\sqrt{\hat b\hat h}&-\hat u&0&\hat u+\sqrt{\hat b\hat h}\\
\hat v&0&\hat b&\hat v\\\hat b&\hat b&0&\hat b\end{pmatrix},\quad
R^{-1}_{\jph,k}=\frac{1}{4}\begin{pmatrix}\hat b+2\hat u\sqrt{\dfrac{\hat b}{\hat h}}&0&-2\sqrt{\dfrac{\hat b}{\hat h}}&1\\[1.8ex]
-2\hat b&0&0&2\\-2\hat v&0&4&-2\dfrac{\hat v}{\hat b}\\[1.0ex]
\hat b-2\hat u\sqrt{\dfrac{\hat b}{\hat h}}&2\sqrt{\dfrac{\hat b}{\hat h}}&0&1\end{pmatrix}.
\end{equation*}

\bibliography{reference}
\bibliographystyle{siamnodash}
\end{document}